\documentclass[12pt, a4paper, reqno]{amsart}
\usepackage{amsmath}
\usepackage{amsfonts}
\usepackage{amssymb}
\usepackage{color}
\usepackage{comment,graphicx}
\usepackage[T1]{fontenc}
\usepackage{url}
\usepackage{ae}
\usepackage{mathtools}
\usepackage{esint}
\usepackage{latexsym}
\usepackage{amscd}

\def\phi{\varphi}
\def\rho{\varrho}
\def\epsilon{\varepsilon}
\numberwithin{equation}{section}
\theoremstyle{plain}
\newtheorem{thm}[equation]{Theoremm}
\newtheorem{lem}[equation]{Lemma}
\newtheorem{prop}[equation]{Proposition}
\newtheorem{cor}[equation]{Corollary}
\theoremstyle{definition}
\newtheorem{defn}[equation]{Definition}

\theoremstyle{remark}
\newtheorem{rem}[equation]{Remark}
\newtheorem{ex}[equation]{Example}

\renewcommand{\leq}{\leqslant}
\renewcommand{\geq}{\geqslant}

\begin{document}
\title[Weigthed Triebel-Lizorkin spaces]{ Triebel-Lizorkin spaces with
general weights}
\author[D. Drihem]{Douadi Drihem}
\address{Douadi Drihem\\
M'sila University\\
Department of Mathematics\\
Laboratory of Functional Analysis and Geometry of Spaces\\
M'sila 28000, Algeria.}
\email{douadidr@yahoo.fr, douadi.drihem@univ-msila.dz}
\thanks{ }
\date{\today }

\begin{abstract}
In this paper, the author introduces Triebel-Lizorkin spaces with general
smoothness. We present the $\varphi $-transform characterization of these
spaces in the sense of Frazier and Jawerth and we prove their Sobolev
embeddings. Also, we establish the smooth atomic and molecular decomposition
of these function spaces. To do these we need a generalization of some
maximal inequality to the case of general weights.
\end{abstract}

\maketitle

\section{Introduction}

This paper is a continuation of \cite{D20}, where the author introduced
Besov spaces with general smoothness and presented some their 
properties, such as the $\varphi $-transform characterization in the sense
of Frazier and Jawerth, the smooth atomic, molecular and wavelet
decomposition, and the characterization of these function spaces in terms of
the difference relations.

The spaces of generalized smoothness are have been introduced by several
authors. We refer, for instance, to Bownik \cite{M07}, Cobos and Fernandez 
\cite{CF88}, Goldman \cite{Go79} and \cite{Go83}, and Kalyabin \cite{Ka83};
see also\ Besov \cite{B03} and \cite{B05}, and Kalyabin and Lizorkin \cite%
{Kl87}. More general Besov spaces with variable smoothness were explicitly
studied by Ansorena and Blasco \cite{AB95} and \cite{AB96}, including
characterizations by differences and atomic decomposition. The wavelet
decomposition (with respect to a compactly supported wavelet basis of
Daubechies type) of nonhomogeneous Besov spaces of generalized smoothness
was achieved by Almeida \cite{Al}. The theory of these spaces had a
remarkable development in part due to its usefulness in applications. For
instance, they appear in the study of trace spaces on fractals, see Edmunds
and Triebel \cite{ET96} and \cite{ET99}, were they introduced the spaces $%
B_{p,q}^{s,\Psi }$, where $\Psi $ is a so-called admissible function,
typically of log-type near $0$. For a complete treatment of these spaces we
refer the reader the work of Moura \cite{Mo01}. Further results on Besov
spaces of variable smoothness are given in \cite{AA16}.

More general function spaces of generalized smoothness can be found in Caetano and Leopold \cite{CL},
Farkas and Leopold \cite{FL06}, and reference therein.

Recently, Dominguez and Tikhonov \cite{DT}  gave a  treatment of function spaces with logarithmic smoothness (Besov, Sobolev, Triebel-Lizorkin), including various new characterizations for Besov norms in terms of different,  sharp estimates for Besov norms of derivatives and potential operators (Riesz and Bessel potentials) in terms of norms of functions themselves and  sharp embeddings between the Besov spaces defined by differences and by
Fourier-analytical decompositions as well as between Besov and Sobolev/Triebel-Lizorkin spaces.

Tyulenev introduced in \cite{Ty15}, \cite{Ty-N-L}  and \cite{Ty-151} a new
family of Besov spaces of variable smoothness which cover many classes of
Besov spaces, where the norm on these spaces was defined\ with the help of
classical differences. Based on this weighted class and \cite{D20} we
introduce Triebel-Lizorkin spaces of variable smoothness, is defined as
follows. Let $\mathcal{S}(\mathbb{R}^{n})$ be the set of all Schwartz
functions $\varphi $ on $\mathbb{R}^{n}$, i.e., $\varphi $ is infinitely
differentiable and 
\begin{equation*}
\big\|\varphi |\mathcal{S}_{M}\big\|=\sup_{\beta \in \mathbb{N}%
	_{0}^{n},|\beta |\leq M}\sup_{x\in \mathbb{R}^{n}}|\partial ^{\beta }\varphi
(x)|(1+|x|)^{n+M+|\beta |}<\infty
\end{equation*}%
for all $M\in \mathbb{N}$. Select a Schwartz function $\varphi $ such that%
\begin{equation*}
\text{supp}(\mathcal{F}(\varphi ))\subset \big\{\xi :1/2\leq |\xi |\leq 2%
\big\}
\end{equation*}%
and 
\begin{equation*}
|\mathcal{F}(\varphi )(\xi )|\geq c\text{\quad if\quad }3/5\leq |\xi |\leq
5/3
\end{equation*}%
where $c>0$ and we put $\varphi _{k}=2^{kn}\varphi (2^{k}\cdot )$, $k\in 
\mathbb{Z}$.\textrm{\ }Here $\mathcal{F}(\varphi )$ denotes the Fourier
transform of $\varphi $, defined by%
\begin{equation*}
\mathcal{F}(\varphi )(\xi )=(2\pi )^{-n/2}\int_{\mathbb{R}^{n}}e^{-ix\cdot
	\xi }\varphi (x)dx,\quad \xi \in \mathbb{R}^{n}.
\end{equation*}%
Let%
\begin{equation*}
\mathcal{S}_{\infty }(\mathbb{R}^{n})=\Big\{\varphi \in \mathcal{S}(\mathbb{R%
}^{n}):\int_{\mathbb{R}^{n}}x^{\beta }\varphi (x)dx=0\text{ for all
	multi-indices }\beta \in \mathbb{N}_{0}^{n}\Big\}.
\end{equation*}%
Following Triebel \cite{T1}, we consider $\mathcal{S}_{\infty }(\mathbb{R}%
^{n})$ as a subspace of $\mathcal{S}(\mathbb{R}^{n})$, including the
topology. Thus, $\mathcal{S}_{\infty }(\mathbb{R}^{n})$ is a complete metric
space. Equivalently, $\mathcal{S}_{\infty }(\mathbb{R}^{n})$ can be defined
as a collection of all $\varphi \in \mathcal{S}(\mathbb{R}^{n})$ such that
semi-norms%
\begin{equation*}
\big\|\varphi \big\|_{M}=\sup_{\left\vert \beta \right\vert \leq M}\sup_{\xi
	\in \mathbb{R}^{n}}|\partial ^{\beta }\varphi (\xi )|\big(\left\vert \xi
\right\vert ^{M}+\left\vert \xi \right\vert ^{-M}\big)<\infty
\end{equation*}%
for all $M\in \mathbb{N}_{0}$, see \cite[Section 3]{BoHo06}. Let $\mathcal{S}%
_{\infty }^{\prime }(\mathbb{R}^{n})$ be the topological dual of $\mathcal{S}%
_{\infty }(\mathbb{R}^{n})$, namely, the set of all continuous linear
functionals on $\mathcal{S}_{\infty }(\mathbb{R}^{n})$. Let $0<p<\infty $
and $0<q\leq \infty $. Let $\{t_{k}\}$ be a $p$-admissible sequence i.e. $%
t_{k}\in L_{p}^{\mathrm{loc}}(\mathbb{R}^{n})$, $k\in \mathbb{Z}$. The
Triebel-Lizorkin space $\dot{F}_{p,q}(\mathbb{R}^{n},\{t_{k}\})$\ is the
collection of all $f\in \mathcal{S}_{\infty }^{\prime }(\mathbb{R}^{n})$\
such that 
\begin{equation*}
\big\|f|\dot{F}_{p,q}(\mathbb{R}^{n},\{t_{k}\})\big\|=\Big\|\Big(%
\sum\limits_{k=-\infty }^{\infty }t_{k}^{q}|\varphi _{k}\ast f|^{q}\Big)%
^{1/q}|L_{p}(\mathbb{R}^{n})\Big\|<\infty
\end{equation*}%
with the usual modifications if $q=\infty $.

We have organized the article in three sections. First we give some
preliminaries and recall some basic facts on the Muckenhoupt classes and the
weighted class of Tyulenev. Also we give some key technical lemmas needed in
the proofs of the main statements. Especially, we present new version of
weighted vector-valued maximal inequality of Fefferman and Stein. In Section
3 several basic properties such as the $\varphi $-transform characterization
are obtained. We extend the well-known Sobolev embeddings to these function
spaces and we give the atomic and molecular decomposition of these function
spaces.

In Section 4 we study the inhomogeneous spaces $F_{p,q}(\mathbb{R}%
^{n},\{t_{k}\})$, and we outline analogous results for these spaces. Other
properties of these function spaces are given in \cite{D20.2}.

\section{Maximal\ inequalities}

Our arguments of this paper essentially rely on the weighted boundedness
of Hardy-Littlewood maximal function. In this paper we will assume that the
weight sequence $\{t_{k}\}$ used to define the space $\dot{F}_{p,q}(\mathbb{R%
}^{n},\{t_{k}\})$ lies in the new weighted class $\dot{X}_{\alpha ,\sigma
	,p} $ (see Definition \ref{Tyulenev-class}). Therefore we need a new\
version of Hardy-Littlewood maximal inequality. Throughout this paper, we
make some notation and conventions.

\subsection{Notation and conventions}

We denote by $\mathbb{R}^{n}$ the $n$-dimensional real Euclidean space, $%
\mathbb{N}$ the collection of all natural numbers and $\mathbb{N}_{0}=%
\mathbb{N}\cup \{0\}$. The letter $\mathbb{Z}$ stands for the set of all
integer numbers.\ The expression $f\lesssim g$ means that $f\leq c\,g$ for
some independent constant $c$ (and non-negative functions $f$ and $g$), and $%
f\approx g$ means $f\lesssim g\lesssim f$.\vskip5pt

By supp($f$) we denote the support of the function $f$, i.e., the closure of
its non-zero set. If $E\subset {\mathbb{R}^{n}}$ is a measurable set, then $%
|E|$ stands for the (Lebesgue) measure of $E$ and $\chi _{E}$ denotes its
characteristic function. By $c$ we denote generic positive constants, which
may have different values at different occurrences. \vskip5pt

A weight is a nonnegative locally integrable function on $\mathbb{R}^{n}$
that takes values in $(0,\infty )$ almost everywhere. For measurable set $%
E\subset \mathbb{R}^{n}$ and a weight $\gamma $, $\gamma (E)$ denotes 
\begin{equation*}
\int_{E}\gamma (x)dx.
\end{equation*}%
Given a measurable set $E\subset \mathbb{R}^{n}$ and $0<p\leq \infty $, we
denote by $L_{p}(E)$ the space of all functions $f:E\rightarrow \mathbb{C}$
equipped with the quasi-norm 
\begin{equation*}
\big\|f|L_{p}(E)\big\|=\Big(\int_{E}\left\vert f(x)\right\vert ^{p}dx\Big)%
^{1/p}<\infty ,
\end{equation*}%
with $0<p<\infty $ and%
\begin{equation*}
\big\|f|L_{\infty }(E)\big\|=\underset{x\in E}{\text{ess-sup}}\left\vert
f(x)\right\vert <\infty .
\end{equation*}%
For a function $f$ in $L_{1}^{\mathrm{loc}}(\mathbb{R}^{n})$, we set 
\begin{equation*}
M_{A}(f)=\frac{1}{|A|}\int_{A}\left\vert f(x)\right\vert dx
\end{equation*}%
for any $A\subset \mathbb{R}^{n}$. Furthermore, we put%
\begin{equation*}
M_{A,p}(f)=\Big(\frac{1}{|A|}\int_{A}\left\vert f(x)\right\vert ^{p}dx\Big)%
^{1/p},
\end{equation*}%
with $0<p<\infty $. Further, given a measurable set $E\subset \mathbb{R}^{n}$
and a weight $\gamma $, we denote the space of all functions $f:\mathbb{R}%
^{n}\rightarrow \mathbb{C}$ with finite quasi-norm 
\begin{equation*}
\big\|f|L_{p}(\mathbb{R}^{n},\gamma )\big\|=\big\|f\gamma |L_{p}(\mathbb{R}%
^{n})\big\|
\end{equation*}%
by $L_{p}(\mathbb{R}^{n},\gamma )$.

If $1\leq p\leq \infty $ and $1/p+1/p^{\prime }=1$, then $p^{\prime }$ is
called the conjugate exponent of $p$.

Let $0<p,q\leq \infty $. The space $L_{p}(\ell _{q})$ is defined to be the
set of all sequences $\{f_{k}\}$ of functions such that%
\begin{equation*}
\big\|\{f_{k}\}|L_{p}(\ell _{q})\big\|=\Big\|\Big(\sum_{k=-\infty }^{\infty
}|f_{k}|^{q}\Big)^{1/q}|L_{p}(\mathbb{R}^{n})\Big\|<\infty
\end{equation*}%
with the usual modifications if $q=\infty $ and if $\{t_{k}\}$ is a sequence
of functions then%
\begin{equation*}
\big\|\{f_{k}\}|L_{p}(\ell _{q},\{t_{k}\})\big\|=\big\|\{t_{k}f_{k}\}|L_{p}(%
\ell _{q})\big\|.
\end{equation*}
In what follows, $Q$ will denote a cube in the
space $\mathbb{R}^{n}$\ with sides parallel to the coordinate axes and $l(Q)$%
\ will denote the side length of the cube $Q$. For all cubes $Q$ and $r>0$,
let $rQ$ be the cube concentric with $Q$ having side length $rl(Q)$. For $%
v\in \mathbb{Z}$ and $m\in \mathbb{Z}^{n}$, denote by $Q_{v,m}$ the dyadic
cube,%
\begin{equation*}
Q_{v,m}=2^{-v}([0,1)^{n}+m).
\end{equation*}%
For the collection of all such cubes we use $\mathcal{Q}=\{Q_{v,m}:v\in 
\mathbb{Z},m\in \mathbb{Z}^{n}\}$. For each cube $Q$, we denote by $x_{v,m}$
the lower left-corner $2^{-v}m$ of $Q=Q_{v,m}$.

\subsection{Muckenhoupt weights}

The purpose of this subsection is to review some known properties of\
Muckenhoupt class.

\begin{defn}
	Let $1<p<\infty $. We say that a weight $\gamma $ belongs to the Muckenhoupt
	class $A_{p}(\mathbb{R}^{n})$ if there exists a constant $C>0$ such that for
	every cube $Q$ the following inequality holds 
	\begin{equation}
	M_{Q}(\gamma )M_{Q,p^{\prime }/p}(\gamma ^{-1})\leq C.  \label{Ap-constant}
	\end{equation}
\end{defn}

The smallest constant $C$ for which $\mathrm{\eqref{Ap-constant}}$ holds is
denoted by $A_{p}(\gamma )$. As an example, we can take%
\begin{equation*}
\gamma (x)=|x|^{\alpha },\quad \alpha \in \mathbb{R}.
\end{equation*}%
Then $\gamma \in A_{p}(\mathbb{R}^{n})$, $1<p<\infty $, if and only if $%
-n<\alpha <n(p-1)$.

For $p=1$ we rewrite the above definition in the following way.

\begin{defn}
	We say that a weight $\gamma $ belongs to the Muckenhoupt class $A_{1}(%
	\mathbb{R}^{n})$ if there exists a constant $C>0$ such that for every cube $%
	Q $\ and for a.e.\ $y\in Q$ the following inequality holds 
	\begin{equation}
	M_{Q}(\gamma )\leq C\gamma (y).  \label{A1-constant}
	\end{equation}
\end{defn}

The smallest constant $C$ for which $\mathrm{\eqref{A1-constant}}$ holds is
denoted by $A_{1}(\gamma )$. The above classes have been first studied by
Muckenhoupt\ \cite{Mu72} and use to characterize the boundedness of the
Hardy-Littlewood maximal function on $L_{p}(\gamma )$, see the monographs 
\cite{GR85} and \cite{L. Graf14}\ for a complete account on the theory of
Muckenhoupt weights.

We recall a few basic properties of the class of $A_{p}(\mathbb{R}^{n})$
weights, see \cite[Chapter 7]{Du01}, \cite[Chapter 7]{L. Graf14} and \cite[%
Chapter 5]{St93}.

\begin{lem}
	\label{Ap-Property}Let $1\leq p<\infty $.\newline
	$\mathrm{(i)}$ If $\gamma \in A_{p}(\mathbb{R}^{n})$, then for any $1\leq
	p<q $, $\gamma \in A_{q}(\mathbb{R}^{n})$.\newline
	$\mathrm{(ii)}$ Let $1<p<\infty $. $\gamma \in A_{p}(\mathbb{R}^{n})$ if and
	only if $\gamma ^{1-p^{\prime }}\in A_{p^{\prime }}(\mathbb{R}^{n})$.\newline
	$\mathrm{(iii)}$ Let $\gamma \in A_{p}(\mathbb{R}^{n})$. There is $C>0$ such
	that for any cube $Q$ and a measurable subset $E\subset Q$%
	\begin{equation*}
	\Big(\frac{|E|}{|Q|}\Big)^{p-1}M_{Q}(\gamma )\leq CM_{E}(\gamma ).
	\end{equation*}%
	$\mathrm{(iv)}$ Suppose that $\gamma \in A_{p}(\mathbb{R}^{n})$ for some $%
	1<p<\infty $. Then there exists a $1<p_{1}<p<\infty $ such that $\gamma \in
	A_{p_{1}}(\mathbb{R}^{n})$.\newline
	$\mathrm{(v)}$ If $\gamma \in A_{p}(\mathbb{R}^{n})$, then for any $%
	0<\varepsilon \leq 1$, $\gamma ^{\varepsilon }\in A_{p}(\mathbb{R}^{n})$.\newline
	$\mathrm{(vi)}$ Let $1\leq p<\infty $ and $\gamma \in A_{p}(\mathbb{R}^{n})$.
	Then there exist $\delta \in (0,1)$ and $C>0$ depending only on $n$, $p$,
	and $A_{p}(\gamma )$ such that for any cube $Q$ and any measurable subset $S$
	of $Q$ we have%
	\begin{equation*}
	\frac{M_{S}(\gamma )}{M_{Q}(\gamma )}\leq C\Big(\frac{|S|}{|Q|}\Big)^{\delta-1}.
	\end{equation*}
\end{lem}

The following theorem gives a useful property of $A_{p}(\mathbb{R}^{n})$
weights (reverse H\"{o}lder inequality), see \cite[Chapter 7]{L. Graf14} or 
\cite[Chapter 1]{LuDiYa07}.

\begin{thm}
	\label{reverse Holder inequality}Let $1\leq p<\infty \ $and $\gamma \in
	A_{p}(\mathbb{R}^{n})$. Then there exist a constants $C>0$ and $\varepsilon
	>0$\ depending only on $p$ and the $A_{p}(\mathbb{R}^{n})$ constant of $%
	\gamma $, such that for every cube $Q$, 
	\begin{equation*}
	M_{Q,1+\varepsilon }(\gamma )\leq CM_{Q}(\gamma ).
	\end{equation*}
\end{thm}

\subsection{The weight class $\dot{X}_{\protect\alpha ,\protect\sigma ,p}$}

Let $0<p\leq \infty $. A weight sequence $\{t_{k}\}$ is called $p$%
-admissible if $t_{k}\in L_{p}^{\mathrm{loc}}(\mathbb{R}^{n})$ for all $k\in 
\mathbb{Z}$. We mention here that 
\begin{equation*}
\int_{E}t_{k}^{p}(x)dx<c(k)
\end{equation*}%
for any\ compact set $E\subset \mathbb{R}^{n}$. For a $p$-admissible weight
sequence $\{t_{k}\}$\ we set%
\begin{equation*}
t_{k,m}=\big\|t_{k}|L_{p}(Q_{k,m})\big\|,\quad k\in \mathbb{Z},m\in \mathbb{Z%
}^{n}.
\end{equation*}

Tyulenev\ \cite{Ty14} introduced the following new weighted class\ and used
it to study Besov spaces of variable smoothness.

\begin{defn}
	\label{Tyulenev-class}Let $\alpha _{1}$, $\alpha _{2}\in \mathbb{R}$, $%
	p,\sigma _{1}$, $\sigma _{2}$ $\in (0,+\infty )$, $\alpha =(\alpha
	_{1},\alpha _{2})$ and let $\sigma =(\sigma _{1},\sigma _{2})$. We let $\dot{%
		X}_{\alpha ,\sigma ,p}=\dot{X}_{\alpha ,\sigma ,p}(\mathbb{R}^{n})$ denote
	the set of $p$-admissible weight sequences $\{t_{k}\}$ satisfying the
	following conditions. There exist numbers $C_{1},C_{2}>0$ such that for any $%
	k\leq j$\ and every cube $Q,$%
	\begin{eqnarray}
	M_{Q,p}(t_{k})M_{Q,\sigma _{1}}(t_{j}^{-1}) &\leq &C_{1}2^{\alpha _{1}(k-j)},
	\label{Asum1} \\
	M_{Q,p}^{-1}(t_{k})M_{Q,\sigma _{2}}(t_{j}) &\leq &C_{2}2^{\alpha _{2}(j-k)}.
	\label{Asum2}
	\end{eqnarray}
\end{defn}

The constants $C_{1},C_{2}>0$ are independent of both the indexes $k$ and $j$%
.

\begin{rem}
	$\mathrm{(i)}$\ We would like to mention that if $\{t_{k}\}$ satisfies $%
	\mathrm{\eqref{Asum1}}$ with $\sigma _{1}=r\left( p/r\right) ^{\prime }$ and 
	$0<r<p<\infty $, then $t_{k}^{p}\in A_{p/r}(\mathbb{R}^{n})$ for any $k\in 
	\mathbb{Z}$.\newline
	$\mathrm{(ii)}$ We say that $t_{k}\in A_{p}(\mathbb{R}^{n})$,\ $k\in \mathbb{%
		Z}$, $1<p<\infty ,$ have the same Muckenhoupt constant if%
	\begin{equation*}
	A_{p}(t_{k})=c,\quad k\in \mathbb{Z},
	\end{equation*}%
	where $c$ is independent of $k$.\newline
	$\mathrm{(iii)}$ Definition \ref{Tyulenev-class} is different from the one
	used in \cite[Definition 2.1]{Ty14} and Definition 2.7 in \cite{Ty15},
	because we used the boundedness of the maximal function on weighted Lebesgue
	spaces.
\end{rem}

\begin{ex}
	\label{Example1}Let $0<r<p<\infty $, a weight $\omega ^{p}\in A_{\frac{p}{r}%
	}(\mathbb{R}^{n})$ and
	\begin{equation*}
	\{s_{k}\}=\{2^{ks}\omega ^{p}(2^{-k}\cdot )\}_{k\in \mathbb{Z}},\quad s\in \mathbb{R}.
	\end{equation*}
	Clearly, $\{s_{k}\}_{k\in \mathbb{Z}}$ lies
	in $\dot{X}_{\alpha ,\sigma ,p}$ for $\alpha _{1}=\alpha _{2}=s$, $\sigma
	=(r(p/r)^{\prime },p)$.
\end{ex}

\begin{rem}
	\label{Tyulenev-class-properties}Let $0<\theta \leq p<\infty $. Let $\alpha
	_{1}$, $\alpha _{2}\in \mathbb{R}$, $\sigma _{1},\sigma _{2}\in (0,+\infty )$%
	, $\sigma _{2}\geq p$, $\alpha =(\alpha _{1},\alpha _{2})$ and let $\sigma
	=(\sigma _{1}=\theta \left( \frac{p}{\theta }\right) ^{\prime },\sigma _{2})$%
	. Let a $p$-admissible weight sequence $\{t_{k}\}\in \dot{X}_{\alpha ,\sigma
		,p}$. Then%
	\begin{equation*}
	\alpha _{2}\geq \alpha _{1},
	\end{equation*}%
	see \cite{D20}.
\end{rem}

%\begin{proof}

Further notation will be properly introduced whenever needed.

\subsection{Auxiliary results}

In this subsection we present some results which are useful for us. Let
recall the vector-valued maximal inequality of Fefferman and Stein \cite%
{FeSt71}. As usual, we put%
\begin{equation*}
\mathcal{M}(f)(x)=\sup_{Q}M_{Q}(f),\quad f\in L_{1}^{\mathrm{loc}}(\mathbb{R}%
^{n}),
\end{equation*}%
where the supremum\ is taken over all cubes with sides parallel to the axis
and $x\in Q$. Also we set 
\begin{equation*}
\mathcal{M}_{\sigma }(f)=\sup_{Q}M_{Q,\sigma }(f),\quad 0<\sigma <\infty .
\end{equation*}%
Observe that $\mathcal{M}_{\sigma }(f)$ can be rewritten as%
\begin{equation*}
\mathcal{M}_{\sigma }(f)=(\mathcal{M}(|f|^{\sigma }))^{1/\sigma },\quad
0<\sigma <\infty .
\end{equation*}

\begin{thm}
	\label{Maximal}Let $1<p\leq \infty $. Then%
	\begin{equation*}
	\big\|\mathcal{M}(f)|L_{p}(\mathbb{R}^{n})\big\|\lesssim \big\|f|L_{p}(%
	\mathbb{R}^{n})\big\|
	\end{equation*}%
	holds for all $f\in L_{p}(\mathbb{R}^{n})$.
\end{thm}

For the proof see \cite[Chapter 7]{L. Graf14}. Now, we state the
vector-valued maximal inequality of Fefferman and Stein \cite{FeSt71}.

\begin{thm}
	\label{FS-inequality}Let $0<p<\infty ,0<q\leq \infty $ and $0<\sigma <\min
	(p,q)$. Then%
	\begin{equation}
	\Big\|\Big(\sum\limits_{k=-\infty }^{\infty }\big(\mathcal{M}_{\sigma
	}(f_{k})\big)^{q}\Big)^{1/q}|L_{p}(\mathbb{R}^{n})\Big\|\lesssim \Big\|\Big(%
	\sum\limits_{k=-\infty }^{\infty }\left\vert f_{k}\right\vert ^{q}\Big)%
	^{1/q}|L_{p}(\mathbb{R}^{n})\Big\|  \label{Fe-St71}
	\end{equation}%
	holds for all sequence of functions $\{f_{k}\}\in L_{p}(\ell _{q})$.
\end{thm}

We shall require the following theorem, that is the Fefferman-Stein's inequality, see \cite%
{FeSt71}.

\begin{lem}
	\label{FS-lemma}Let $1<p<\infty $. Given a non-negative real valued
	functions $f$ and $g$. We have%
	\begin{equation*}
	\int_{\mathbb{R}^{n}}\left( \mathcal{M}(f)(x)\right) ^{p}g(x)dx\leq c\int_{%
		\mathbb{R}^{n}}\left( f(x)\right) ^{p}\mathcal{M}(g)(x)dx,
	\end{equation*}%
	with $c$ independent of $f$ and $g$.
\end{lem}

For the proof see \cite[Chapter 7]{L. Graf14}. We need the following version
of the Calder\'{o}n-Zygmund covering lemma, see \cite[Lemma 3.3]{CP99}, \cite%
[Appendix A]{CMP11} and \cite[Chapter 7]{Ca99}.

\begin{lem}
	\label{CZ-lemma}Let $f$ be a measurable function such that $%
	M_{Q}(f)\rightarrow 0$ as $|Q|\rightarrow \infty $ and given\ a positive
	number $a$ such that $a>2^{n+1}$. For each $i\in \mathbb{Z}$\ there exists a
	disjoint collection of maximal dyadic cubes $\{Q^{i,h}\}_{h}$ such that for
	each $h$,%
	\begin{equation*}
	a^{i}\leq M_{Q^{i,h}}(f)\leq 2^{n}a^{i}
	\end{equation*}%
	and%
	\begin{equation*}
	\Omega _{i}=\{x\in \mathbb{R}^{n}:\mathcal{M}(f)(x)>4^{n}a^{i}\}\subset \cup
	_{h}3Q^{i,h}.
	\end{equation*}%
	Let 
	\begin{equation*}
	E^{i}=\cup _{h}Q^{i,h}
	\end{equation*}%
	and 
	\begin{equation*}
	E^{i,h}=Q^{i,h}\backslash (Q^{i,h}\cap E^{i+1}). 
	\end{equation*}
	Then $E^{i,h}\subset
	Q^{i,h},$ there exists a constant $\beta >1$, depending only on $a$, such\ $%
	|Q^{i,h}|\leq \beta |E^{i,h}|$\ and the sets $E^{i,h}$ are pairwise disjoint
	for all $i$ and $h$.
\end{lem}

Next, we recall the following Hadamard's three line theorem for subharmonic
functions, see \cite[Theorem 14.15]{M83}.

\begin{thm}
	\label{hadamard}Let $f$ be a nonnegative, bounded function on the strip $%
	0\leq \Re (z)\leq 1$ of the complex plane such that $\log f(z)$ is
	subharmonic in the open strip $0<\Re (z)<1$ and continuous on the strip $%
	0\leq \Re (z)\leq 1$. If there are positive constants $M_{1},M_{2}$, such
	that $f(0+iy)\leq M_{1}$ and $f(1+iy)\leq M_{2}$ for every real $y$, then 
	\begin{equation*}
	f(\theta +iy)\leq M_{1}^{1-\theta }M_{2}^{\theta }
	\end{equation*}%
	for every $\theta \in \lbrack 0,1]$ and any real $y$.
\end{thm}

Let $S$ be the linear space consisting of all sequences $\{f_{k}\}$ with $%
f_{k}$ simple function on $\mathbb{R}^{n}$ for each $k,$ and $f_{k}=0$ for $%
|k|$ large enough. Then $S$ is dense in $L_{p}(\ell _{q}),1<p,q<\infty $,
see \cite{BP61}.

The main aim of the following lemma is to extend an interpolation result for
sublinear operators on Lebesgue spaces\ obtained by Calder\'{o}n and Zygmund
in \cite{CZ56} to similar operators acting on $L_{p}(\ell _{q},\{t_{k}\})$ spaces.

Consider a mapping function $T$, which maps measurable functions on $\mathbb{R}^{n}$ into measurable functions on $\mathbb{R}^{n}$. We say that $T$ is sublinear if it is satisfies\\
(i) $T(f)$ is defined (uniquely) if $f=f_{1}+f_{2}$, $T(f_{1})$ and $%
T(f_{2})$ are defined, and
\begin{equation*}
|T(f_{1}+f_{2})|\leq |T(f_{1})|+|T(f_{2})|.
\end{equation*}
(ii) For any constant $c$, $T(cf)$ is defined if $T(f)$ is defined and $%
|T(cf)|=|c||T(f)|$.

\begin{lem}
	\label{Calderon-Zygmund}Let $t_{k},k\in \mathbb{Z}$ be locally integrable
	functions on $\mathbb{R}^{n}$ and $1<q_{i}<\infty $, $1<p_{i}<\infty
	,M_{i}>0,i=0,1$. Suppose that $T$ is a sublinear operator satisfying 
	\begin{equation*}
	\big\|\{t_{k}T(f_{k})\}|L_{p_{i}}(\ell _{q_{i}})\big\|\leq M_{i}\big\|%
	\{t_{k}f_{k}\}|L_{p_{i}}(\ell _{q_{i}})\big\|
	\end{equation*}%
	for any $\{t_{k}f_{k}\}\in L_{p_{i}}(\ell _{q_{i}}),i=0,1$. Then $T$ can be
	extended to a bounded operator:%
	\begin{equation}
	\big\|\{t_{k}T(f_{k})\}|L_{p}(\ell _{q})\big\|\leq M_{1}^{1-\theta
	}M_{2}^{\theta }\big\|\{t_{k}f_{k}\}|L_{p}(\ell _{q})\big\|  \label{CZ1956}
	\end{equation}%
	for any $\{t_{k}f_{k}\}\in L_{p}(\ell _{q})$, where%
	\begin{equation*}
	\frac{1}{p}=\frac{1-\theta }{p_{0}}+\frac{\theta }{p_{1}}\quad \text{and}%
	\quad \frac{1}{q}=\frac{1-\theta }{q_{0}}+\frac{\theta }{q_{1}},\quad 0\leq
	\theta \leq 1.
	\end{equation*}
\end{lem}

\begin{proof}
	We will do the proof into two steps.
	
	\textit{Step 1. }We prove that \eqref{CZ1956} holds for any sequence of
	simple functions $\{t_{k}f_{k}\}\in S$.
	
	\textit{Substep 1.1. Preparation}. Assume that%
	\begin{equation*}
	\big\|\{t_{k}f_{k}\}|L_{p}(\ell _{q})\big\|=1
	\end{equation*}%
	and we put%
	\begin{equation*}
	I=\sum\limits_{k=-\infty }^{\infty }\int_{\mathbb{R}%
		^{n}}t_{k}(x)|T(f_{k})(x)|g_{k}(x)dx,
	\end{equation*}%
	where $\{g_{k}\}\in S$ such that%
	\begin{equation*}
	\big\|\{g_{k}\}|L_{p^{\prime }}(\ell _{q^{\prime }})\big\|\leq 1.
	\end{equation*}%
	From our assumption on the sequence of simple functions $\{t_{k}f_{k}\}$ and 
	$\{g_{k}\}$, we clearly have that%
	\begin{equation*}
	t_{k}f_{k}=\sum_{l=1}^{N_{1}}|a_{l,k}|c_{l,k}\chi _{E_{l,k}}\quad \text{and}%
	\quad g_{k}=\sum_{j=1}^{N_{2}}b_{j,k}\chi _{E_{j,k}^{\prime }},
	\end{equation*}%
	where $|c_{l,k}|=1,b_{j,k}>0,l=1,...,N_{1},j=1,...,N_{2},N_{1},N_{2}\in 
	\mathbb{N}$ and 
	\begin{equation*}
	t_{k}f_{k}=g_{k}=0
	\end{equation*}%
	for sufficiently large $k,|k|\geq N,N\in \mathbb{N}$. For $k$ fixed, the
	sets $E_{l,k}$ are\ disjoint. The same property holds for the sets $%
	E_{j,k}^{\prime }$. For any $z\in \mathbb{C}$ we set%
	\begin{equation*}
	\frac{1}{q(z)}=\frac{1-z}{q_{0}}+\frac{z}{q_{1}},\quad \frac{1}{p(z)}=\frac{%
		1-z}{p_{0}}+\frac{z}{p_{1}}
	\end{equation*}%
	and 
	\begin{align*}
	&\Phi (z)\\
	&=\sum\limits_{|k|\leq N}\int_{\mathbb{R}^{n}}t_{k}(x)\big|T\big(%
	\omega (\cdot ,z)|f_{k}|^{q/q(z)}\mathrm{sign}f_{k}\big)(x)\big|%
	(g_{k}(x))^{\beta (\Re (z))}\vartheta (x,\Re (z))dx,
	\end{align*}%
	where%
	\begin{equation*}
	\omega (\cdot ,z)=A^{\tau (z)}(\cdot )t_{k}^{-\alpha (z)}(\cdot ),\quad A(\cdot )=%
	\big\|\{t_{k}f_{k}\}_{|k|\leq N}|\ell _{q}\big\|
	\end{equation*}%
	and 
	\begin{equation*}
	\vartheta (\cdot ,z)=B^{\kappa (z)-\beta (z)}(\cdot ),\quad B(\cdot )=\big\|%
	\{g_{k}\}_{|k|\leq N}|\ell _{q^{\prime }}\big\|,
	\end{equation*}%
	with%
	\begin{equation*}
	\tau (z)=\frac{p}{p(z)}-\frac{q}{q(z)},\quad \alpha (z)=1-\frac{q}{q(z)}%
	,\quad \kappa (z)=\frac{1-\frac{1}{p(z)}}{1-\frac{1}{p}},\quad \beta (z)=%
	\frac{1-\frac{1}{q(z)}}{1-\frac{1}{q}}.
	\end{equation*}%
	The non-negative function $\Phi $ reduces to $I$ for $z=\theta $. Using the
	fact that $\{g_{k}\}$ is a sequence of simple functions we have 
	\begin{equation*}
	\Phi (z)=\sum_{j=1}^{N_{2}}\sum\limits_{|k|\leq N}\int_{E_{j,k}^{\prime
	}}\vartheta (x,\Re (z))t_{k}(x)\left\vert T(\psi _{z,j,k})(x)\right\vert dx,
	\end{equation*}%
	where%
	\begin{equation*}
	\psi _{z,j,k}=\omega (\cdot ,z)(b_{j,k})^{\beta (\Re
		(z))}\sum_{l=1}^{N_{1}}|t_{k}^{-1}a_{l,k}|^{q/q(z)}\chi _{E_{l,k}}.
	\end{equation*}%
	We put%
	\begin{equation*}
	\Psi _{j,k}(z)=\int_{E_{j,k}^{\prime }}\vartheta (x,\Re
	(z))t_{k}(x)\left\vert T(\psi _{z,j,k})(x)\right\vert dx,\quad |k|\leq
	N,j=1,...,N_{2}.
	\end{equation*}%
	
	\textit{Substep 1.2. }We prove that $\Phi $ is continuous in the strip $%
	0\leq \Re (z)\leq 1$. First we prove that $\Psi _{j,k}$ is well defined. Let
	us consider the integrals%
	\begin{equation*}
	I_{z,j,k}^{1}=\int_{E_{j,k}^{\prime }}t_{k}^{p_{0}}(x)\big|T\big(\omega
	(\cdot ,z)|f_{k}|^{q/q(z)}\mathrm{sign}f_{k}\big)(x)\big|^{p_{0}}dx
	\end{equation*}%
	and 
	\begin{equation*}
	I_{z,j,k}^{2}=\int_{E_{j,k}^{\prime }}|\vartheta (x,\Re (z))|^{p_{0}^{\prime
	}}(g_{k}(x))^{p_{0}^{\prime }\beta (\Re (z))}dx.
	\end{equation*}%
	An easy calculation shows that the integral%
	\begin{equation*}
	J=\int_{\mathbb{R}^{n}}\Big(\sum\limits_{|k|\leq N}t_{k}^{q_{0}}(x)\big|%
	\omega (x,z)|f_{k}(x)|^{q/q(z)}\big|^{q_{0}}\Big)^{p_{0}/q_{0}}dx
	\end{equation*}%
	is just%
	\begin{equation*}
	\int_{\mathbb{R}^{n}}A^{p_{0}\tau (\Re (z))}(x)\Big(\sum\limits_{|k|\leq N}%
	\big(t_{k}(x)|f_{k}(x)|\big)^{(1-\alpha (\Re (z)))q_{0}}\Big)%
	^{p_{0}/q_{0}}dx.
	\end{equation*}%
	We may assume without loss of generality that $q_{1}<q_{0}$. As a
	consequence we see that%
	\begin{equation*}
	\frac{1}{q}(1-\alpha (\Re (z)))=\frac{1-\Re (z)}{q_{0}}+\frac{\Re (z)}{q_{1}}%
	>\frac{1}{q_{0}}.
	\end{equation*}%
	This yields%
	\begin{align*}
	&\Big(\sum\limits_{|k|\leq N}t_{k}^{(1-\alpha (\Re
		(z)))q_{0}}(x)|f_{k}(x)|^{(1-\alpha (\Re (z)))q_{0}}\Big)^{1/q_{0}}\\
	&\leq \Big(\sum\limits_{|k|\leq N}t_{k}^{q}(x)|f_{k}(x)|^{q}\Big)^{(1-\alpha (\Re
		(z)))/q},
	\end{align*}%
	since $\ell _{r}\hookrightarrow \ell _{s},r\leq s$. Hence%
	\begin{equation*}
	J\leq \big\|\{t_{k}f_{k}\}|L_{pp_{0}/p(\Re (z))}(\ell _{q})\big\|%
	^{pp_{0}/p(\Re (z))}<\infty ,
	\end{equation*}%
	because of $\{t_{k}f_{k}\}$ is a sequence of simple functions. Consequently,
	by H\"{o}lder's inequality,%
	\begin{align}
	\Psi _{j,k}(z) &\leq (I_{z,j,k}^{1})^{1/p_{0}}(I_{z,j,k}^{2})^{1/p_{0}^{\prime }}
	\notag \\
	&\leq \big\|\{t_{k}T\big(\omega (\cdot ,z)|f_{k}|^{q/q(z)}\mathrm{sign}f_{k}%
	\big)\}|L_{p_{0}}(\ell _{q_{0}})\big\|(I_{z,j,k}^{2})^{1/p_{0}^{\prime
	}}.  \label{firstterm}
	\end{align}%
	Our assumption on $T$ yields that the first term in $\mathrm{%
		\eqref{firstterm}}$ is finite.\ $I_{z,j,k}^{2}$ is finite since $%
	\{g_{k}\}\in S$. Hence $\Phi (z)$ exists for each $z$ in the strip $0\leq
	\Re (z)\leq 1$.
	
	Now we prove the continuity of $\Phi $ in the strip $0\leq \Re (z)\leq 1$.
	We clearly have that%
	\begin{equation}
	\Psi _{j,k}(z+\Delta z)-\Psi _{j,k}(z)=J_{1,j,k}(z,\Delta
	z)+J_{2,j,k}(z,\Delta z),  \label{limit}
	\end{equation}%
	where%
	\begin{equation*}
	J_{1,j,k}(z,\Delta z)=\int_{E_{j,k}^{\prime }}\big(\vartheta (x,\Re
	(z+\Delta z))-\vartheta (x,\Re (z))\big)t_{k}(x)\left\vert T(\psi
	_{z,j,k})(x)\right\vert dx
	\end{equation*}%
	and%
	\begin{align*}
	&	J_{2,j,k}(z,\Delta z)\\
	&=\int_{E_{j,k}^{\prime }}\vartheta (x,\Re (z+\Delta
	z))t_{k}(x)\big(\left\vert T(\psi _{z+\Delta z,j,k})(x)\right\vert
	-\left\vert T(\psi _{z,j,k})(x)\right\vert \big)dx.
	\end{align*}%
	We want to find the limit of $\mathrm{\eqref{limit}}$ as $\Delta z$
	approaches $0$, so we may assume that $|\Delta z=\Re (\Delta z)+i\Im (\Delta
	z)|\leq 1$. After a simple calculation we find that%
	\begin{align*}
	&\big|\vartheta (x,\Re (z+\Delta z))-\vartheta (x,\Re (z))\big| \\
	&=\big|%
	B^{\kappa (\Re (z))-\beta (\Re (z))}(x)\big(B^{d+h\Re (\Delta z)}(x)-1\big)%
	\big| \\
	&\leq B^{\kappa (\Re (z))-\beta (\Re (z))}(x)\big(B^{d+h\Re (\Delta z)}(x)+1%
	\big) \\
	&\leq B^{\kappa (\Re (z))-\beta (\Re (z))}(x)\big(B^{d}(x)\max
	(1,B^{h}(x),B^{-h}(x))+1\big),
	\end{align*}%
	where%
	\begin{equation*}
	h=\frac{\frac{1}{q_{1}}-\frac{1}{q_{0}}}{1-\frac{1}{q}}+\frac{\frac{1}{p_{0}}%
		-\frac{1}{p_{1}}}{1-\frac{1}{p}},\quad d=\frac{1-\frac{1}{p_{0}}}{1-\frac{1}{p}}-\frac{1-%
		\frac{1}{q_{0}}}{1-\frac{1}{q}}.
	\end{equation*}%
	Recall that 
	\begin{equation*}
	\big|T(\psi _{z,j,k})(x)\big|=(b_{j,k})^{\beta (\Re (z))}\big|T\big(\omega
	(\cdot ,z)|f_{k}|^{\frac{q}{q(z)}}\mathrm{sign}f_{k}\big)(x)\big|.
	\end{equation*}%
	Therefore the function 
	\begin{equation*}
	x\longrightarrow \big|\vartheta (x,\Re (z+\Delta z))-\vartheta (x,\Re (z))%
	\big|t_{k}(x)\left\vert T(\psi _{z,j,k})(x)\right\vert
	\end{equation*}%
	is integrable. Dominated convergence theorem yields that $J_{1,j,k}(z,\Delta
	z)$ tends to zero as $\Delta z$ tends to $0$. Now%
	\begin{equation*}
	|J_{2,j,k}(z,\Delta z)|\leq \int_{E_{j,k}^{\prime }}\vartheta (x,\Re
	(z+\Delta z))t_{k}(x)\big|T(\psi _{z+\Delta z,j,k}-\psi _{z,j,k})(x)\big|dx,
	\end{equation*}%
	since%
	\begin{equation*}
	|T(\psi _{z+\Delta z,j,k})|-|T(\psi _{z,j,k})|\leq |T(\psi _{z+\Delta
		z,j,k}-\psi _{z,j,k})|,\quad j\in \{1,...,N_{2}\},|k|\leq N.
	\end{equation*}%
	Hence%
	\begin{equation*}
	|J_{2,j,k}(z,\Delta z)|\leq \Lambda _{j,k}\times B_{j,k},\quad j\in
	\{1,...,N_{2}\},|k|\leq N,
	\end{equation*}%
	where%
	\begin{equation*}
	\Lambda _{j,k}=(b_{j,k})^{\beta (\Re (z))}\Big(\int_{E_{j,k}^{\prime
	}}\vartheta ^{p_{0}^{\prime }}(x,\Re (z+\Delta z))dx\Big)^{1/p_{0}^{\prime }}
	\end{equation*}%
	and%
	\begin{equation*}
	B_{j,k}=\Big(\int_{E_{j,k}^{\prime }}t_{k}^{p_{0}}(x)\big|T\big((b_{j,k})^{-\beta
		(\Re (z))}(\psi _{z+\Delta z,j,k}-\psi _{z,j,k})\big)(x)\big|^{p_{0}}dx\Big)%
	^{1/p_{0}}.
	\end{equation*}%
	Clearly, $\Lambda _{j,k},j\in \{1,...,N_{2}\},|k|\leq N$ is bounded. To
	estimate $B_{j,k}$ we have%
	\begin{align}
	& \big|\psi _{z+\Delta z,j,k}\big|  \notag \\
	& =\big|\omega (\cdot ,z+\Delta z)\big|(b_{j,k})^{\beta (\Re (z+\Delta
		z))}\sum_{l=1}^{N_{1}}\Big|\big|t_{k}^{-1}a_{l,k}\big|^{q/q(z+\Delta z)}\Big|%
	\chi _{E_{l,k}}  \notag \\
	& =\big|\psi _{z,j,k}\big|A^{\sigma \Re (\Delta z)}(b_{j,k})^{q^{\prime
		}(1/q_{0}-1/q_{1})\Re (\Delta z)}\sum_{l=1}^{N_{1}}\big|a_{l,k}\big|%
	^{q(1/q_{1}-1/q_{0})\Re (\Delta z)}\chi _{E_{l,k}},  \label{simple}
	\end{align}%
	where%
	\begin{equation*}
	\sigma =p\Big(\frac{1}{p_{1}}-\frac{1}{p_{0}}\Big)+q\Big(\frac{1}{q_{0}}-%
	\frac{1}{q_{1}}\Big).
	\end{equation*}%
	There are two cases: $0\leq \Re (\Delta z)\leq 1$ and $-1\leq \Re (\Delta
	z)<0$. In the first case we obtain the estimates%
	\begin{equation}
	\sum_{l=1}^{N_{1}}|a_{l,k}|^{q(1/q_{1}-1/q_{0})\Re (\Delta z)}\chi
	_{E_{l,k}}=A^{q(1/q_{1}-1/q_{0})\Re (\Delta z)}  \label{ref1}
	\end{equation}%
	and%
	\begin{equation}
	(b_{j,k})^{q^{\prime }(1/q_{0}-1/q_{1})\Re (\Delta z)}\leq \max
	\big(1,(b_{j,k})^{q^{\prime }(1/q_{0}-1/q_{1})}\big).  \label{ref2}
	\end{equation}%
	Indeed, the left-hand side of $\mathrm{\eqref{ref1}}$ is equal to 
	\begin{equation*}
	\Big(\sum_{l=1}^{N_{1}}|a_{l,k}c_{l,k}|^{q}\chi _{E_{l,k}}\Big)%
	^{(1/q_{1}-1/q_{0})\Re (\Delta z)}=A^{q(1/q_{1}-1/q_{0})\Re (\Delta z)}.
	\end{equation*}%
	Now, if $b_{j,k}\geq 1,j\in \{1,...,N_{2}\},|k|\leq N$, then%
	\begin{equation*}
	(b_{j,k})^{q^{\prime }(1/q_{0}-1/q_{1})\Re (\Delta z)}\leq 1,
	\end{equation*}%
	since $q_{1}<q_{0}$ and $0\leq \Re (\Delta z)\leq 1$. While if $0<b_{j,k}<1$%
	, then%
	\begin{align*}
	&	(b_{j,k})^{q^{\prime }(1/q_{0}-1/q_{1})\Re (\Delta z)}\notag\\
	&=(b_{j,k})^{q^{\prime
		}(1/q_{0}-1/q_{1})(\Re (\Delta z)-1)}(b_{j,k})^{q^{\prime
		}(1/q_{0}-1/q_{1})}\\
	&\leq (b_{j,k})^{q^{\prime }(1/q_{0}-1/q_{1})}.
	\end{align*}%
	Thus, we obtain $\mathrm{\eqref{ref2}}$\textrm{. }Assume that $-1\leq \Re
	(\Delta z)<0$. We obtain%
	\begin{equation*}
	(b_{j,k})^{q^{\prime }(1/q_{0}-1/q_{1})\Re (\Delta z)}\leq \max
	(1,(b_{j,k})^{q^{\prime }(1/q_{1}-1/q_{0})}),
	\end{equation*}
	\begin{equation*}
	|a_{l,k}|^{q(1/q_{1}-1/q_{0})\Re (\Delta z)}\leq \max
	(1,|a_{l,k}|^{q(1/q_{0}-1/q_{1})})
	\end{equation*}%
	and%
	\begin{equation*}
	A^{q(1/q_{0}-1/q_{1})\Re (\Delta z)}\leq \max
	(1,A^{q(1/q_{1}-1/q_{0})}),
	\end{equation*}
	\begin{equation*}
	A^{\sigma \Re (\Delta z)}\leq \max(1,A^{\sigma },A^{-\sigma }).
	\end{equation*}%
	Substituting these estimations in $\mathrm{\eqref{simple}}$ we find that%
	\begin{align*}
	&	|\psi _{z+\Delta z,j,k}|\\
	&	\leq |\psi _{z,j,k}|\max \Big(1,\max
	(1,A^{q(1/q_{1}-1/q_{0})}),\max (1,A^{\sigma },A^{-\sigma })\Big)\mu _{j,k},
	\notag
	\end{align*}%
	where $\{\mu _{j,k}\}\in S$. This estimate together with 
	\begin{equation*}
	\{t_{k}(b_{j,k})^{-\beta (\Re (z))}\psi _{z,j,k}\}_{|k|\leq N}\in
	L_{p_{0}}(\ell _{q_{0}})
	\end{equation*}%
	guarantees that the function 
	\begin{align*}
	&	x\longrightarrow t_{k}(x)(b_{j,k})^{-\beta (\Re (z))}|\psi _{z,j,k}(x)|\\
	& \times\max \Big(1,\max (1,A^{q(1/q_{1}-1/q_{0})}(x)),\max (1,A^{\sigma }(x),A^{-\sigma }(x))\Big)%
	\mu _{j,k}(x)
	\end{align*}%
	belongs to $L_{p_{0}}(\ell _{q_{0}})$, $|k|\leq N$, $j\in \{1,...,N_{2}\}$. Hence%
	\begin{align*}
	& \Big(\int_{E_{j,k}^{\prime }}t_{k}^{p_{0}}(x)|T((b_{j,k})^{-\beta (\Re
		(z))}(\psi _{z+\Delta z,j,k}-\psi _{z,j,k}))(x)|^{p_{0}}dx\Big)^{1/
		p_{0}}  \notag \\
	& \lesssim \Big(\int_{\mathbb{R}^{n}}\Big(\sum\limits_{|k|\leq N}t_{k}^{q_{0}}(x)%
	\big|(b_{j,k})^{-\beta (\Re (z))}(\psi _{z+\Delta z,j,k}(x)-\psi _{z,j,k}(x))%
	\big|^{q_{0}}\Big)^{p_{0}/q_{0}}dx\Big)^{1/p_{0}}.
	\end{align*}%
	Dominated convergence theorem yields that $J_{2,j,k}(z,\Delta z)$ tends to
	zero as $\Delta z$ tends to $0$\ and $\Psi _{j,k}$ is continuous at $z$.
	Hence $\Phi $ is continuous in the strip $0\leq \Re (z)\leq 1$.
	
	\textit{Substep 1.3.} Here we prove that $\Phi $ is bounded in the strip $%
	0\leq \Re (z)\leq 1$. From Substep 1.2, we need only to prove that $%
	I_{z,j,k}^{2}\leq M$, for a suitable constant $M$ independent of $z$. Since 
	\begin{equation*}
	B^{\kappa (\Re (z))-\beta (\Re (z))}(x)\leq \max (1,B^{p^{\prime
	}}(x))(g_{k}(x))^{-\beta (\Re (z))}
	\end{equation*}%
	for any $x\in \mathbb{R}^{n}$ and $|k|\leq N$, we find that%
	\begin{equation*}
	I_{z,j,k}^{2}\leq \int_{E_{j,k}^{\prime }}\max (1,B^{p^{\prime }}(x))dx\leq M,
	\end{equation*}%
	we are done.
	
	\textit{Substep 1.4. }In this step we prove that $\log \Psi _{j,k}(z)$ is
	subharmonic. We fix a harmonic function $h(z)$, and denote by $H(z)$ the
	analytic function whose real part is $h(z)$. We need only to prove that $%
	\Psi _{j,k}(z)e^{h(z)}$ is subharmonic for every harmonic $h(z)$. Since the
	problem is local, we may consider $h$ and $H$ in a given circle. Put%
	\begin{equation*}
	\psi _{z,j,k}^{\star }=\psi _{z,j,k}e^{H(z)}\quad \text{and}\quad \Psi
	_{j,k}^{\star }(z)=\Psi _{j,k}(z)e^{h(z)}.
	\end{equation*}%
	We fix $z$, take a $\varrho >0$, and denote by $z_{1},...,z_{r}$ be a system
	of points equally spaced over the circumference of the circle with center $z$
	and radius $\varrho $. Our estimate use partially some techniques already
	used in \cite{Kr72}. First we prove that $\log \big|T(\psi _{z,j,k})\big|$
	is subharmonic. We want to show that%
	\begin{equation*}
	\big|T(\psi _{z,j,k}^{\star })\big|=e^{h(z)}\big|T(\psi _{z,j,k})\big|\leq 
	\frac{1}{2\pi }\int_{0}^{2\pi }\big|T(\psi _{z+\varrho e^{it},j,k}^{\ast })%
	\big|dt.
	\end{equation*}%
	Let us calculate the limit of 
	\begin{equation}
	\big\|t_{k}T\big(\psi _{z,j,k}^{\star }-\frac{1}{r}\sum_{v=1}^{r}\psi
	_{z_{v},j,k}^{\star }\big)\big\|_{p_{1}}.  \label{estimate}
	\end{equation}%
	as $r$ tends to infinity. We estimate $|\psi _{z_{v},j,k}^{\star
	}|,v=1,...,r $. Clearly%
	\begin{equation*}
	|\psi _{z_{v},j,k}^{\star }|\leq C_{2}t_{k}^{-1}A^{\tau (\Re
		(z_{v}))}\sum_{l=1}^{N_{1}}\chi _{E_{l,k}}.
	\end{equation*}%
	We consider separately the possibilities $\tau (\Re (z_{v}))\geq 0$\ and $%
	\tau (\Re (z_{v}))<0$. In the first case%
	\begin{equation*}
	A^{\tau (\Re (z_{v}))}\leq C_{3} \big\|\{|t_{k}f_{k}|^{\tau (\Re
		(z_{v}))}\}_{|k|\leq N}|\ell _{q}\big\|\leq C_{4}\big\|\big\{%
	\sum_{l=1}^{N_{1}}\chi _{E_{l,k}}\big\}_{|k|\leq N}|\ell _{q}\big\|=C_{4}K.
	\end{equation*}%
	If $\tau (\Re (z_{v}))<0$, then 
	\begin{equation*}
	A^{\tau (\Re (z_{v}))}\leq |t_{0}f_{0}|^{\tau (\Re (z_{v}))}\leq
	C_{5}\sum_{l=1}^{N_{1}}\chi _{E_{l,0}}.
	\end{equation*}%
	Therefore%
	\begin{equation*}
	|\psi _{z_{v},j,k}^{\star }|\lesssim \max (1,K)t_{k}^{-1}g_{j,k}
	\end{equation*}%
	where $\{g_{j,k}\}_{k\in \mathbb{Z}}\in S$. Notice that the implicit
	constant independent of $r,v,j$ and $k$. This guarantees that $\mathrm{%
		\eqref{estimate}}$ can be estimated by%
	\begin{equation*}
	C\big\|\psi _{z,j,k}^{\star }-\frac{1}{r}\sum_{v=1}^{r}\psi
	_{z_{v},j,k}^{\star }\big\|_{p_{1}}.
	\end{equation*}%
	In view of the fact that 
	\begin{equation*}
	\psi _{z,j,k}^{\star }(x)=\lim_{r\rightarrow \infty }\frac{1}{r}%
	\sum_{v=1}^{r}\psi _{z_{v},j,k}^{\star }(x),
	\end{equation*}%
	which is valid since $\psi _{z,j,k}^{\star }$ is analytic for each $x$,
	dominated convergence theorem yields that the last norm tends to zero as $r$
	tends to infinity. Therefore $\mathrm{\eqref{estimate}}$ tends to zero as $r$
	tends to infinity. There exist subsequences 
	\begin{equation*}
	\Big\{T\big(\psi _{z,j,k}^{\star }-\frac{1}{r_{l}}\sum_{v=1}^{r_{l}}\psi
	_{z_{v},j,k}^{\star }\big)\Big\}_{l}
	\end{equation*}%
	converging to zero as $r_{l}$ tends to infinity. Hence%
	\begin{align*}
	\big|T(\psi _{z,j,k}^{\star })\big| &\leq \lim_{r_{l}\rightarrow \infty }%
	\big|T(\frac{1}{r_{l}}\sum_{v=1}^{r_{l}}\psi _{z_{v},j,k}^{\star })\big| \\
	&\leq \lim_{r_{l}\rightarrow \infty }\frac{1}{r_{l}}\sum_{v=1}^{r_{l}}\big|%
	T(\psi _{z_{v},j,k}^{\star })\big| \\
	&=\frac{1}{2\pi }\int_{0}^{2\pi }\big|T(\psi _{z+\varrho e^{it},j,k}^{\star
	})\big|dt.
	\end{align*}%
	In addition, the mapping $z\rightarrow \log (\vartheta (\cdot ,\Re
	(z))t_{k}(\cdot ))$ is subharmonic. Consequently, 
	\begin{equation*}
	z\rightarrow \log (\vartheta (\cdot ,\Re (z))t_{k}(\cdot )\big|T(\psi
	_{z,j,k})\big|)
	\end{equation*}%
	is subharmonic. Then%
	\begin{align*}
	&	e^{h(z)}\vartheta (x,\Re (z))t_{k}(x)\big|T(\psi _{z,j,k})(x)\big|\\
	&\leq\frac{1}{2\pi }\int_{0}^{2\pi }\vartheta (x,\Re (z+\varrho e^{it}))t_{k}(x)\big|T%
	\big(\psi _{z+\varrho e^{it},j,k}^{\star }\big)(x)\big|dt.
	\end{align*}%
	We integrate with respect to $x$ over the set $E_{j,k}^{\prime }$, we apply
	Fubini's theorem and observe that%
	\begin{equation*}
	\Psi _{j,k}(z)e^{h(z)}=\int_{E_{j,k}^{\prime }}\vartheta (x,\Re (z))t_{k}(x)%
	\big|T(\psi _{z,j,k}^{\star })(x)\big|dx,
	\end{equation*}%
	we obtain
	
	\begin{equation*}
	\Psi _{j,k}^{\star }(z)\leq \frac{1}{2\pi }\int_{0}^{2\pi }\Psi
	_{j,k}^{\star }(z+\varrho e^{it})dt.
	\end{equation*}%
	Hence $\Psi _{j,k}^{\star }$ is subharmonic in the strip $0\leq \Re (z)\leq
	1 $.
	
	\textit{Substep 1.5. }We prove that \eqref{CZ1956} holds for any sequence of
	simple functions $\{t_{k}f_{k}\}_{k\in \mathbb{Z}}\in S$. Let $z=\Re
	(z)+i\Im (z),0\leq \Re (z)\leq 1$ and $\Im (z)\in \mathbb{R}$. From Substeps
	1.2-1.4, $\Phi $ is continuous, bounded\ its logarithm is subharmonic in the
	strip $0\leq \Re (z)\leq 1$.\ Let us show that $\Phi (i\Im (z)\leq M_{1}$
	for any $\Im (z)\in \mathbb{R}.$ We have%
	\begin{equation*}
	D^{\tau (0)}=D^{p/p_{0}-q/q_{0}},\quad D^{\kappa (0)-\beta (0)}=D^{p^{\prime
		}/p_{0}^{\prime }-q^{\prime }/q_{0}^{\prime }},\quad |D^{\alpha (i\Im
		(z))}|=D^{1-q/q_{0}}
	\end{equation*}%
	for any $D\geq 0$. By H\"{o}lder's inequality $\Phi (i\Im (z))$ can be
	estimated by 
	\begin{align*}
	&\Big(\int_{\mathbb{R}^{n}}\Big(\sum\limits_{|k|\leq N}t_{k}^{q_{0}}(x)\big|T%
	\big(\omega (\cdot ,i\Im (z))|f_{k}|^{q/q(i\Im (z))}\mathrm{sign}f_{k}\big)%
	(x)\big|^{q_{0}}\Big)^{p_{0}/q_{0}}dx\Big)^{1/p_{0}} \\
	&\times \Big(\int_{\mathbb{R}^{n}}\Big(\sum\limits_{|k|\leq N}\big(%
	(g_{k}(x))^{\beta (0)}\vartheta (x,0)\big)^{q_{0}^{\prime }}\Big)%
	^{p_{0}^{\prime }/q_{0}^{\prime }}dx\Big)^{1/p_{0}^{\prime }}.
	\end{align*}%
	The second integral is just%
	\begin{equation*}
	\Big(\int_{\mathbb{R}^{n}}\Big(\sum\limits_{|k|\leq N}g_{k}^{q^{\prime }}(x)%
	\Big)^{p^{\prime }/q^{\prime }}dx\Big)^{1/p_{0}^{\prime }}\leq 1,
	\end{equation*}%
	since $\kappa (0)-\beta (0)=p^{\prime }/p_{0}^{\prime }-q^{\prime
	}/q_{0}^{\prime }$ and $\beta (0)=q^{\prime }/q_{0}^{\prime }$. The
	boundedness of $T$ on $L_{p_{0}}(\ell _{q_{0}},\{t_{k}\})$ yields that the first integral does not exceed%
	\begin{equation}
	M_{1}\Big(\int_{\mathbb{R}^{n}}\Big(\sum\limits_{|k|\leq N}t_{k}^{q_{0}}(x)\big|%
	\omega (x,i\Im (z))|f_{k}(x)|^{q/q(i\Im (z))}\mathrm{sign}f_{k}\big|^{q_{0}}%
	\Big)^{p_{0}/q_{0}}dx\Big)^{1/p_{0}},  \label{caseRez=0}
	\end{equation}%
	but \eqref{caseRez=0} is just%
	\begin{equation*}
	M_{1}\big\|\{t_{k}f_{k}\}|L_{p}(\ell _{q})\big\|^{p/p_{0}}= M_{1}.
	\end{equation*}%
	Similarly we obtain 
	\begin{equation*}
	\Phi (1+i\Im (z))\leq M_{2}.
	\end{equation*}%
	From our previous\ substeps $\log \Phi (z)$ is continuous bounded above and
	subharmonic in the strip $0\leq \Re (z)\leq 1$. In addition $\log \Phi (z)$
	does not exceed $\log M_{1}$ and $\log M_{2}$ on the lines $\Re (z)=0$ and $%
	\Re (z)=1$, respectively. Therefore we apply the Hadamard three-line theorem
	for subharmonic functions, see Theorem \ref{hadamard}, we obtain that %
	\begin{equation*}
	\log \Phi (z)\leq (1-\theta )\log M_{1}+\theta \log M_{2},\quad 0<\Re (z)<1.
	\end{equation*}%
	In particular%
	\begin{equation*}
	I=\Phi (\theta )\leq M_{1}^{1-\theta }M_{2}^{\theta }.
	\end{equation*}%
	Finally we have proved that%
	\begin{equation*}
	\big\|\{t_{k}T(f_{k})\}|L_{p}(\ell _{q})\big\|\leq M_{1}^{1-\theta
	}M_{2}^{\theta }\big\|\{t_{k}f_{k}\}|L_{p}(\ell _{q})\big\|
	\end{equation*}%
	for each sequence of simple functions $\{t_{k}f_{k}\}$.
	
	\textit{Step 2. }We prove \eqref{CZ1956}. Assume $q_{1}<q_{0}$ and let $%
	\{f_{k}\}\in L_{p}(\ell _{q},\{t_{k}\})$. Then $\{t_{k}f_{k}\}\in L_{p}(\ell
	_{q})$. There exists a sequence $\{g_{n}\}_{n\in \mathbb{N}}=\{\{g_{n}^{k}\}_{k\in 
		\mathbb{Z}}\}_{n\in \mathbb{N}}\subset S\subset L_{p}(\ell _{q})$ such that $%
	\{g_{n}\}_{n\in \mathbb{N}}$ converges to $\{t_{k}f_{k}\}$ in $L_{p}(\ell
	_{q})$. Hence%
	\begin{equation*}
	\Big\|\{t_{k}^{-1}g_{n}^{k}-f_{k}\}|L_{p}(\ell _{q},\{t_{k}\})\Big\|
	\end{equation*}%
	tends to zero as $n$ tends to infinity. Therefore%
	\begin{equation*}
	\Big\|\{t_{k}^{-1}g_{n}^{k}-f_{k}\}|L_{p}(\ell _{q_{0}},\{t_{k}\})\Big\|
	\end{equation*}%
	tends to zero as $n$ tends to infinity, since $L_{p}(\ell
	_{q})\hookrightarrow L_{p}(\ell _{q_{0}})$. Observe that 
	\begin{equation*}
	\big|\big|T(f_{k})\big|-\big|T(t_{k}^{-1}g_{n}^{k})\big|\big|\leq \big|%
	T(f_{k}-t_{k}^{-1}g_{n}^{k})\big|
	\end{equation*}%
	for any $k\in \mathbb{Z}$ and any $n\in \mathbb{N}$. Therefore 
	\begin{equation*}
	\Big\|\{\big|T(f_{k})\big|\}-\{\big|T(t_{k}^{-1}g_{n}^{k})\big|\}|L_{p}(\ell
	_{q_{0}},\{t_{k}\})\Big\|  \label{step2}
	\end{equation*}%
	can be estimated by%
	\begin{equation*}
	\Big\|\{T(f_{k}-t_{k}^{-1}g_{n}^{k})\}|L_{p}(\ell _{q_{0}},\{t_{k}\})\Big\|%
	\lesssim \Big\|\{t_{k}^{-1}g_{n}^{k}-f_{k}\}|L_{p}(\ell _{q_{0}},\{t_{k}\})%
	\Big\|,
	\end{equation*}%
	where this inequality follows by the fact that $(t_{k}^{-1}g_{n}^{k})\}\in
	L_{p}(\ell _{q_{0}},\{t_{k}\})$ since%
	\begin{equation*}
	\Big\|\{t_{k}^{-1}g_{n}^{k})\}|L_{p}(\ell _{q_{0}},\{t_{k}\})\Big\|=\Big\|%
	\{g_{n}^{k})\}|L_{p}(\ell _{q_{0}})\Big\|<\infty .
	\end{equation*}%
	This yields that%
	\begin{equation*}
	\Big\|\{T(f_{k}-t_{k}^{-1}g_{n}^{k})\}|L_{p}(\ell _{q_{0}},\{t_{k}\})\Big\|
	\end{equation*}%
	tends to zero as $n$ tends to infinity. Therefore $\big\{\{\big|%
	T(t_{k}^{-1}g_{n}^{k})\big|\}\big\}_{n\in \mathbb{N}}$ converges to $\big\{\big|%
	T(f_{k})\big|\big\}$ in $L_{p}(\ell _{q_{0}},\{t_{k}\})$-norm, which yields that 
	$\big\{\{t_{k}\big|T(t_{k}^{-1}g_{n}^{k})\big|\}\big\}_{n\in \mathbb{N}}$ converges
	to $\big\{t_{k}\big|T(f_{k})\big|\big\}$ in $L_{p}(\ell _{q_{0}})$-norm. Hence $%
	\big\{t_{k}\big|T(t_{k}^{-1}g_{n}^{k})\big|\big\}_{n\in \mathbb{N}}$ converges to $t_{k}\big|T(f_{k})\big|$ in $L_{p}$ for every $k\in \mathbb{Z}$. The Cantor diagonal technique gives  an increasing sequence $\{\varpi _{k}(n)\}_{n\in \mathbb{N}}$ in ${\mathbb{N}}$ such that $\big\{t_{k}\big|T(t_{k}^{-1}g_{\varpi _{k}(n)}^{k})\big|\big\}_{n\in 
		\mathbb{N}}$ converges to $t_{k}\big|T(f_{k})\big|$  for every $k\in \mathbb{Z}$. We have%
	\begin{align*}
	\Big\|\{T(f_{k}\}|L_{p}(\ell _{q},\{t_{k}\})\Big\| &=\Big\|%
	\{t_{k}T(f_{k})\}|L_{p}(\ell _{q})\Big\| \\
	&=\Big\|\Big(\sum\limits_{k=-\infty }^{\infty }t_{k}^{q}\big|T(f_{k})\big|%
	^{q}\Big)^{1/q}\Big\|_{p} \\
	&=\Big\|\Big(\sum\limits_{k=-\infty }^{\infty }\lim_{n\rightarrow \infty }%
	\big(t_{k}\big|T(t_{k}^{-1}g_{\varpi _{k}(n)}^{k})\big|\big)^{q}\Big)^{1/q}\Big\|_{p} \\
	&\leq \underset{n\longrightarrow \infty }{\lim \inf }\Big\|\Big(%
	\sum\limits_{k=-\infty }^{\infty }\big(t_{k}\big|T(t_{k}^{-1}g_{\varpi
		_{k}(n)}^{k})\big|\big)^{q}\Big)^{1/q}\Big\|_{p}
	\end{align*}%
	by Fatou Lemma. We have proved that the last norm is bounded by 
	\begin{align*}
	& M_{1}^{1-\theta }M_{2}^{\theta }\Big\|\{t_{k}^{-1}g_{\varpi
		_{k}(n)}^{k}\}|L_{p}(\ell _{q},\{t_{k}\})\Big\| \\
	& \leq M_{1}^{1-\theta }M_{2}^{\theta }\Big\|\{t_{k}^{-1}g_{\varpi
		_{k}(n)}^{k}-f_{k}\}|L_{p}(\ell _{q},\{t_{k}\})\Big\|+M_{1}^{1-\theta
	}M_{2}^{\theta }\Big\|\{f_{k}\}|L_{p}(\ell _{q},\{t_{k}\})\Big\|.
	\end{align*}%
	Consequently%
	\begin{equation*}
	\Big\|\{T(f_{k}\}|L_{p}(\ell _{q},\{t_{k}\})\Big\|\leq M_{1}^{1-\theta
	}M_{2}^{\theta }\Big\|\{f_{k}\}|L_{p}(\ell _{q},\{t_{k}\})\Big\|
	\end{equation*}%
	for any $\{t_{k}f_{k}\}\in L_{p}(\ell _{q})$ and the lemma is proved.
\end{proof}

Now we state the main result of this subsection.

\begin{lem}
	\label{key-estimate1}Let $1<\theta \leq p<\infty $\ and $1<q<\infty $. Let $%
	\{t_{k}\}$\ be a $p$-admissible\ weight\ sequence\ such that $t_{k}^{p}\in
	A_{\frac{p}{\theta }}(\mathbb{R}^{n})$, $k\in \mathbb{Z}$. Assume that $%
	t_{k}^{p}$,\ $k\in \mathbb{Z}$ have the same Muckenhoupt constant, $A_{\frac{%
			p}{\theta }}(t_{k})=c,k\in \mathbb{Z}$. Then%
	\begin{equation}
	\Big\|\Big(\sum\limits_{k=-\infty }^{\infty }t_{k}^{q}\big(\mathcal{M}(f_{k})%
	\big)^{q}\Big)^{1/q}|L_{p}(\mathbb{R}^{n})\Big\|\lesssim \Big\|\Big(%
	\sum\limits_{k=-\infty }^{\infty }t_{k}^{q}\left\vert f_{k}\right\vert ^{q}%
	\Big)^{1/q}|L_{p}(\mathbb{R}^{n})\Big\|  \label{key-est}
	\end{equation}%
	holds for all sequences of functions $\{t_{k}f_{k}\}\in L_{p}(\ell _{q})$.
\end{lem}

\begin{proof}
	We will do the proof into two steps.
	
	\textit{Step 1.} We prove $\mathrm{\eqref{key-est}}$ with $1<q\leq p<\infty $%
	. The case $p=q$ follows, e.g., by \cite{D20}, so we assume that $%
	1<q<p<\infty $. Let $j\in \mathbb{N}$ be such that $q^{j}\leq p<q^{j+1}$.
	
	\textit{Substep 1.1. }We consider the case $q^{j}<p<q^{j+1}$. In order to
	prove we additionally do it into the two steps Substeps 1.1.1 and 1.1.2.
	
	\textit{Substep 1.1.1.} We consider the case $1<q\leq \theta ^{\frac{1}{j}}$%
	. From Lemma \ref{Ap-Property}/(i), $t_{k}^{p}\in A_{p/q^{j}}(\mathbb{R}%
	^{n}) $. By duality the left-hand side of $\mathrm{\eqref{key-est}}$ is
	bounded by 
	\begin{equation*}
	\sup \sum\limits_{k=-\infty }^{\infty }\int_{\mathbb{R}^{n}}t_{k}(x)\mathcal{%
		M}(f_{k})(x)|g_{k}(x)|dx=\sup \sum\limits_{k=-\infty }^{\infty }T_{k},
	\end{equation*}%
	where the supremum is taken over all sequence of functions $\{g_{k}\}\in
	L_{p^{\prime }}(\ell _{q^{\prime }})$ with 
	\begin{equation*}
	\big\|\{g_{k}\}|L_{p^{\prime }}(\ell _{q^{\prime }})\big\|\leq 1,
	\end{equation*}%
	where $p^{\prime }$ and $q^{\prime }$ are the conjugate exponent of $p$ and $%
	q$, respectively. Let $Q$ be a cube. By H\"{o}lder's inequality,%
	\begin{equation*}
	M_{Q}(f_{k})\leq \frac{1}{|Q|}\big\|t_{k}f_{k}|L_{p}(Q)\big\|\big\|%
	t_{k}^{-1}|L_{p^{\prime }}(Q)\big\|\leq \frac{c}{|Q|}\big\|%
	t_{k}^{-1}|L_{p^{\prime }}(Q)\big\|,\quad k\in \mathbb{Z}
	\end{equation*}%
	with $c>0$ is independent of $k$. Since $t_{k}^{p}\in A_{p/\theta }(\mathbb{R%
	}^{n})$, $k\in \mathbb{Z}$, by Lemma \ref{Ap-Property}/(i),\ $t_{k}^{p}\in
	A_{p}(\mathbb{R}^{n})$, $k\in \mathbb{Z}$\ and 
	\begin{equation*}
	\frac{1}{|Q|}\big\|t_{k}^{-1}|L_{p^{\prime }}(Q)\big\|\leq c\big\|%
	t_{k}|L_{p}(Q)\big\|^{-1}.
	\end{equation*}%
	Moreover $\big\|t_{k}|L_{p}(Q)\big\|\rightarrow \infty $ as $|Q|\rightarrow
	\infty $ for any $k\in \mathbb{Z}$. Hence we can apply Lemma \ref{CZ-lemma}.
	Let 
	\begin{equation*}
	\Omega _{k}^{i}=\{x\in \mathbb{R}^{n}:\mathcal{M}(f_{k})(x)>4^{n}\lambda
	^{i}\},\quad k,i\in \mathbb{Z}
	\end{equation*}%
	with $\lambda >2^{n+1}$\ and%
	\begin{equation*}
	H_{k}^{i}=\{x\in \mathbb{R}^{n}:4^{n}\lambda ^{i}<\mathcal{M}(f_{k})(x)\leq
	4^{n}\lambda ^{i+1}\},\quad k,i\in \mathbb{Z}.
	\end{equation*}%
	We have%
	\begin{equation*}
	T_{k}=\sum_{i=-\infty }^{\infty }\int_{H_{k}^{i}}t_{k}(x)\mathcal{M}%
	(f_{k})(x)|g_{k}(x)|dx\leq 4^{n}\sum_{i=-\infty }^{\infty }\lambda
	^{i+1}\int_{\Omega _{k}^{i}}t_{k}(x)|g_{k}(x)|dx.
	\end{equation*}%
	Let\ $\{Q_{k}^{i,h}\}_{h}$ be the collection of maximal dyadic cubes as in
	Lemma \ref{CZ-lemma} with%
	\begin{equation*}
	\Omega _{k}^{i}\subset \cup _{h}3Q_{k}^{i,h},
	\end{equation*}%
	which implies that%
	\begin{equation}
	T_{k}\leq 4^{n}\sum_{i=-\infty }^{\infty }\sum_{h=0}^{\infty }\lambda
	^{i+1}\int_{3Q_{k}^{i,h}}t_{k}(x)|g_{k}(x)|dx,\quad k\in \mathbb{Z}.
	\label{Key-est-Tk}
	\end{equation}%
	Applying H\"{o}lder's inequality,%
	\begin{align*}
	\int_{3Q_{k}^{i,h}}t_{k}(x)|g_{k}(x)|dx &\leq \Big(%
	\int_{3Q_{k}^{i,h}}t_{k}^{\tau }(x)dx\Big)^{1/\tau }\Big(%
	\int_{3Q_{k}^{i,h}}|g_{k}(x)|^{\tau ^{\prime }}dx\Big)^{1/\tau ^{\prime }} \\
	&=|3Q_{k}^{i,h}|M_{3Q_{k}^{i,h},\tau }(t_{k})M_{3Q_{k}^{i,h},\tau ^{\prime
	}}(g_{k})
	\end{align*}%
	with $\tau >1$. Put $\tau =p(1+\varepsilon )$ with $\varepsilon $ as in
	Theorem \ref{reverse Holder inequality}, which is possible since $%
	t_{k}^{p}\in A_{p/\theta }(\mathbb{R}^{n})$, $k\in \mathbb{Z}$. Obviously,
	we have%
	\begin{equation*}
	M_{3Q_{k}^{i,h},\tau }(t_{k})=M_{3Q_{k}^{i,h},p(1+\varepsilon )}(t_{k})\leq c%
	\text{ }M_{3Q_{k}^{i,h},p}(t_{k}),\quad k\in \mathbb{Z}.
	\end{equation*}%
	Since $t_{k}^{p}$,\ $k\in \mathbb{Z}$ have the same Muckenhoupt constant and
	from the proof of Theorem 7.2.2 in \cite{L. Graf14} the constant $c$ is
	independent of $k$. Therefore, 
	\begin{equation*}
	\int_{3Q_{k}^{i,h}}t_{k}(x)|g_{k}(x)|dx\lesssim
	|Q_{k}^{i,h}|M_{3Q_{k}^{i,h},p}(t_{k})M_{3Q_{k}^{i,h},\tau ^{\prime
	}}(g_{k}).
	\end{equation*}%
	We deduce from the above that%
	\begin{equation*}
	\lambda ^{i}\int_{3Q_{k}^{i,h}}t_{k}(x)|g_{k}(x)|dx\lesssim
	|Q_{k}^{i,h}|M_{3Q_{k}^{i,h},p}(t_{k})M_{Q_{k}^{i,h}}(f_{k})M_{3Q_{k}^{i,h},%
		\tau ^{\prime }}(g_{k}).
	\end{equation*}%
	By H\"{o}lder's inequality,%
	\begin{equation*}
	M_{Q_{k}^{i,h}}(f_{k})\leq M_{3Q_{k}^{i,h},s^{\prime
	}}(t_{k}^{-1})M_{3Q_{k}^{i,h},s}(t_{k}f_{k}),\quad s=\frac{p}{q^{j}}
	\end{equation*}%
	and, with the help of the fact that $q^{j}>1$, 
	\begin{equation*}
	M_{3Q_{k}^{i,h},s^{\prime }}(t_{k}^{-1})\leq \big(M_{3Q_{k}^{i,h},\frac{%
			s^{\prime }}{s}}(t_{k}^{-p})\big)^{1/p}.
	\end{equation*}%
	Hence%
	\begin{equation*}
	\lambda ^{i}\int_{3Q_{k}^{i,h}}t_{k}(x)|g_{k}(x)|dx\lesssim
	|Q_{k}^{i,h}|M_{3Q_{k}^{i,h},s}(t_{k}f_{k})M_{3Q_{k}^{i,h},\tau ^{\prime
	}}(g_{k}),
	\end{equation*}%
	because of $t_{k}^{p}$ $ \in $ $A_{s}(\mathbb{R}^{n})$,\ $k\in \mathbb{Z}$  . Since $|Q_{k}^{i,h}|\leq \beta |E_{k}^{i,h}|$, with $E_{k}^{i,h}=Q_{k}^{i,h}%
	\backslash (Q_{k}^{i,h}\cap (\cup _{h}Q_{k}^{i+1,h}))$ and the family $%
	E_{k}^{i,h}$ are pairwise disjoint, the last expression is bounded by%
	\begin{align*}
	&	c\int_{E_{k}^{i,h}}M_{3Q_{k}^{i,h},s}(t_{k}f_{k})M_{3Q_{k}^{i,h},\tau
		^{\prime }}(g_{k})dx\\
	&\lesssim \int_{\mathbb{R}^{n}}\mathcal{M}%
	_{s}(t_{k}f_{k})(x)\mathcal{M}_{\tau ^{\prime }}(g_{k})(x)\chi
	_{E_{k}^{i,h}}(x)dx.
	\end{align*}%
	Therefore, the right-hand side of $\mathrm{\eqref{Key-est-Tk}}$ does not
	exceed 
	\begin{align*}
	&	c\sum_{i=-\infty }^{\infty }\sum_{h=0}^{\infty }\int_{\mathbb{R}^{n}}%
	\mathcal{M}_{s}(t_{k}f_{k})(x)\mathcal{M}_{\tau ^{\prime }}(g_{k})(x)\chi
	_{E_{k}^{i,h}}(x)dx\\
	&\lesssim \int_{\mathbb{R}^{n}}\mathcal{M}_{s}(t_{k}f_{k})(x)\mathcal{M}_{\tau ^{\prime }}(g_{k})(x)dx.
	\end{align*}%
	This implies that%
	\begin{equation*}
	\sum\limits_{k=-\infty }^{\infty }T_{k}\lesssim \int_{\mathbb{R}%
		^{n}}\sum\limits_{k=0}^{\infty }\mathcal{M}_{s}(t_{k}f_{k})(x)\mathcal{M}%
	_{\tau ^{\prime }}(g_{k})(x)dx.
	\end{equation*}%
	By H\"{o}lder's inequality the term inside the integral is bounded by%
	\begin{equation*}
	\Big(\sum\limits_{k=-\infty }^{\infty }\big(\mathcal{M}_{s}\big(t_{k}f_{k}%
	\big)(x)\big)^{q}\Big)^{1/q}\Big(\sum\limits_{k=-\infty }^{\infty }\big(%
	\mathcal{M}_{\tau ^{\prime }}(g_{k})(x)\big)^{q^{\prime }}\Big)^{1/q^{\prime
	}}
	\end{equation*}%
	for any $x\in \mathbb{R}^{n}$. Again by H\"{o}lder's inequality $%
	\sum\limits_{k=-\infty }^{\infty }T_{k}$ can be estimated by 
	\begin{align*}
	& c\Big\|\Big(\sum\limits_{k=-\infty }^{\infty }\big(\mathcal{M}_{s}\big(%
	t_{k}f_{k}\big)\big)^{q}\Big)^{1/q}|L_{p}(\mathbb{R}^{n})\Big\|\Big\|\Big(%
	\sum\limits_{k=-\infty }^{\infty }\big(\mathcal{M}_{\tau ^{\prime }}(g_{k})%
	\big)^{q^{\prime }}\Big)^{1/q^{\prime }}|L_{p^{\prime }}(\mathbb{R}^{n})%
	\Big\| \\
	& \lesssim \Big\|\Big(\sum\limits_{k=-\infty }^{\infty }t_{k}^{q}\left\vert
	f_{k}\right\vert ^{q}\Big)^{1/q}|L_{p}(\mathbb{R}^{n})\Big\|,
	\end{align*}%
	where in the second inequality we used the vector-valued maximal inequality
	of Fefferman and Stein $\mathrm{\eqref{Fe-St71}}$, because of $0<s<q<p$, and
	the fact that 
	\begin{align*}
	\Big\|\Big(\sum\limits_{k=-\infty }^{\infty }\big(\mathcal{M}_{\tau ^{\prime
	}}(g_{k})\big)^{q^{\prime }}\Big)^{1/q^{\prime }}|L_{p^{\prime }}(\mathbb{R}%
	^{n})\Big\| &\lesssim \big\|\{g_{k}\}|L_{p^{\prime }}(\ell _{q^{\prime }})%
	\big\| \\
	&\lesssim 1,
	\end{align*}%
	since $\tau =p(1+\varepsilon )>p>q$\textrm{. }
	
	\textit{Substep 1.1.2.} We consider the case $\theta ^{\frac{1}{j}}<q<p$.
	Let 
	\begin{equation*}
	\frac{1}{q}=\frac{1-\lambda }{p}+\frac{\lambda }{\theta ^{\frac{1}{j}}}%
	,\quad 0<\lambda <1.
	\end{equation*}%
	Applying Lemma \ref{Calderon-Zygmund} we obtain the desired estimate.
	
	\textit{Substep 1.2. }We have proved $\mathrm{\eqref{key-est}}$ in case of $%
	q^{j}<p<q^{j+1}$. Now we study the case $p=q^{j},j\in \mathbb{N}$. Then $%
	j\geq 2$. Since $t_{k}^{p}\in A_{p}(\mathbb{R}^{n}),k\in \mathbb{Z}$, from
	Theorem \ref{reverse Holder inequality} there exists a number $\gamma _{1}>0$
	such%
	\begin{equation*}
	M_{Q,1+\gamma _{1}}(t_{k}^{p})\lesssim M_{Q}(t_{k}^{p}),\quad k\in \mathbb{Z}%
	,
	\end{equation*}%
	where the implicit constant is independent of $k$, see, e.g., \cite[Theorem
	7.2.5 and Corollary 7.2.6.]{L. Graf14}. In addition by Lemma \ref%
	{Ap-Property}/(i), $t_{k}^{-p^{\prime }}\in A_{p^{\prime }}(\mathbb{R}^{n})$
	and again by Theorem \ref{reverse Holder inequality} there exists a number $%
	\gamma _{2}>0$ such%
	\begin{equation*}
	M_{Q,1+\gamma _{2}}(t_{k}^{-p^{\prime }})\lesssim M_{Q}(t_{k}^{-p^{\prime
	}}),\quad k\in \mathbb{Z}.
	\end{equation*}%
	Let 
	\begin{equation*}
	0<\gamma <\min \Big(\gamma _{1},\gamma _{2},\frac{p-p/q}{p/q-1}\Big).
	\end{equation*}%
	Then 
	\begin{equation*}
	M_{Q,1+\gamma }(t_{k}^{p})\big(M_{Q,1+\gamma }(t_{k}^{-p^{\prime }})\big)^{%
		\frac{p}{p^{\prime }}}\lesssim 1,\quad k\in \mathbb{Z},
	\end{equation*}%
	where the implicit constant is independent of $k$. Hence $t_{k}^{p(1+\gamma
		)}\in A_{p}(\mathbb{R}^{n}),k\in \mathbb{Z}$. Therefore $t_{k}^{p}\in
	A_{p_{1}}(\mathbb{R}^{n}),k\in \mathbb{Z}$, where 
	\begin{equation*}
	p_{1}=\frac{p+\gamma }{1+\gamma },
	\end{equation*}%
	see again \cite[Exercise 7.1.3]{L. Graf14}. Observe that 
	\begin{equation*}
	q^{j-1}<p_{1}<p=q^{j}
	\end{equation*}%
	and $t_{k}^{p_{1}}\in A_{p_{1}}(\mathbb{R}^{n}),k\in \mathbb{Z}$, see Lemma %
	\ref{Ap-Property}/(v). Substep 1.1 gives%
	\begin{equation*}
	\Big\|\Big(\sum\limits_{k=-\infty }^{\infty }t_{k}^{q}\big(\mathcal{M}(f_{k})%
	\big)^{q}\Big)^{1/q}|L_{p_{1}}(\mathbb{R}^{n})\Big\|\lesssim \Big\|\Big(%
	\sum\limits_{k=-\infty }^{\infty }t_{k}^{q}\left\vert f_{k}\right\vert ^{q}%
	\Big)^{1/q}|L_{p_{1}}(\mathbb{R}^{n})\Big\|.
	\end{equation*}%
	The same procedure yields that $t_{k}^{p(1+\nu )}\in A_{p}(\mathbb{R}%
	^{n})\subset A_{p(1+\nu )}(\mathbb{R}^{n}),k\in \mathbb{Z}$ with 
	\begin{equation*}
	0<\nu <\min \Big(\gamma _{1},\gamma _{2},q-1\Big)
	\end{equation*}%
	and then $q^{j}=p<p(1+\nu )<q^{j+1}$. Also, Substep 1.1. gives 
	\begin{equation*}
	\Big\|\Big(\sum\limits_{k=-\infty }^{\infty }t_{k}^{q}\big(\mathcal{M}(f_{k})%
	\big)^{q}\Big)^{1/q}|L_{p(1+\nu )}(\mathbb{R}^{n})\Big\|\lesssim \Big\|\Big(%
	\sum\limits_{k=-\infty }^{\infty }t_{k}^{q}\left\vert f_{k}\right\vert ^{q}%
	\Big)^{1/q}|L_{p(1+\nu )}(\mathbb{R}^{n})\Big\|.
	\end{equation*}%
	Again, by Lemma \ref{Calderon-Zygmund} we obtain the desired estimate.
	
	\textit{Step 2.} We shall prove $\mathrm{\eqref{key-est}}$ with $%
	1<p<q<\infty $. Let $1<\varrho <\theta <\infty $. Again, by duality the
	left-hand side of $\mathrm{\eqref{key-est}}$, raised to the\ power $\varrho $
	is just 
	\begin{equation*}
	\sup \sum\limits_{k=-\infty }^{\infty }\int_{\mathbb{R}^{n}}t_{k}^{\varrho
	}(x)\left( \mathcal{M}(f_{k})(x)\right) ^{\varrho }|g_{k}(x)|dx=\sup V,
	\end{equation*}%
	where the supremum is taken over all sequence of functions $\{g_{k}\}$ such that
	\begin{equation}
	\{g_{k}\}\in L_{(p/\varrho )^{\prime }}(\ell _{(q/\varrho )^{\prime }}),\quad \big\|\{g_{k}\}|L_{(p/\varrho )^{\prime }}(\ell _{(q/\varrho )^{\prime }})\big\|\leq 1.  \label{last}
	\end{equation}%
	By Lemma \ref{FS-lemma} and H\"{o}lder's inequality, $V$ is bounded by%
	\begin{align*}
	&	\ c\int_{\mathbb{R}^{n}}\sum\limits_{k=-\infty }^{\infty }\big|f_{k}(x)\big|%
	^{\varrho }\mathcal{M}(t_{k}^{\varrho }g_{k})(x)dx \\
	& \lesssim \Big\|\Big(\sum\limits_{k=-\infty }^{\infty }t_{k}^{q}|f_{k}|^{q}%
	\Big)^{1/q}|L_{p}(\mathbb{R}^{n})\Big\|^{\varrho }\\
	&\times	\Big\|\Big(%
	\sum\limits_{k=-\infty }^{\infty }t_{k}^{-\varrho (q/\varrho )^{\prime }}(%
	\mathcal{M}(t_{k}^{\varrho }g_{k}))^{(q/\varrho )^{\prime }}\Big)%
	^{1/(q/\varrho )^{\prime }}|L_{(p/\varrho )^{\prime }}(\mathbb{R}^{n})\Big\|.
	\end{align*}%
	By Step 1, the second term is bounded because of $(q/\varrho )^{\prime
	}<(p/\varrho )^{\prime }$, 
	\begin{equation*}
	t_{k}^{-\varrho (p/\varrho )^{\prime }}\in A_{(p/\varrho )^{\prime }}(%
	\mathbb{R}^{n}),\quad k\in \mathbb{Z},
	\end{equation*}%
	and from Lemma \ref{Ap-Property}/(iv) there exists a $1<\delta <(p/\varrho
	)^{\prime }<\infty $ such that 
	\begin{equation*}
	t_{k}^{-\varrho (p/\varrho )^{\prime }}\in A_{(p/\varrho )^{\prime }/\delta
	}(\mathbb{R}^{n}),\quad k\in \mathbb{Z}.
	\end{equation*}%
	Thus, the desired estimate follows by Step 1 and $\mathrm{\eqref{last}}$.
	The proof is complete.%\end{proof}
\end{proof}

\begin{rem}
	\label{r-estimates}$\mathrm{(i)}$ We would like to mention that the result
	of this lemma is true if we assume that $t_{k}\in A_{p/\theta }(\mathbb{R}%
	^{n})$,\ $k\in \mathbb{Z}$, $1<p<\infty $ with%
	\begin{equation*}
	A_{p/\theta }(t_{k})\leq c,\quad k\in \mathbb{Z},
	\end{equation*}%
	where $c>0$ independent of $k$.$\newline
	\mathrm{(ii)}$ The proof of Lemma \ref{key-estimate1}\ for $t_{k}^{p}=\omega 
	$, $k\in \mathbb{Z}$ is given in \cite{AJ80} and \cite{Kok78}. $\newline
	\mathrm{(iii)}$ Lemma \ref{key-estimate1} with $t_{k}^{p}=\omega $, $k\in 
	\mathbb{Z}$, see e.g., \cite{Ca02}, can be obtained by using the
	extrapolation theory of J. Garcia-Cuerva and J.L. Rubio de Francia, \cite%
	{GR85}\ or by the theory of vector-valued singular integral with
	operator-valued kernel, see \cite{RRT86}. $\newline
	\mathrm{(iv)}$ To circumvent the drawbacks of dealing with general weights,
	we use different techniques than those using the papers \cite{AJ80}, \cite%
	{GR85}, \cite{Kok78} and \cite{RRT86}. $\newline
	\mathrm{(v)}$ In view of Lemma \ref{Ap-Property}/(iv) we can assume that $%
	t_{k}^{p}\in A_{p}(\mathbb{R}^{n})$,\ $k\in \mathbb{Z}$, $1<p<\infty $ with%
	\begin{equation*}
	A_{p}(t_{k}^{p})\leq c,\quad k\in \mathbb{Z},
	\end{equation*}%
	where $c>0$ independent of $k$.
\end{rem}

We need the following lemma, which is a discrete convolution inequality.

\begin{lem}
	\label{lq-inequality1}\textit{Let }$0<a<1,1\leq p\leq \infty ,1\leq r\leq
	\infty $ \textit{and }$0<q<\infty $\textit{. Let }$\left\{ f_{k}\right\} $%
	\textit{\ and }$\left\{ g_{k}\right\} $ \textit{be two sequences of
		positive\ real\ functions\ and denote}%
	\begin{equation*}
	\delta _{k}=\sum_{j=-\infty }^{k+v}a^{k-j}\big\|g_{k}f_{j}|L_{1}(\mathbb{R}%
	^{n})\big\|^{1/q},\quad k, v\in \mathbb{Z}
	\end{equation*}%
	and\textit{\ }%
	\begin{equation*}
	\eta _{k}=\sum_{j=k+v}^{\infty }a^{j-k}\big\|g_{k}f_{j}|L_{1}(\mathbb{R}^{n})%
	\big\|^{1/q},\quad k, v\in \mathbb{Z}.
	\end{equation*}%
	Then there exists a constant $c>0\ $\textit{depending only on }$a$\textit{\
		and }$q$ such that%
	\begin{equation}
	\sum\limits_{k=-\infty }^{\infty }\delta _{k}^{q}+\sum\limits_{k=-\infty
	}^{\infty }\eta _{k}^{q}\leq c\Big\|\Big(\sum\limits_{k=-\infty }^{\infty
	}f_{k}^{r}\Big)^{1/r}|L_{p}(\mathbb{R}^{n})\Big\|\Big\|\Big(%
	\sum\limits_{k=-\infty }^{\infty }g_{k}^{r^{\prime }}\Big)^{1/r^{\prime
	}}|L_{p^{\prime }}(\mathbb{R}^{n})\Big\|. \label{estimate1}
	\end{equation}
\end{lem}

\begin{proof}
	As the proof for $\{\eta _{k}\}$ is similar, we only consider $\{\delta
	_{k}\}$. We will do the proof in two steps.
	
	\textit{Step 1.} We prove our estimate under the restriction $0<q\leq 1$. We
	have 
	\begin{align*}
	\sum\limits_{k=-\infty }^{\infty }\delta _{k}^{q} &\leq
	\sum\limits_{k=-\infty }^{\infty }\sum_{j=-\infty }^{k+v}a^{(k-j)q}\big\|%
	g_{k}f_{j}|L_{1}(\mathbb{R}^{n})\big\| \\
	&=\sum\limits_{k=-\infty }^{\infty }\sum_{i=-v}^{\infty }a^{iq}\big\|%
	g_{k}f_{k-i}|L_{1}(\mathbb{R}^{n})\big\| \\
	&=\sum\limits_{i=-v}^{\infty }\sum\limits_{k=-\infty }^{\infty }a^{iq}\big\|%
	g_{k}f_{k-i}|L_{1}(\mathbb{R}^{n})\big\|.
	\end{align*}%
	Applying H\"{o}lder's inequality to estimate 
	\begin{equation*}
	\sum\limits_{k=-\infty }^{\infty }\big\|g_{k}f_{k-i}|L_{1}(\mathbb{R}^{n})%
	\big\|
	\end{equation*}%
	by%
	\begin{equation*}
	\Big\|\Big(\sum\limits_{k=-\infty }^{\infty }f_{k}^{r}\Big)^{1/r}|L_{p}(%
	\mathbb{R}^{n})\Big\|\Big\|\Big(\sum\limits_{k=-\infty }^{\infty
	}g_{k}^{r^{\prime }}\Big)^{1/r^{\prime }}|L_{p^{\prime }}(\mathbb{R}^{n})%
	\Big\|.
	\end{equation*}%
	By the fact that $\sum\limits_{i=-v}^{\infty }a^{iq}\lesssim 1$ we obtain the
	desired estimate.
	
	\textit{Step 2.} We consider the case $1<q<\infty $. By duality,%
	\begin{equation*}
	\Big(\sum\limits_{k=-\infty }^{\infty }\delta _{k}^{q}\Big)^{1/q}=\sup
	\sum\limits_{k=-\infty }^{\infty }\sum_{j=-\infty }^{k+v}a^{k-j}\big\|%
	g_{k}f_{j}|L_{1}(\mathbb{R}^{n})\big\|^{1/q}h_{k}=\sup T,
	\end{equation*}%
	where the supremum is taken over all sequence of positive\ real\ numbers $%
	\{h_{k}\}\in \ell _{q^{\prime }}$ with 
	\begin{equation*}
	\Big(\sum\limits_{k=-\infty }^{\infty }h_{k}^{q^{\prime }}\Big)^{1/q^{\prime
	}}\leq 1.
	\end{equation*}%
	Again by H\"{o}lder's inequality,%
	\begin{align*}
	T &=\sum\limits_{k=-\infty }^{\infty }\sum\limits_{i=-v}^{\infty }a^{i}\big\|%
	g_{k}f_{k-i}|L_{1}(\mathbb{R}^{n})\big\|^{1/q}h_{k} \\
	&\leq \sum\limits_{i=-v}^{\infty }a^{i}\Big(\sum\limits_{k=-\infty }^{\infty }%
	\big\|g_{k}f_{k-i}|L_{1}(\mathbb{R}^{n})\big\|\Big)^{1/q}.
	\end{align*}%
	As in Step 1, we obtain the desired estimate. The proof is complete.%
	%\end{proof}
\end{proof}

Using the same type of arguments as in Lemma \ref{lq-inequality1 copy(1)} it
is easy to prove the following lemma.

\begin{lem}
	\label{lq-inequality1 copy(1)}\textit{Let }$a>0,1\leq p\leq \infty ,1\leq
	r\leq \infty $ \textit{and }$0<q<\infty $\textit{. Let }$\left\{
	f_{k}\right\} $\textit{\ and }$\left\{ g_{k}\right\} $ \textit{be two
		sequences of positive\ real\ functions\ and denote}%
	\begin{equation*}
	\delta _{k}=\sum_{j=k}^{k+v}a^{k-j}\big\|g_{k}f_{j}|L_{1}(\mathbb{R}^{n})%
	\big\|^{1/q},\quad k\in \mathbb{Z},v\in \mathbb{N}
	\end{equation*}%
	and\textit{\ }%
	\begin{equation*}
	\eta _{k}=\sum_{j=k+l}^{k}a^{j-k}\big\|g_{k}f_{j}|L_{1}(\mathbb{R}^{n})\big\|%
	^{1/q},\quad k\in \mathbb{Z},l\leq 0.
	\end{equation*}%
	Then there exists a constant $c>0\ $\textit{depending only on }$a,v,l$%
	\textit{\ and }$q$ such that \eqref{estimate1} holds.
\end{lem}

The next lemmas are important for the study of our function spaces.

\begin{lem}
	\label{key-estimate1.1}Let $v\in \mathbb{Z}$, $K\geq 0,1<\theta \leq p<\infty ,1<q<\infty \ $%
	and $\alpha =(\alpha _{1},\alpha _{2})\in \mathbb{R}^{2}$. Let $\{t_{k}\}\in 
	\dot{X}_{\alpha ,\sigma ,p}$ be a $p$-admissible weight sequence with $%
	\sigma =(\sigma _{1}=\theta \left( p/\theta \right) ^{\prime },\sigma
	_{2}\geq p)$. Then for all sequence of functions $\{t_{k}f_{k}\}\in
	L_{p}(\ell _{q})$,%
	\begin{align}
	&\Big\|\Big(\sum\limits_{k=-\infty }^{\infty }t_{k}^{q}\Big(\sum_{j=-\infty
	}^{k+v}2^{(j-k)K}\mathcal{M}(f_{j})\Big)^{q}\Big)^{1/q}|L_{p}(\mathbb{R}^{n})%
	\Big\|\notag \\
	&\lesssim \Big\|\Big(\sum\limits_{k=-\infty }^{\infty
	}t_{k}^{q}\left\vert f_{k}\right\vert ^{q}\Big)^{1/q}|L_{p}(\mathbb{R}^{n})%
	\Big\|  \label{key-est1.1.1}
	\end{align}%
	if $K>\alpha _{2}$ and 
	\begin{align}
	&	\Big\|\Big(\sum\limits_{k=-\infty }^{\infty }t_{k}^{q}\Big(%
	\sum_{j=k+v}^{\infty }2^{(j-k)K}\mathcal{M}(f_{j})\Big)^{q}\Big)^{1/q}|L_{p}(%
	\mathbb{R}^{n})\Big\|\notag \\
	&\lesssim \Big\|\Big(\sum\limits_{k=-\infty }^{\infty
	}t_{k}^{q}\left\vert f_{k}\right\vert ^{q}\Big)^{1/q}|L_{p}(\mathbb{R}^{n})%
	\Big\|  \label{key-est1.1.2.}
	\end{align}%
	if $K<\alpha _{1}$.
\end{lem}

\begin{proof}
	We divide the proof into two steps.
	
	\textit{Step 1.} We will present the proof of\ $\mathrm{\eqref{key-est1.1.1}}
	$\textrm{. }We separate this step into\ two distinct cases: $1<q\leq
	p<\infty $ and $1<p<q<\infty $.
	
	\textit{Substep 1.1.} We consider the case $1<q\leq p<\infty $. By duality
	the left-hand side of $\mathrm{\eqref{key-est1.1.1}}$ is just 
	\begin{equation*}
	\sup \sum\limits_{k=-\infty }^{\infty }\int_{\mathbb{R}^{n}}t_{k}(x)%
	\sum_{j=-\infty }^{k+v}2^{(j-k)K}\mathcal{M}(f_{j})(x)|g_{k}(x)|dx=\sup
	\sum\limits_{k=-\infty }^{\infty }S_{k},
	\end{equation*}%
	where the supremum is taking over all sequence of functions $\{g_{k}\}\in
	L_{p^{\prime }}(\ell _{q^{\prime }})$\ with 
	\begin{equation}
	\big\|\{g_{k}\}|L_{p^{\prime }}(\ell _{q^{\prime }})\big\|\leq 1.
	\label{con-est1}
	\end{equation}%
	We easily find that%
	\begin{equation*}
	S_{k}=\sum_{j=-\infty }^{k+v}2^{(j-k)K}\int_{\mathbb{R}^{n}}\mathcal{M}%
	(f_{j})(x)t_{k}(x)|g_{k}(x)|dx=\sum_{j=-\infty }^{k+v}2^{(j-k)K}D_{k,j}
	\end{equation*}%
	for any $k\in \mathbb{Z}$. As in Lemma \ref{key-estimate1} we find that $%
	M_{Q}(f_{j})\rightarrow 0$ as $|Q|\rightarrow \infty $ for any $j\in \mathbb{%
		Z}$. Therefore,\ we can apply Lemma \ref{CZ-lemma}. Let 
	\begin{equation*}
	\Omega _{j}^{i}=\{x\in \mathbb{R}^{n}:\mathcal{M}(f_{j})(x)>4^{n}\lambda
	^{i}\},\quad j,i\in \mathbb{Z}
	\end{equation*}%
	with $\lambda >2^{n+1}$\ and%
	\begin{equation*}
	H_{j}^{i}=\{x\in \mathbb{R}^{n}:4^{n}\lambda ^{i}<\mathcal{M}(f_{j})(x)\leq
	4^{n}\lambda ^{i+1}\},\quad j,i\in \mathbb{Z}.
	\end{equation*}%
	Let\ $\{Q_{j}^{i,h}\}_{h}$ be the collection of maximal dyadic cubes as in
	Lemma \ref{CZ-lemma} with%
	\begin{equation*}
	\Omega _{j}^{i}\subset \cup _{h}3Q_{j}^{i,h}.
	\end{equation*}%
	We find that%
	\begin{align*}
	D_{k,j}& =\sum_{i=-\infty }^{\infty }\int_{H_{j}^{i}}t_{k}(x)\mathcal{M}%
	(f_{j})(x)|g_{k}(x)|dx \\
	& \lesssim \sum_{i=-\infty }^{\infty }\lambda ^{i}\int_{\Omega
		_{j}^{i}}t_{k}(x)|g_{k}(x)|dx \\
	& \lesssim \sum_{i=-\infty }^{\infty }\sum_{h=0}^{\infty }\lambda
	^{i}\int_{3Q_{j}^{i,h}}t_{k}(x)|g_{k}(x)|dx \\
	& \lesssim \sum_{i=-\infty }^{\infty }\sum_{h=0}^{\infty }\frac{1}{%
		|Q_{j}^{i,h}|}\int_{Q_{j}^{i,h}}|f_{j}(x)|dx%
	\int_{3Q_{j}^{i,h}}t_{k}(x)|g_{k}(x)|dx
	\end{align*}%
	for any $j\leq k+v$. Notice that, by the H\"{o}lder inequality, we find that
	for all $j\leq k+v$, 
	\begin{equation*}
	\frac{1}{|3Q_{j}^{i,h}|}\int_{3Q_{j}^{i,h}}t_{k}(x)|g_{k}(x)|dx\leq
	M_{3Q_{j}^{i,h},\delta ^{\prime }}(g_{k})M_{3Q_{j}^{i,h},\delta }(t_{k}),
	\end{equation*}%
	with $\delta =p(1+\varepsilon )$ and $\varepsilon $ as in Theorem \ref%
	{reverse Holder inequality}, which is possible since $t_{k}^{p}\in
	A_{p/\theta }(\mathbb{R}^{n})$ for any $k\in \mathbb{Z}$. We shall
	distinguish two cases. Let $j\mathbb{\in Z}$ be such that $j\leq k$. We
	obtain by $\mathrm{\eqref{Asum2}}$ 
	\begin{equation*}
	M_{3Q_{j}^{i,h},\delta }(t_{k})\lesssim M_{3Q_{j}^{i,h},p}(t_{k})\lesssim
	2^{\alpha _{2}(k-j)}M_{3Q_{j}^{i,h},p}(t_{j}).
	\end{equation*}%
	By H\"{o}lder's inequality,%
	\begin{equation*}
	1=\Big(\frac{1}{|3Q_{j}^{i,h}|}\int_{3Q_{j}^{i,h}}t_{j}^{-\eta
	}(x)t_{j}^{\eta }(x)dx\Big)^{1/\eta }\leq M_{3Q_{j}^{i,h},\varrho
	}(t_{j}^{-1})M_{3Q_{j}^{i,h},\tau }(t_{j})
	\end{equation*}%
	for any $\eta >0$ and any $\varrho ,\tau >0$ with $1/\eta =1/\varrho +1/\tau 
	$. Taking any $0<\varrho <\sigma _{1}$ and any $0<\tau <q<\infty $ we obtain%
	\begin{equation}
	1\leq M_{3Q_{j}^{i,h},\sigma _{1}}(t_{j}^{-1})M_{3Q_{j}^{i,h},\tau }(t_{j}),
	\label{estimate2}
	\end{equation}%
	which together with $\mathrm{\eqref{Asum1}}$ implies that%
	\begin{equation*}
	M_{3Q_{j}^{i,h},\delta }(t_{k})\lesssim 2^{\alpha
		_{2}(k-j)}M_{3Q_{j}^{i,h},\tau }(t_{j}),
	\end{equation*}%
	which further implies that 
	\begin{align*}
	M_{3Q_{j}^{i,h},\delta }(t_{k})\int_{Q_{j}^{i,h}}|f_{j}(x)|dx& \lesssim
	2^{\alpha _{2}(k-j)}M_{3Q_{j}^{i,h},\tau
	}(t_{j})\int_{3Q_{j}^{i,h}}|f_{j}(x)|dx \\
	& \lesssim 2^{\alpha _{2}(k-j)}|Q_{j}^{i,h}|M_{3Q_{j}^{i,h},\tau }(t_{j}%
	\mathcal{M}(f_{j})).
	\end{align*}%
	Thus,%
	\begin{equation*}
	\frac{1}{|3Q_{j}^{i,h}|}\int_{Q_{j}^{i,h}}|f_{j}(x)|dx%
	\int_{3Q_{j}^{i,h}}t_{k}(x)|g_{k}(x)|dx
	\end{equation*}%
	is bounded by%
	\begin{equation*}
	c\text{ }2^{\alpha _{2}(k-j)}|Q_{j}^{i,h}|M_{3Q_{j}^{i,h},\delta ^{\prime
	}}(g_{k})M_{3Q_{j}^{i,h},\tau }(t_{j}\mathcal{M}(f_{j})),
	\end{equation*}%
	where the positive constant $c$ is independent of $j,h$ and $k$. Let $j%
	\mathbb{\in Z}$ be such that $k<j\leq k+v$. By $\mathrm{\eqref{estimate2}}$
	and $\mathrm{\eqref{Asum1}}$ we obtain%
	\begin{equation*}
	M_{3Q_{j}^{i,h},\delta }(t_{k})\lesssim M_{3Q_{j}^{i,h},p}(t_{k})\lesssim
	2^{\alpha _{1}(k-j)}M_{3Q_{j}^{i,h},p}(t_{j}).
	\end{equation*}%
	Consequently,%
	\begin{equation*}
	D_{k,j}\lesssim \Lambda _{k,j}\sum_{i=-\infty }^{\infty }\sum_{h=0}^{\infty
	}|Q_{j}^{i,h}|M_{3Q_{j}^{i,h},\delta ^{\prime }}(g_{k})M_{3Q_{j}^{i,h},\tau
	}(t_{j}\mathcal{M}(f_{j}))
	\end{equation*}%
	for any $j\leq k+v$, where%
	\begin{equation*}
	\Lambda _{k,j}=\left\{ 
	\begin{array}{ccc}
	2^{\alpha _{2}(k-j)}, & \text{if} & j\leq k, \\ 
	2^{\alpha _{1}(k-j)} & \text{if} & k<j\leq k+v.%
	\end{array}%
	\right. 
	\end{equation*}%
	Since $|Q_{j}^{i,h}|\leq \beta |E_{j}^{i,h}|$, with $E_{j}^{i,h}=Q_{j}^{i,h}%
	\backslash (Q_{j}^{i,h}\cap (\cup _{h}Q_{j}^{i+1,h}))$ and the family $%
	E_{j}^{i,h}$ are pairwise, we find that 
	\begin{align*}
	D_{k,j}& \lesssim \Lambda _{k,j}\sum_{i=-\infty }^{\infty
	}\sum_{h=0}^{\infty }|E_{j}^{i,h}|M_{3Q_{j}^{i,h},\delta ^{\prime
	}}(g_{k})M_{3Q_{j}^{i,h},\tau }(t_{j}\mathcal{M}(f_{j})) \\
	& =c\text{ }\Lambda _{k,j}\sum_{i=-\infty }^{\infty }\sum_{h=0}^{\infty
	}\int_{E_{j}^{i,h}}M_{3Q_{j}^{i,h},\delta ^{\prime
	}}(g_{k})M_{3Q_{j}^{i,h},\tau }(t_{j}\mathcal{M}(f_{j}))dx \\
	& \lesssim \Lambda _{k,j}\sum_{i=-\infty }^{\infty }\sum_{h=0}^{\infty
	}\int_{E_{j}^{i,h}}\mathcal{M}_{\delta ^{\prime }}(g_{k})(x)\mathcal{M}%
	_{\tau }(t_{j}\mathcal{M}(f_{j}))(x)dx \\
	& \lesssim \Lambda _{k,j}\int_{\mathbb{R}^{n}}\mathcal{M}_{\delta ^{\prime
	}}(g_{k})(x)\mathcal{M}_{\tau }(t_{j}\mathcal{M}(f_{j}))(x)dx
	\end{align*}%
	for any $j\leq k+v$. Since $K>\alpha _{2}$, applying Lemmas \ref%
	{lq-inequality1} and \ref{lq-inequality1 copy(1)}, $\sum\limits_{k=-\infty
	}^{\infty }S_{k}$ can be estimated by
	
	\begin{align*}
	& c\Big\|\Big(\sum\limits_{k=-\infty }^{\infty }\big(\mathcal{M}_{\tau
	}(t_{k}\mathcal{M}(f_{k}))\big)^{q}\Big)^{1/q}|L_{p}(\mathbb{R}^{n})\Big\|%
	\Big\|\Big(\sum\limits_{k=-\infty }^{\infty }\big(\mathcal{M}_{\delta
		^{\prime }}(g_{k})\big)^{q^{\prime }}\Big)^{1/q^{\prime }}|L_{p^{\prime }}(%
	\mathbb{R}^{n})\Big\| \\
	& \lesssim \Big\|\Big(\sum\limits_{k=-\infty }^{\infty }t_{k}^{q}|f_{k}|^{q}%
	\Big)^{1/q}|L_{p}(\mathbb{R}^{n})\Big\|,
	\end{align*}%
	where we used the vector-valued maximal inequality of Fefferman and Stein $%
	\mathrm{\eqref{Fe-St71}}$, Lemma \ref{key-estimate1} and $\mathrm{%
		\eqref{con-est1}}$.
	
	\textit{Substep 1.2.} We consider the case $1<p<q<\infty $. Again by duality
	the left-hand side of $\mathrm{\eqref{key-est1.1.1}}$, raised to the power\ $%
	\theta $, is just 
	\begin{equation*}
	\sup \sum\limits_{k=-\infty }^{\infty }\int_{\mathbb{R}^{n}}t_{k}^{\theta
	}(x)\Big(\sum_{j=-\infty }^{k+v}2^{(j-k)K}\mathcal{M(}f_{j})(x)\Big)^{\theta
	}|g_{k}(x)|dx=\sup \sum\limits_{k=-\infty }^{\infty }S_{k},
	\end{equation*}%
	where the supremum is taking over all sequence of functions $\{g_{k}\}$  such that 
	\begin{equation}
	\{g_{k}\}\in
	L_{(p/\theta )^{\prime }}(\ell _{(q/\theta )^{\prime }}), \quad \big\|\{g_{k}\}|L_{(p/\theta )^{\prime }}(\ell _{(q/\theta )^{\prime }})%
	\big\|\leq 1.  \label{con-est2}
	\end{equation}%
	Notice that, by the Minkowski inequality, we obtain that%
	\begin{align*}
	S_{k}& =\Big\|t_{k}\sum_{j=-\infty }^{k+v}2^{(j-k)K}\mathcal{M}%
	(f_{j})|g_{k}|^{1/\theta }|L_{\theta }(\mathbb{R}^{n})\Big\|^{\theta } \\
	& \leq \Big(\sum_{j=-\infty }^{k+v}2^{(j-k)K}\Big(\int_{\mathbb{R}%
		^{n}}t_{k}^{\theta }(x)\big(\mathcal{M}(f_{j})(x)\big)^{\theta }|g_{k}(x)|dx%
	\Big)^{1/\theta }\Big)^{\theta }
	\end{align*}%
	for any $k\in \mathbb{Z}$. Using Lemma\ \ref{FS-lemma}, we deduce that, for
	any $j\leq k+v$, 
	\begin{equation*}
	\int_{\mathbb{R}^{n}}t_{k}^{\theta }(x)\big(\mathcal{M}(f_{j})(x)\big)%
	^{\theta }|g_{k}(x)|dx\lesssim \int_{\mathbb{R}^{n}}|f_{j}(x)|^{\theta }%
	\mathcal{M}(t_{k}^{\theta }g_{k})(x)dx=c\text{ }D_{k,j},
	\end{equation*}%
	where the positive constant $c$ is independent of $k$ and $j$. Let $Q$ be a
	cube. By H\"{o}lder's inequality,%
	\begin{equation*}
	\frac{1}{|Q|}\int_{Q}t_{k}^{\theta }(x)|g_{k}(x)|dx\leq \frac{1}{|Q|}\big\|%
	g_{k}|L_{\left( p/\theta \right) ^{\prime }}(Q)\big\|\big\|t_{k}^{\theta
	}|L_{p/\theta }(Q)\big\|\leq \frac{1}{|Q|}\big\|t_{k}|L_{p}(Q)\big\|^{\theta
	},
	\end{equation*}%
	with $c>0$ is independent of $k$. Since $\{t_{k}\}$ is a $p$-admissible
	sequence satisfying $\mathrm{\eqref{Asum1}}$ with $\sigma _{1}=\theta \left(
	p/\theta \right) ^{\prime }$, we find that 
	\begin{equation*}
	\frac{1}{|Q|^{\frac{1}{\theta }}}\big\|t_{k}|L_{p}(Q)\big\|\leq
	C|Q|^{1/p+1/\sigma _{1}-1/\theta }\big\|t_{k}^{-1}|L_{\sigma _{1}}(Q)\big\|%
	^{-1}=C\big\|t_{k}^{-1}|L_{\sigma _{1}}(Q)\big\|^{-1},
	\end{equation*}%
	with $C>0$ is independent of $Q$, $\theta $ and $k$. Since $t_{k}^{p}\in
	A_{p/\theta }(\mathbb{R}^{n})$ it follows by Lemma \ref{Ap-Property}/(ii), 
	\begin{equation*}
	t_{k}^{-\sigma _{1}}\in A_{(p/\theta )^{\prime }}(\mathbb{R}^{n}),\quad k\in 
	\mathbb{Z}.
	\end{equation*}%
	Hence $\big\|t_{k}^{-1}|L_{\sigma _{1}}(Q)\big\|^{-1}\rightarrow 0$, $%
	|Q|\rightarrow \infty $, $k\in \mathbb{Z}$. Therefore,\ we can apply Lemma %
	\ref{CZ-lemma}. Using the same arguments as in proof of Lemma \ref%
	{key-estimate1}, we get%
	\begin{align*}
	D_{k,j}& \lesssim \sum_{i=-\infty }^{\infty }\sum_{h=0}^{\infty }\lambda
	^{i}\int_{3Q_{k}^{i,h}}|f_{j}(x)|^{\theta }dx \\
	& \lesssim \sum_{i=-\infty }^{\infty }\sum_{h=0}^{\infty }\frac{1}{%
		|Q_{k}^{i,h}|}\int_{3Q_{k}^{i,h}}|f_{j}(x)|^{\theta
	}dx\int_{Q_{k}^{i,h}}t_{k}^{\theta }(x)|g_{k}(x)|dx,
	\end{align*}%
	where we have used the inequality 
	\begin{equation*}
	\lambda ^{i}\leq \frac{1}{|Q_{k}^{i,h}|}\int_{Q_{k}^{i,h}}t_{k}^{\theta
	}(x)|g_{k}(x)|dx\leq 2^{n}\lambda ^{i}.
	\end{equation*}%
	On the other hand, by H\"{o}lder's inequality, we obtain%
	\begin{equation*}
	\int_{3Q_{k}^{i,h}}|f_{j}(x)|^{\theta }dx\leq
	|Q_{k}^{i,h}|M_{3Q_{k}^{i,h},\nu ^{\prime }}(t_{j}^{-\theta
	})M_{3Q_{k}^{i,h},\nu }(t_{j}^{\theta }|f_{j}|^{\theta })
	\end{equation*}%
	with $\nu ^{\prime }=\left( p/\theta \right) ^{\prime }(1+\varepsilon )$ and 
	$\varepsilon $ as in Theorem\ \ref{reverse Holder inequality} since, again, $%
	t_{k}^{-\sigma _{1}}\in A_{(p/\theta )^{\prime }}(\mathbb{R}^{n})$ for any $%
	k\in \mathbb{Z}$. Therefore, 
	\begin{equation*}
	M_{3Q_{k}^{i,h},\nu ^{\prime }}(t_{j}^{-\theta })=\big(M_{3Q_{k}^{i,h},%
		\sigma _{1}(1+\varepsilon )}(t_{j}^{-1})\big)^{\theta }\lesssim \big(%
	M_{3Q_{k}^{i,h},\sigma _{1}}(t_{j}^{-1})\big)^{\theta }.
	\end{equation*}%
	As before, by H\"{o}lder's inequality,%
	\begin{equation}
	1=\Big(\frac{1}{|3Q_{k}^{i,h}|}\int_{3Q_{k}^{i,h}}t_{k}^{-\varrho
	}(x)t_{k}^{\varrho }(x)dx\Big)^{1/\varrho }\leq M_{3Q_{k}^{i,h},\theta
	}(t_{k}^{-1})M_{3Q_{k}^{i,h},p}(t_{k})  \label{estimate6}
	\end{equation}%
	for any $\varrho >0$ with $1/\varrho =1/\theta +1/p$. Again, we distinguish
	two cases. Let $j\mathbb{\in Z}$ be such that $j\leq k$. From $\mathrm{%
		\eqref{Asum2}}$ we obtain%
	\begin{equation*}
	M_{3Q_{k}^{i,h},p}(t_{k})\lesssim 2^{\alpha
		_{2}(k-j)}M_{3Q_{k}^{i,h},p}(t_{j}),\quad j\leq k.
	\end{equation*}%
	Therefore,%
	\begin{equation*}
	1\lesssim 2^{\alpha _{2}(k-j)}M_{3Q_{k}^{i,h},\theta
	}(t_{k}^{-1})M_{3Q_{k}^{i,h},p}(t_{j}),\quad j\leq k.
	\end{equation*}%
	Multiplying by $M_{3Q_{k}^{i,h},\sigma _{1}}(t_{j}^{-1})$ and using $\mathrm{%
		\eqref{Asum1}}$ we get%
	\begin{equation*}
	M_{3Q_{k}^{i,h},\sigma _{1}}(t_{j}^{-1})\lesssim 2^{\alpha
		_{2}(k-j)}M_{3Q_{k}^{i,h},\theta }(t_{k}^{-1}).
	\end{equation*}%
	Now, Let $j\mathbb{\in Z}$ be such that $k<j\leq k+v$. From $\mathrm{%
		\eqref{estimate6}}$ and $\mathrm{\eqref{Asum1}}$ we obtain%
	\begin{equation*}
	M_{3Q_{k}^{i,h},\sigma _{1}}(t_{j}^{-1})\lesssim 2^{\alpha
		_{1}(k-j)}M_{3Q_{k}^{i,h},\theta }(t_{k}^{-1}).
	\end{equation*}%
	Hence%
	\begin{equation*}
	|Q_{k}^{i,h}|M_{3Q_{k}^{i,h}}(|f_{j}|^{\theta
	})M_{Q_{k}^{i,h}}(t_{k}^{\theta }g_{k}),
	\end{equation*}%
	can be estimated by%
	\begin{align*}
	& c\Lambda _{k,j}^{\theta }|Q_{k}^{i,h}|M_{Q_{k}^{i,h}}(t_{k}^{\theta
	}g_{k})M_{3Q_{k}^{i,h},\nu }(t_{j}^{\theta }|f_{j}|^{\theta
	})M_{3Q_{k}^{i,h}}(t_{k}^{-\theta }) \\
	& \lesssim \Lambda _{k,j}^{\theta
	}|Q_{k}^{i,h}|M_{3Q_{k}^{i,h}}(t_{k}^{-\theta }\mathcal{M}(t_{k}^{\theta
	}g_{k}))M_{3Q_{k}^{i,h},\nu }(t_{j}^{\theta }|f_{j}|^{\theta }).
	\end{align*}%
	Consequently,%
	\begin{equation*}
	D_{k,j}\lesssim \Lambda _{k,j}^{\theta }\sum_{i=-\infty }^{\infty
	}\sum_{h=0}^{\infty }M_{3Q_{k}^{i,h},\nu }(t_{j}^{\theta }|f_{j}|^{\theta
	})\int_{3Q_{k}^{i,h}}t_{k}^{-\theta }(x)\mathcal{M}(t_{k}^{\theta
	}g_{k})(x)dx
	\end{equation*}%
	for any $j\leq k+v$. Since $|Q_{k}^{i,h}|\leq \beta |E_{k}^{i,h}|$, with $%
	E_{k}^{i,h}=Q_{k}^{i,h}\backslash (Q_{k}^{i,h}\cap (\cup _{h}Q_{k}^{i+1,h}))$
	and of the family $E_{k}^{i,h}$ are pairwise, as before we find that 
	\begin{equation*}
	D_{k,j}\lesssim \Lambda _{k,j}^{\theta }\int_{\mathbb{R}^{n}}\mathcal{M}%
	(t_{k}^{-\theta }\mathcal{M}(t_{k}^{\theta }g_{k}))(x)\mathcal{M}_{\nu
	}(t_{j}^{\theta }|f_{j}|^{\theta })(x)dx.
	\end{equation*}%
	Since $K>\alpha _{2}$ and $\nu <p/\theta \leq q/\theta $, applying Lemma \ref%
	{lq-inequality1}, we obtain $\sum\limits_{k=-\infty }^{\infty }S_{k}$ can be
	estimated by 
	\begin{align*}
	& c\Big\|\Big(\sum\limits_{k=-\infty }^{\infty }\big(\mathcal{M}_{\nu
	}(t_{k}^{\theta }|f_{k}|^{\theta })\big)^{q/\theta }\Big)^{\theta
		/q}|L_{p/\theta }(\mathbb{R}^{n})\Big\| \\
	& \times \Big\|\Big(\sum\limits_{k=-\infty }^{\infty }\big(\mathcal{M}%
	(t_{k}^{-\theta }\mathcal{M}(t_{k}^{\theta }g_{k}))\big)^{(q/\theta
		)^{\prime }}\Big)^{1/(q/\theta )^{\prime }}|L_{(p/\theta )^{\prime
	}}(\mathbb{R}^{n})\Big\| \\
	& \lesssim \Big\|\Big(\sum\limits_{k=-\infty }^{\infty }t_{k}^{q}|f_{k}|^{q}%
	\Big)^{1/q}|L_{p}(\mathbb{R}^{n})\Big\|^{\theta } \\
	& \times \Big\|\Big(\sum\limits_{k=-\infty }^{\infty }\big(t_{k}^{-\theta }%
	\mathcal{M}(t_{k}^{\theta }g_{k})\big)^{(q/\theta )^{\prime }}\Big)^{1/(q/\theta )^{\prime }}|L_{(p/\theta )^{\prime }}(\mathbb{R}^{n})\Big\| \\
	& \lesssim \Big\|\Big(\sum\limits_{k=-\infty }^{\infty }t_{k}^{q}|f_{k}|^{q}%
	\Big)^{1/q}|L_{p}(\mathbb{R}^{n})\Big\|^{\theta },
	\end{align*}%
	where we used the vector-valued maximal inequality of Fefferman and Stein $%
	\mathrm{\eqref{Fe-St71}}$, Lemma \ref{key-estimate1} and $\mathrm{%
		\eqref{con-est2}}$.
	
	This proves the first part of the Lemma.
	
	\textit{Step 2.} We will present the proof of $\mathrm{\eqref{key-est1.1.2.}}
	$. In order to prove we additionally do it into the two steps Substeps 2.1
	and 2.2.
	
	\textit{Substep 2.1.} We consider the case $1<q\leq p<\infty $. We employ
	the same notation as in Substep 1.1. We obtain 
	\begin{align*}
	D_{k,j}& \lesssim \sum_{i=-\infty }^{\infty }\sum_{h=0}^{\infty }\lambda
	^{i}\int_{3Q_{j}^{i,h}}t_{k}(x)|g_{k}(x)|dx \\
	& \lesssim \sum_{i=-\infty }^{\infty }\sum_{h=0}^{\infty }\frac{1}{%
		|Q_{j}^{i,h}|}\int_{Q_{j}^{i,h}}|f_{j}(x)|dx%
	\int_{3Q_{j}^{i,h}}t_{k}(x)|g_{k}(x)|dx.
	\end{align*}%
	Let $j\geq k+v$. Recall that%
	\begin{equation*}
	M_{3Q_{j}^{i,h},\delta }(t_{k})\leq M_{3Q_{j}^{i,h},p}(t_{k})\quad \text{and}%
	\quad 1\leq M_{3Q_{j}^{i,h},\sigma _{1}}(t_{j}^{-1})M_{3Q_{j}^{i,h},\tau
	}(t_{j}).
	\end{equation*}%
	Using $\mathrm{\eqref{Asum1}}$, we find%
	\begin{equation*}
	M_{3Q_{j}^{i,h},\delta }(t_{k})\lesssim 2^{\alpha
		_{1}(k-j)}M_{3Q_{j}^{i,h},p}(t_{j}),\quad j\geq k.
	\end{equation*}%
	Now, assume that $k+v\leq j<k$. From $\mathrm{\eqref{Asum2}}$, we get 
	\begin{equation*}
	M_{3Q_{j}^{i,h},p}(t_{k})\lesssim 2^{\alpha
		_{2}(k-j)}M_{3Q_{j}^{i,h},p}(t_{j}).
	\end{equation*}%
	Hence%
	\begin{align*}
	M_{3Q_{j}^{i,h},\delta }(t_{k})\int_{Q_{j}^{i,h}}|f_{j}(x)|dx& \lesssim
	\digamma _{k,j}M_{3Q_{j}^{i,h},\tau }(t_{j})\int_{Q_{j}^{i,h}}|f_{j}(x)|dx \\
	& \lesssim \digamma _{k,j}|Q_{j}^{i,h}|M_{3Q_{j}^{i,h},\tau }(t_{j}\mathcal{M%
	}(f_{j})),
	\end{align*}%
	where%
	\begin{equation*}
	\digamma _{k,j}=\left\{ 
	\begin{array}{ccc}
	2^{\alpha _{1}(k-j)}, & \text{if} & j\geq k, \\ 
	2^{\alpha _{2}(k-j)} & \text{if} & k+v\leq j<k.%
	\end{array}%
	\right. 
	\end{equation*}%
	Repeating the same arguments of Substep 1.1, we obtain the desired estimate.
	
	\textit{Substep 2.2.} We show $\mathrm{\eqref{key-est1.1.2.}}$ under the
	assumption $1<p<q<\infty $. We employ the same notation as in Substep 1.2.
	From $\mathrm{\eqref{estimate6}}$ we get 
	\begin{equation*}
	M_{3Q_{k}^{i,h},\sigma _{1}}(t_{j}^{-1})\lesssim 2^{\alpha
		_{1}(k-j)}M_{3Q_{j}^{i,h},\theta }(t_{k}^{-1}),\quad j\geq k.
	\end{equation*}%
	We omit the proof since he is essentially similar to the Substep 2.1 and
	Substep 1.2, respectively.
	
	The proof of lemma is complete.
\end{proof}

\begin{rem}
	\label{r-estimates copy(1)}Let $i\in \mathbb{Z},1<\theta \leq p<\infty ,1<q<\infty \ 
	$and $\alpha =(\alpha _{1},\alpha _{2})\in \mathbb{R}^{2}$. Let $%
	\{t_{k}\}\in \dot{X}_{\alpha ,\sigma ,p}$ be a $p$-admissible weight
	sequence with $\sigma =(\sigma _{1}=\theta \left( p/\theta \right) ^{\prime
	},\sigma _{2}\geq p)$. From Lemma \ref{key-estimate1.1} we easily obtain 
	\begin{equation*}
	\Big\|\Big(\sum\limits_{k=-\infty }^{\infty }t_{k}^{q}\big(\mathcal{M}%
	(f_{k+i})\big)^{q}\Big)^{1/q}|L_{p}(\mathbb{R}^{n})\Big\|\lesssim \Big\|\Big(%
	\sum\limits_{k=-\infty }^{\infty }t_{k}^{q}\left\vert f_{k}\right\vert ^{q}%
	\Big)^{1/q}|L_{p}(\mathbb{R}^{n})\Big\|
	\end{equation*}%
	holds for all sequence of functions $\{t_{k}f_{k}\}\in L_{p}(\ell _{q})$,
	where the implicit constant depends on $i$. Indeed, we have%
	\begin{equation*}
	\mathcal{M}(f_{k+i})\leq \sum\limits_{j=-\infty }^{k+i}2^{(j-k-i)M}\mathcal{M%
	}(f_{j}),\quad M>\alpha _{2},k\in \mathbb{Z}.
	\end{equation*}%
	Lemma \ref{key-estimate1.1} yields the desired result.
\end{rem}

\section{The spaces\textbf{\ }$\dot{F}_{p,q}(\mathbb{R}^{n},\{t_{k}\})$}

In this section we\ present the Fourier analytical definition of
Triebel-Lizorkin spaces of variable smoothness and we prove their basic
properties in analogy to the classical Triebel-Lizorkin spaces.

\subsection{ The $\varphi $-transform characterization }

Select a pair of Schwartz functions $\varphi $ and $\psi $ satisfy%
\begin{equation}
\text{supp}(\mathcal{F}(\varphi ))\cup \text{supp}(\mathcal{F}(\psi
))\subset \big\{\xi :1/2\leq |\xi |\leq 2\big\},  \label{Ass1}
\end{equation}%
\begin{equation}
|\mathcal{F}(\varphi )(\xi )|,|\mathcal{F}(\psi )(\xi )|\geq c\quad \text{if}%
\quad 3/5\leq |\xi |\leq 5/3  \label{Ass2}
\end{equation}%
and 
\begin{equation}
\sum_{k=-\infty }^{\infty }\overline{\mathcal{F}(\varphi )(2^{-k}\xi )}%
\mathcal{F}(\psi )(2^{-k}\xi )=1\quad \text{if}\quad \xi \neq 0,
\label{Ass3}
\end{equation}%
where $c>0$. Throughout the paper we put $\tilde{\varphi}(x)=\overline{%
	\varphi (-x)},x\in \mathbb{R}^{n}$. Let $\varphi \in \mathcal{S}(\mathbb{R}%
^{n})$ be a function satisfying $\mathrm{\eqref{Ass1}}$-$\mathrm{\eqref{Ass2}%
}$. We recall that there exists a function $\psi \in \mathcal{S}(\mathbb{R}%
^{n})$ satisfying $\mathrm{\eqref{Ass1}}$-$\mathrm{\eqref{Ass3}}$, see \cite[%
Lemma (6.9)]{FrJaWe01}.

\begin{rem}
	Let $\dot{F}_{p,q}(\mathbb{R}^{n},\{t_{k}\})$ the spaces under
	consideration. We would like to mention that the elements of the spaces $%
	\dot{F}_{p,q}(\mathbb{R}^{n},\{t_{k}\})$ are not distributions but
	equivalence classes of distributions. Observe that $\dot{F}_{p,p}(\mathbb{R}%
	^{n},\{t_{k}\})$ is just the space $\dot{B}_{p,p}(\mathbb{R}^{n},\{t_{k}\})$%
	, where the space $\dot{B}_{p,q}(\mathbb{R}^{n},\{t_{k}\}),0<p,q\leq \infty $%
	, is defined to be the set of all $f\in \mathcal{S}_{\infty }^{\prime }(%
	\mathbb{R}^{n})$\ such that 
	\begin{equation*}
	\big\|f|\dot{B}_{p,q}(\mathbb{R}^{n},\{t_{k}\})\big\|=\Big(%
	\sum\limits_{k=-\infty }^{\infty }\big\|t_{k}(\varphi _{k}\ast f)|L_{p}(%
	\mathbb{R}^{n})\big\|^{q}\Big)^{1/q}<\infty ,
	\end{equation*}%
	which studied in detail in \cite{D20}.
\end{rem}

Using the system $\{\varphi _{k}\}$ we can define the quasi-norms%
\begin{equation*}
\big\|f|\dot{F}_{p,q}^{s}(\mathbb{R}^{n})\big\|=\big\|\Big(%
\sum\limits_{k=-\infty }^{\infty }2^{ksq}|\varphi _{k}\ast f|^{q}\Big)%
^{1/q}|L_{p}(\mathbb{R}^{n})\big\|
\end{equation*}%
for constants $s\in \mathbb{R}$, $0<p<\infty $ and $0<q\leq \infty $. The
Triebel-Lizorkin space $\dot{F}_{p,q}^{s}(\mathbb{R}^{n})$\ consist of all
distributions $f\in \mathcal{S}_{\infty }^{\prime }(\mathbb{R}^{n})$ for
which 
\begin{equation*}
\big\|f|\dot{F}_{p,q}^{s}(\mathbb{R}^{n})\big\|<\infty .
\end{equation*}%
It is well-known that these spaces do not depend on the choice of the system 
$\{\varphi _{k}\}$ (up to equivalence of quasi-norms). Further details on
the classical theory of these spaces, included the nonhomogeneous case, can
be found in \cite{FJ86}, \cite{FJ90}, \cite{FrJaWe01}, \cite{SchTr87}, \cite%
{T1}, \cite{T2} and \cite{T3}.

One recognizes immediately that if $\{t_{k}\}=\{2^{sk}\}$, $s\in \mathbb{R}$%
, then 
\begin{equation}
\dot{F}_{p,q}(\mathbb{R}^{n},\{2^{sk}\})=\dot{F}_{p,q}^{s}(\mathbb{R}%
^{n})\notag . 
\end{equation}
Moreover, for $\{t_{k}\}=\{2^{sk}w\}$, $s\in \mathbb{R}$ with a
weight $w$ we re-obtain the weighted Triebel-Lizorkin spaces; we refer, in
particular, to the papers \cite{Bui82}, \cite{IzSa12}, \cite{Ry01}, \cite%
{Sch981} and \cite{Sch982} for a comprehensive treatment of the weighted
spaces.

A basic tool to study\ the above\ function\ spaces is the following Calder%
\'{o}n reproducing formula, see \cite[Lemma 2.1]{YY2}.

\begin{lem}
	\label{DW-lemma1}Suppose that $\varphi $, $\psi \in \mathcal{S}(\mathbb{R}%
	^{n})$\ satisfy $\mathrm{\eqref{Ass1}}$ through $\mathrm{\eqref{Ass3}}$%
	\textrm{. }If\textrm{\ }$f\in \mathcal{S}_{\infty }^{\prime }(\mathbb{R}%
	^{n}) $, then%
	\begin{equation}
	f=\sum_{k=-\infty }^{\infty }2^{-kn}\sum_{m\in \mathbb{Z}^{n}}\widetilde{%
		\varphi }_{k}\ast f(2^{-k}m)\psi _{k}(\cdot -2^{-k}m).  \label{proc2}
	\end{equation}
\end{lem}

Let $\varphi $, $\psi \in \mathcal{S}(\mathbb{R}^{n})$ satisfying $\mathrm{%
	\eqref{Ass1}}$\ through\ $\mathrm{\eqref{Ass3}}$. Recall that the $\varphi $%
-transform $S_{\varphi }$ is defined by setting 
\begin{equation*}
(S_{\varphi}f)_{k,m}=\langle f,\varphi _{k,m}\rangle,
\end{equation*}
where $\varphi _{k,m}(x)=2^{k%
	\frac{n}{2}}\varphi (2^{k}x-m)$, $m\in \mathbb{Z}^{n}$ and $k\in \mathbb{Z}$%
. The inverse $\varphi $-transform $T_{\psi }$ is defined by 
\begin{equation*}
T_{\psi }\lambda =\sum_{k=-\infty }^{\infty }\sum_{m\in \mathbb{Z}%
	^{n}}\lambda _{k,m}\psi _{k,m},  
\end{equation*}%
where $\lambda =\{\lambda _{k,m}\}_{k\in \mathbb{Z},m\in \mathbb{Z}%
	^{n}}\subset \mathbb{C}$, see \cite{FJ90}.

Now we introduce the corresponding sequence spaces of $\dot{F}_{p,q}(\mathbb{%
	R}^{n},\{t_{k}\})$.

\begin{defn}
	\label{sequence-space}Let $0<p<\infty $ and $0<q\leq \infty $. Let $%
	\{t_{k}\} $ be a $p$-admissible weight sequence. Then for all complex valued
	sequences $\lambda =\{\lambda _{k,m}\}_{k\in \mathbb{Z},m\in \mathbb{Z}%
		^{n}}\subset \mathbb{C}$ we define%
	\begin{equation*}
	\dot{f}_{p,q}(\mathbb{R}^{n},\{t_{k}\})=\Big\{\lambda :\big\|\lambda |\dot{f}%
	_{p,q}(\mathbb{R}^{n},\{t_{k}\})\big\|<\infty \Big\},
	\end{equation*}%
	where%
	\begin{equation*}
	\big\|\lambda |\dot{f}_{p,q}(\mathbb{R}^{n},\{t_{k}\})\big\|=\Big\|\Big(%
	\sum_{k=-\infty }^{\infty }\sum\limits_{m\in \mathbb{Z}%
		^{n}}2^{knq/2}t_{k}^{q}|\lambda _{k,m}|^{q}\chi _{k,m}\Big)^{1/q}|L_{p}(%
	\mathbb{R}^{n})\Big\|,
	\end{equation*}%
	with the usual modifications if $q=\infty $.
\end{defn}

Allowing the smoothness $t_{k}$, $k\in \mathbb{Z}$ to vary from point to
point will raise extra difficulties\ to study these function spaces. But by
the following lemma the problem can be reduced to the case of fixed
smoothness, see \cite{D20.2}.

\begin{prop}
	\label{Equi-norm1}Let $0<\theta \leq p<\infty $, $0<q<\infty $ and $0<\delta
	\leq 1$. Assume that $\{t_{k}\}$ satisfying $\mathrm{\eqref{Asum1}}$ with $%
	\sigma _{1}=\theta \left( \frac{p}{\theta }\right) ^{\prime }$ and $j=k$.
	Then%
	\begin{align*}
	&	\big\|\lambda |\dot{f}_{p,q,\delta }(\mathbb{R}^{n},\{t_{k}\})\big\|^{\ast }\\
	&	=\Big\|\Big(\sum_{k=-\infty }^{\infty }\sum\limits_{m\in \mathbb{Z}%
		^{n}}2^{knq(1/2+1/\delta p)}t_{k,m,\delta }^{q}|\lambda _{k,m}|^{q}\chi
	_{k,m}\Big)^{1/q}|L_{p}(\mathbb{R}^{n})\Big\|,
	\end{align*}%
	is an equivalent quasi-norm in $\dot{f}_{p,q}(\mathbb{R}^{n},\{t_{k}\})$,
	where%
	\begin{equation*}
	t_{k,m,\delta }=\big\|t_{k}|L_{\delta p}(Q_{k,m})\big\|,\quad k\in \mathbb{Z}%
	,m\in \mathbb{Z}^{n}.
	\end{equation*}
\end{prop}

The following important properties of the sequence spaces will be required
in what follows.

\begin{lem}
	\label{Lamda-est}Let $0<\theta \leq p<\infty $ and $0<q<\infty $. Let $%
	\{t_{k}\}$ be a $p$-admissible weight sequence satisfying $\mathrm{%
		\eqref{Asum1}}$ with $\sigma _{1}=\theta \left( p/\theta \right) ^{\prime }$
	and $j=k$. Let $k\in \mathbb{Z},m\in \mathbb{Z}^{n}$ and $\lambda \in \dot{f}%
	_{p,q}(\mathbb{R}^{n},\{t_{k}\})$.\ Then there exists $c>0$ independent of $%
	k $ and $m$ such that%
	\begin{equation*}
	|\lambda _{k,m}|\leq c\text{ }2^{-kn/2}t_{k,m}^{-1}\big\|\lambda |\dot{f}%
	_{p,q}(\mathbb{R}^{n},\{t_{k}\})\big\|.
	\end{equation*}
\end{lem}

\begin{proof}
	Let $\lambda \in \dot{f}_{p,q}(\mathbb{R}^{n},\{t_{k}\}),k\in \mathbb{Z}\ $%
	and $m\in \mathbb{Z}^{n}$. Since $\{t_{k}\}$ is a $p$-admissible sequence
	satisfying $\mathrm{\eqref{Asum1}}$ with $\sigma _{1}=\theta \left( p/\theta
	\right) ^{\prime }$, we get by H\"{o}lder's inequality 
	\begin{align*}
	|\lambda _{k,m}|& =\Big(\frac{1}{|Q_{k,m}|}\int_{Q_{k,m}}|\lambda
	_{k,m}|^{\theta }dy\Big)^{1/\theta } \\
	& \leq M_{Q_{k,m},p}(\lambda _{k,m}t_{k})M_{Q_{k,m},\sigma _{1}}(t_{k}^{-1})
	\\
	& \leq c\text{ }2^{-\frac{kn}{2}}t_{k,m}^{-1}\big\|\lambda |\dot{f}_{p,q}(%
	\mathbb{R}^{n},\{t_{k}\})\big\|,
	\end{align*}%
	where $c>0$ is independent of $k\in \mathbb{Z}\ $and $m\in \mathbb{Z}^{n}$.
\end{proof}

The following lemma is a slight variant of \cite{D20}. For the convenience
of the reader, we give some details.

\begin{lem}
	\label{Inv-phi-trans}Let $\alpha =(\alpha _{1},\alpha _{2})\in \mathbb{R}%
	^{2},0<\theta \leq p<\infty $ and $0<q\leq \infty $. Let $\{t_{k}\}\in \dot{X%
	}_{\alpha ,\sigma ,p}$ be a $p$-admissible weight sequence with $\sigma
	=(\sigma _{1}=\theta \left( p/\theta \right) ^{\prime },\sigma _{2}\geq p)$.
	Let\ $\psi \in \mathcal{S}(\mathbb{R}^{n})$ satisfying $\mathrm{\eqref{Ass1}}
	$ and $\mathrm{\eqref{Ass2}}$. Then for all $\lambda \in \dot{f}_{p,q}(%
	\mathbb{R}^{n},\{t_{k}\})$%
	\begin{equation*}
	T_{\psi }\lambda =\sum_{k=-\infty }^{\infty }\sum_{m\in \mathbb{Z}%
		^{n}}\lambda _{k,m}\psi _{k,m},
	\end{equation*}%
	converges in $\mathcal{S}_{\infty }^{\prime }(\mathbb{R}^{n})$; moreover, $%
	T_{\psi }:\dot{f}_{p,q}(\mathbb{R}^{n},\{t_{k}\})\rightarrow \mathcal{S}%
	_{\infty }^{\prime }(\mathbb{R}^{n})$ is continuous.
\end{lem}

\begin{proof}
	Let $\lambda \in \dot{f}_{p,q}(\mathbb{R}^{n},\{t_{k}\})$ and $\varphi \in 
	\mathcal{S}_{\infty }(\mathbb{R}^{n})$. We see that 
	\begin{equation*}
	\sum_{k=-\infty }^{\infty }\sum_{m\in \mathbb{Z}^{n}}|\lambda
	_{k,m}||\langle \psi _{k,m},\varphi \rangle |=I_{1}+I_{2},
	\end{equation*}%
	where 
	\begin{equation*}
	I_{1}=\sum_{k=-\infty }^{0}\sum_{m\in \mathbb{Z}^{n}}|\lambda
	_{k,m}||\langle \psi _{k,m},\varphi \rangle |\quad \text{and}\quad
	I_{2}=\sum_{k=1}^{\infty }\sum_{m\in \mathbb{Z}^{n}}|\lambda _{k,m}||\langle
	\psi _{k,m},\varphi \rangle |.
	\end{equation*}%
	It suffices to show that both $I_{1}$ and $I_{2}$ are dominated by 
	\begin{equation*}
	c\big\|\lambda |\dot{f}_{p,q}(\mathbb{R}^{n},\{t_{k}\})\big\|.
	\end{equation*}
	
	\textit{Estimate of }$I_{1}$. Let us recall the following estimate, see
	(3.18) in \cite{M07}. For any $L>0$, there exists a positive constant $M\in 
	\mathbb{N}$ such that for all $\varphi $, $\psi \in \mathcal{S}_{\infty }(%
	\mathbb{R}^{n}),i,k\in \mathbb{Z}$ and $m,h\in \mathbb{Z}^{n}$, 
	\begin{align*}
	&\left\vert \langle \varphi _{k,m},\psi _{i,h}\rangle \right\vert \\
	& \lesssim \big\|\varphi \big\|_{\mathcal{S}_{M+1}}\big\|\psi \big\|_{\mathcal{S}_{M+1}}%
	\Big(1+\frac{|2^{-k}m-2^{-i}h|^{n}}{\max (2^{-kn},2^{-in})}\Big)^{-L}\min
	\left( 2^{(i-k)nL},2^{(k-i)nL}\right) .
	\end{align*}%
	Therefore, 
	\begin{equation*}
	\left\vert \langle \psi _{k,m},\varphi \rangle \right\vert \lesssim \big\|%
	\varphi \big\|_{\mathcal{S}_{M+1}}\big\|\psi \big\|_{\mathcal{S}_{M+1}}\Big(%
	1+\frac{|2^{-k}m|^{n}}{\max (1,2^{-kn})}\Big)^{-L}2^{-|k|nL}.
	\end{equation*}%
	Our estimate employs partially some decomposition techniques already used in 
	\cite{FJ90} and \cite{Ky03}. {For each $j$}$\in \mathbb{N}${\ we define } 
	\begin{equation*}
	{\Omega _{j}=\{m\in \mathbb{Z}^{n}:2^{j-1}<|m|\leq 2^{j}\}\quad }\text{and}{%
		\quad \Omega _{0}=\{m\in \mathbb{Z}^{n}:|m|\leq 1\}.}
	\end{equation*}%
	Thus, 
	\begin{align*}
	I_{1}& \lesssim \sum_{k=-\infty }^{0}2^{knL}\sum_{m\in \mathbb{Z}^{n}}\frac{%
		|\lambda _{k,m}|}{\big(1+|m|\big)^{nL}} \\
	& =c\sum_{k=-\infty }^{0}2^{knL}\sum\limits_{j=0}^{\infty }\sum\limits_{m\in
		\Omega _{j}}\frac{|\lambda _{k,m}|}{\big(1+|m|\big)^{nL}} \\
	& \lesssim \sum_{k=-\infty }^{0}2^{knL}\sum\limits_{j=0}^{\infty
	}2^{-nLj}\sum\limits_{m\in \Omega _{j}}|\lambda _{k,m}|.
	\end{align*}%
	Let $0<\varrho <\min (1,\theta )$ be such that $1/\varrho =1/\tau +1/\sigma
	_{1}$\ with $0<\tau <\min \big(\frac{1}{1-1/\sigma _{1}},p\big)$.\ We have 
	\begin{align*}
	&I_{1} \lesssim \sum_{k=-\infty }^{0}\sum\limits_{j=0}^{\infty
	}2^{-nL(j-k)}\Big(\sum\limits_{m\in \Omega _{j}}|\lambda _{k,m}|^{\varrho }\Big)%
	^{1/\varrho } \\
	& =c\sum_{k=-\infty }^{0}\sum\limits_{j=0}^{\infty }2^{(1/\varrho
		-L)nj+knL}\Big(2^{(k-j)n}\int_{\cup _{z\in \Omega _{j}}Q_{k,z}}\sum\limits_{m\in
		\Omega _{j}}|\lambda _{k,m}|^{\varrho }\chi _{k,m}(y)dy\Big)^{1/\varrho }.
	\end{align*}%
	Let $y\in \cup _{z\in \Omega _{j}}Q_{k,z}$ and $x\in Q_{0,0}$. Then $y\in
	Q_{k,z}$ for some $z\in \Omega _{j}$ and ${2^{j-1}<|z|\leq 2^{j}}$. From
	this it follows that 
	\begin{align*}
	\left\vert y-x\right\vert & \leq |y-2^{-k}z|+|x-2^{-k}z| \\
	& \leq \sqrt{n}\text{ }2^{-k}+\left\vert x\right\vert +2^{-k}\left\vert
	z\right\vert \\
	& \leq 2^{j-k+\delta _{n}},\quad \delta _{n}\in \mathbb{N},
	\end{align*}%
	which implies that $y$ is located in the ball $B(x,2^{j-k+\delta _{n}})$. In
	addition, from the fact that 
	\begin{equation*}
	\left\vert y\right\vert \leq \left\vert y-x\right\vert +\left\vert
	x\right\vert \leq 2^{j-k+\delta _{n}}+1\leq 2^{j-k+c_{n}},\quad c_{n}\in 
	\mathbb{N},
	\end{equation*}%
	we have that $y$ is located in the ball $B(0,2^{j-k+c_{n}})$. Therefore, by H{\"{o}}lder's inequality 
	\begin{align*}
	& \Big(2^{(k-j)n}\int_{\cup _{z\in \Omega _{j}}Q_{k,z}}\sum\limits_{m\in
		\Omega _{j}}|\lambda _{k,m}|^{\varrho }\chi _{k,m}(y)dy\Big)^{1/\varrho } \\
	& \leq \Big(2^{(k-j)n}\int_{B(x,2^{j-k+c_{n}})}\sum\limits_{m\in \Omega
		_{j}}|\lambda _{k,m}|^{\tau }t_{k}^{\tau }\chi _{k,m}(y)dy\Big)^{1/\tau }M_{B(0,2^{j-k+c_{n}}),\sigma _{1}}(t_{k}^{-1}) \\
	& \lesssim \mathcal{M}_{\tau }\big(\sum\limits_{m\in \mathbb{Z}%
		^{n}}t_{k}\left\vert \lambda _{k,m}\right\vert \chi _{k,m}\big)%
	(x)M_{B(0,2^{j-k+c_{n}}),\sigma _{1}}(t_{k}^{-1}).
	\end{align*}%
	Since $t_{k}^{-\sigma _{1}}\in A_{(p/\theta )^{\prime }}(\mathbb{R}^{n})$, $%
	k\in \mathbb{Z}$, by Lemma {\ref{Ap-Property}/(iii)}, $\mathrm{\eqref{Asum1}}
	$ and $\mathrm{\eqref{Asum2}}$ we obtain 
	\begin{align*}
	M_{B(0,2^{j-k+c_{n}}),\sigma _{1}}(t_{k}^{-1})& \lesssim 2^{(j-k)\frac{n}{p}%
	}M_{B(0,1),\sigma _{1}}(t_{k}^{-1}) \\
	& \lesssim 2^{(j-k)\frac{n}{p}}\left( M_{B(0,1),p}(t_{k})\right) ^{-1} \\
	& \lesssim 2^{(j-k)\frac{n}{p}-k\alpha _{2}}\left( M_{B(0,1),\sigma
		_{2}}(t_{0})\right) ^{-1}
	\end{align*}%
	for any $k\leq 0$ and any $j\in \mathbb{N}_{0}$. Hence, for any $L$ large
	enough, 
	\begin{equation*}
	I_{1}\lesssim \sum_{k=-\infty }^{0}2^{k(nL-\alpha _{2}-n/p)}\mathcal{M}%
	_{\tau }\big(\sum\limits_{m\in \mathbb{Z}^{n}}t_{k}|\lambda _{k,m}|\chi
	_{k,m}\big)(x),\quad x\in Q_{0,0}.
	\end{equation*}%
	The last term is bounded in the $L_{p}(Q_{0,0})$-quasi-norm by $c\big\|%
	\lambda |\dot{f}_{p,q}(\mathbb{R}^{n},\{t_{k}\})\big\|$ with the help of
	Theorem {\ref{Maximal}.}
	
	\textit{Estimate of }$I_{2}$. We have 
	\begin{equation*}
	\left\vert \langle \psi _{k,m},\varphi \rangle \right\vert \lesssim 2^{-knL}%
	\big\|\varphi \big\|_{\mathcal{S}_{M+1}}\big\|\psi \big\|_{\mathcal{S}_{M+1}}%
	\big(1+2^{-kn}|m|^{n}\big)^{-L},\quad k\geq 1.
	\end{equation*}%
	{For each $j,k$}$\in \mathbb{N}${, define } 
	\begin{equation*}
	{\Omega _{j,k}=\{m\in \mathbb{Z}^{n}:2^{j+k-1}<|m|\leq 2^{j+k}\}}\text{\quad
		and\quad }{\Omega _{0,k}=\{m\in \mathbb{Z}^{n}:|m|\leq 2^{k}\}.}
	\end{equation*}%
	Then we find 
	\begin{align*}
	&I_{2} \lesssim \sum_{k=1}^{\infty }2^{-knL}\sum_{m\in \mathbb{Z}^{n}}\frac{%
		|\lambda _{k,m}|}{\big(1+2^{-k}|m|\big)^{nL}} \\
	& =c\sum_{k=1}^{\infty }2^{-knL}\sum\limits_{j=0}^{\infty }\sum\limits_{m\in
		\Omega _{j,k}}\frac{|\lambda _{k,m}|}{\big(1+2^{-k}|m|\big)^{nL}} \\
	& \leq c\sum_{k=1}^{\infty }2^{-knL}\sum\limits_{j=0}^{\infty
	}2^{-nLj}\sum\limits_{m\in {\Omega _{j,k}}}|\lambda _{k,m}|.
	\end{align*}%
	Let $0<\varrho <\min (1,\theta )$ be such that $1/\varrho =1/\tau +1/\sigma
	_{1}$ with $0<\tau <p$. Using the embedding $\ell _{\varrho }\hookrightarrow
	\ell _{1}$ we find that  $I_{2}$ does not exceed
	\begin{equation*}
	c \sum_{k=1}^{\infty }2^{-knL}\sum\limits_{j=0}^{\infty
	}2^{-nLj}\big(\sum\limits_{m\in {\Omega _{j,k}}}|\lambda _{k,m}|^{\varrho }%
	\big)^{1/\varrho }
	\end{equation*}
	which is just the term
	\begin{equation*}
	c\sum_{k=1}^{\infty }2^{-knL}\sum\limits_{j=0}^{\infty }2^{(n/\varrho
		-nL)j}\Big(2^{(k-j)n}\int_{\cup _{z\in {\Omega _{j,k}}}Q_{k,z}}\sum%
	\limits_{m\in {\Omega _{j,k}}}|\lambda _{k,m}|^{\varrho }\chi _{k,m}(y)dy%
	\Big)^{1/\varrho }.
	\end{equation*}%
	Let $y\in \cup _{z\in {\Omega _{j,k}}}Q_{k,z}$ and $x\in Q_{0,0}$. Then $%
	y\in Q_{k,z}$ for some $z\in {\Omega _{j,k}}$ and ${2^{j-1}<2^{-k}|z|\leq
		2^{j}}$. From this it follows that 
	\begin{equation*}
	|y-x|\leq |y-2^{-k}z|+|x-2^{-k}|\leq \sqrt{n}\text{ }2^{-k}+|x|+2^{-k}|z|%
	\leq 2^{j+\delta _{n}},\quad \delta _{n}\in \mathbb{N},
	\end{equation*}%
	which implies that $y$ is located in the ball $B\left( x,2^{j+\delta
		_{n}}\right) $. In addition, from the fact that 
	\begin{equation*}
	\left\vert y\right\vert \leq \left\vert y-x\right\vert +\left\vert
	x\right\vert \leq 2^{j+\delta _{n}}+1\leq 2^{j+c_{n}},\quad c_{n}\in \mathbb{%
		N},
	\end{equation*}%
	we have that $y$ is located in the ball $B(0,2^{j+c_{n}})$. Therefore, 
	\begin{align*}
	& \Big(2^{(k-j)n}\int_{\cup _{z\in {\Omega _{j,k}}}Q_{k,z}}\sum\limits_{m\in 
		{\Omega _{j,k}}}|\lambda _{k,m}|^{\varrho }\chi _{k,m}(y)dy\Big)^{1/\varrho }
	\\
	& \leq 2^{kn/\varrho }\Big(2^{-jn}\int_{B(x,2^{j+\delta
			_{n}})}\sum\limits_{m\in {\Omega _{j,k}}}|\lambda _{k,m}|^{\tau }t_{k}^{\tau
	}\chi _{k,m}(y)dy\Big)^{1/\tau }M_{B(0,2^{j+c_{n}}),\sigma
		_{1}}(t_{k}^{-1}) \\
	& \lesssim 2^{kn/\varrho }\mathcal{M}_{\tau }\big(\sum\limits_{m\in \mathbb{Z%
		}^{n}}t_{k}\lambda _{k,m}\chi _{k,m}\big)(x)M_{B(0,2^{j+c_{n}}),\sigma
		_{1}}(t_{k}^{-1}).
	\end{align*}%
	By \eqref{Asum1} and Lemma {\ref{Ap-Property}/(vi)} we obtain 
	\begin{align*}
	&	M_{B(0,2^{j+c_{n}}),\sigma _{1}}(t_{k}^{-1})\\
	&\lesssim 2^{-k\alpha _{1}}\big(%
	M_{B(0,2^{j+c_{n}}),p}(t_{0})\big)^{-1}\\
	&\lesssim 2^{j(n/p-n\delta /p)-k\alpha
		_{1}}\big(M_{B(0,1),p}(t_{0})\big)^{-1}.
	\end{align*}%
	Therefore 
	\begin{equation}
	I_{2}\lesssim \sum_{k=1}^{\infty }2^{-k(nL-n/\varrho +\alpha _{1})}\mathcal{M%
	}_{\tau }\big(t_{k}\sum\limits_{m\in \mathbb{Z}^{n}}\lambda _{k,m}\chi _{k,m}%
	\big)(x),\quad x\in Q_{0,0}  \label{I-two}
	\end{equation}%
	for any $L$ large enough. Now we take the $L_{p}(Q_{0,0})$-quasi-norm of
	both sides of $\mathrm{\eqref{I-two}}$ and then use Theorem {\ref{Maximal}},
	we obtain 
	\begin{equation*}
	I_{2}\lesssim {\big\|}\lambda |\dot{f}_{p,q}(\mathbb{R}^{n},\{t_{k}\}){\big\|%
	}.
	\end{equation*}%
	The proof is complete.
\end{proof}

For a sequence $\lambda =\{\lambda _{k,m}\}_{k\in \mathbb{Z},m\in \mathbb{Z}%
	^{n}}\subset \mathbb{C},0<r\leq \infty $ and a fixed $d>0$, set%
\begin{equation*}
\lambda _{k,m,r,d}^{\ast }=\Big(\sum_{h\in \mathbb{Z}^{n}}\frac{|\lambda
	_{k,h}|^{r}}{(1+2^{k}|2^{-k}h-2^{-k}m|)^{d}}\Big)^{1/r}
\end{equation*}%
and $\lambda _{r,d}^{\ast }:=\{\lambda _{k,m,r,d}^{\ast }\}_{k\in \mathbb{Z}%
	,m\in \mathbb{Z}^{n}}\subset \mathbb{C}$.

\begin{lem}
	\label{lamda-equi}Let $\alpha =(\alpha _{1},\alpha _{2})\in \mathbb{R}%
	^{2},0<\theta \leq p<\infty ,0<q<\infty ,\gamma \in \mathbb{Z}$ and $d>n$.
	Let $\{t_{k}\}$ be a $p$-admissible weight sequence satisfying $\mathrm{%
		\eqref{Asum1}}$ with $\sigma _{1}=\theta \left( p/\theta \right) ^{\prime }$
	and $\alpha _{1}\in \mathbb{R}$. Then%
	\begin{equation}
	\big\|\lambda _{\min(p,q),d}^{\ast }|\dot{f}_{p,q}(\mathbb{R}^{n},\{t_{k-\gamma }\})%
	\big\|\approx \big\|\lambda |\dot{f}_{p,q}(\mathbb{R}^{n},\{t_{k-\gamma }\})%
	\big\|.  \label{First}
	\end{equation}%
	In addition if $\{t_{k}\}$ satisfies $\mathrm{\eqref{Asum2}}$ with $\sigma
	_{2}\geq p$ and $\alpha _{2}\in \mathbb{R}$, then 
	\begin{equation}
	\big\|\lambda _{\min(p,q),d}^{\ast }|\dot{f}_{p,q}(\mathbb{R}^{n},\{t_{k-\gamma }\})%
	\big\|\lesssim \big\|\lambda |\dot{f}_{p,q}(\mathbb{R}^{n},\{t_{k}\})\big\|.
	\label{Second}
	\end{equation}%

\end{lem}

\begin{proof}
	The proof is similar to that given in {\cite{D20}}. First we prove $\mathrm{%
		\eqref{First}}$. Obviously, 
	\begin{equation*}
	\big\|\lambda |\dot{f}_{p,q}(\mathbb{R}^{n},\{t_{k-\gamma }\})\big\|\leq %
	\big\|\lambda _{p,d}^{\ast }|\dot{f}_{\min(p,q),q}(\mathbb{R}^{n},\{t_{k-\gamma }\})%
	\big\|.
	\end{equation*}%
	Let $n\min(p,q)/d<a<\min(p,q),{j}\in \mathbb{N}$ and $m\in \mathbb{Z}^{n}$. {Define }%
	\begin{equation*}
	{\Omega _{j,m}=\{h\in \mathbb{Z}^{n}:2^{j-1}<|h-m|\leq 2^{j}\},\quad \Omega
		_{0,m}=\{h\in \mathbb{Z}^{n}:|h-m|\leq 1\}}.
	\end{equation*}%
	Then%
	\begin{align*}
	\sum_{h\in \mathbb{Z}^{n}}\frac{|\lambda _{k,h}|^{\min(p,q)}}{\big(1+|h-m|\big)^{d}}%
	& =\sum\limits_{j=0}^{\infty }\sum\limits_{h\in {\Omega _{j,m}}}\frac{%
		|\lambda _{k,h}|^{\min(p,q)}}{\big(1+|h-m|\big)^{d}} \\
	& \lesssim \sum\limits_{j=0}^{\infty }2^{-dj}\sum\limits_{h\in {\Omega _{j,m}%
	}}|\lambda _{k,h}|^{\min(p,q)} \\
	& \lesssim c\sum\limits_{j=0}^{\infty }2^{-dj}\Big(\sum\limits_{h\in {\Omega
			_{j,m}}}|\lambda _{k,h}|^{a}\Big)^{\min(p,q)/a}.
	\end{align*}%
	The last expression can be rewritten as 
	\begin{equation}
	c\sum\limits_{j=0}^{\infty }2^{(n\min(p,q)/a-d)j}\Big(2^{(k-j)n}\int_{\cup _{z\in {%
				\Omega _{j,m}}}Q_{k,z}}\sum\limits_{h\in {\Omega _{j,m}}}|\lambda
	_{k,h}|^{a}\chi _{k,h}(y)dy\Big)^{\min(p,q)/a}.  \label{estimate-lamda}
	\end{equation}%
	Let $y\in \cup _{z\in {\Omega _{j,m}}}Q_{k,z}$ and $x\in Q_{k,m}$. Then $%
	y\in Q_{k,z}$ for some $z\in {\Omega _{j,m}}$ which implies that ${%
		2^{j-1}<|z-m|\leq 2^{j}}$. From this it follows that%
	\begin{align*}
	\left\vert y-x\right\vert & \leq \left\vert y-2^{-k}z\right\vert +\left\vert
	x-2^{-k}z\right\vert \\
	& \leq \sqrt{n}\text{ }2^{-k}+\left\vert x-2^{-k}m\right\vert
	+2^{-k}\left\vert z-m\right\vert \\
	& \leq 2^{j-k+\delta _{n}},\quad \delta _{n}\in \mathbb{N},
	\end{align*}%
	which implies that $y$ is located in the ball $B(x,2^{j-k+\delta _{n}})$.
	Therefore, $\mathrm{\eqref{estimate-lamda}}$ can be estimated by%
	\begin{equation*}
	c\mathcal{M}_{a}\big(\sum\limits_{h\in \mathbb{Z}^{n}}\lambda _{k,h}\chi
	_{k,h}\big)(x),
	\end{equation*}%
	where the positive constant $c$ is independent of $k$. Consequently%
	\begin{equation}
	\big\|\lambda _{\min(p,q),d}^{\ast }|\dot{f}_{p,q}(\mathbb{R}^{n},\{t_{k-\gamma }\})%
	\big\|  \label{estimate-lamda1}
	\end{equation}%
	does not exceed 
	\begin{equation*}
	c\Big\|\Big(\sum_{k=-\infty }^{\infty }2^{knq/2}\big(t_{k-\gamma }\mathcal{M}%
	_{a}\big(\sum\limits_{h\in \mathbb{Z}^{n}}\lambda _{k,h}\chi _{k,h}\big)\big)%
	^{q}\Big)^{1/q}|L_{p}(\mathbb{R}^{n})\Big\|.
	\end{equation*}%
	Applying Lemma {\ref{key-estimate1}} we estimate \eqref{estimate-lamda1} by 
	\begin{equation*}
	c\Big\|\Big(\sum_{k=-\infty }^{\infty }2^{knq/2}\sum\limits_{h\in \mathbb{Z}%
		^{n}}|\lambda _{k,h}|^{q}\chi _{k,h}\Big)^{1/q}|L_{p}(\mathbb{R}%
	^{n},t_{k-\gamma })\Big\|=c\big\|\lambda |\dot{f}_{p,q}(\mathbb{R}%
	^{n},\{t_{k-\gamma }\})\big\|.
	\end{equation*}%
	To prove $\mathrm{\eqref{Second}}$ we use again Lemma {\ref{key-estimate1}
		combined with Remark \ref{r-estimates copy(1)}.}
\end{proof}

Now we have the following result which is called the $\varphi $-transform
characterization in the sense of Frazier and Jawerth. It will play an
important role in the rest of the paper.

\begin{thm}
	\label{phi-tran}Let $\alpha =(\alpha _{1},\alpha _{2})\in \mathbb{R}%
	^{2},0<\theta \leq p<\infty $ and$\ 0<q<\infty $. Let $\{t_{k}\}\in \dot{X}%
	_{\alpha ,\sigma ,p}$ be a $p$-admissible weight sequence with $\sigma
	=(\sigma _{1}=\theta \left( p/\theta \right) ^{\prime },\sigma _{2}\geq p)$%
	.\ Let $\varphi $, $\psi \in \mathcal{S}(\mathbb{R}^{n})$ satisfying $%
	\mathrm{\eqref{Ass1}}$\ through\ $\mathrm{\eqref{Ass3}}$. The operators 
	\begin{equation*}
	S_{\varphi }:\dot{F}_{p,q}(\mathbb{R}^{n},\{t_{k}\})\rightarrow \dot{f}%
	_{p,q}(\mathbb{R}^{n},\{t_{k}\})
	\end{equation*}%
	and 
	\begin{equation*}
	T_{\psi }:\dot{f}_{p,q}(\mathbb{R}^{n},\{t_{k}\})\rightarrow \dot{F}_{p,q}(%
	\mathbb{R}^{n},\{t_{k}\})
	\end{equation*}%
	are bounded. Furthermore, $T_{\psi }\circ S_{\varphi }$ is the identity on $%
	\dot{F}_{p,q}(\mathbb{R}^{n},\{t_{k}\})$.
\end{thm}

\begin{proof}
	The proof is a straightforward adaptation of {\cite[Theorem 2.2]{FJ90} and 
		\cite{D20}. For any $f\in \mathcal{S}_{\infty }^{\prime }(\mathbb{R}^{n})$
		we put $\sup (f):=\{\sup_{k,m}(f)\}_{k\in \mathbb{Z},m\in \mathbb{Z}^{n}}$
		where 
		\begin{equation*}
		\sup_{k,m}(f)=2^{-kn/2}\sup_{y\in Q_{k,m}}|\widetilde{\varphi _{k}}\ast
		f(y)|,\quad k{\in \mathbb{Z},m\in \mathbb{Z}^{n}.}
		\end{equation*}%
		For any $\gamma $}${\in }\mathbb{N}_{0}${, we define the sequence $%
		\inf_{\gamma }(f)=\{\inf_{k,m,\gamma }(f)\}_{k\in \mathbb{Z},m\in \mathbb{Z}%
			^{n}}$ by setting 
		\begin{equation*}
		\text{inf}_{k,m,\gamma }(f)=2^{-kn/2}\max \{\inf_{y\in \tilde{Q}}|\widetilde{%
			\varphi _{k}}\ast f(y)|:\tilde{Q}\subset Q_{k,m},l(\tilde{Q})=2^{-k-\gamma
		}\},
		\end{equation*}%
		where $k\in \mathbb{Z}$, $ m\in \mathbb{Z}^{n} $ and $\widetilde{\varphi _{k}}=2^{kn}\overline{\varphi (-2^{k}\cdot )}$, $k\in \mathbb{Z}$. }
	
	\textit{Step 1. }{In this step we prove that 
		\begin{equation*}
		\big\|\text{inf}_{\gamma }(f)|\dot{f}_{p,q}(\mathbb{R}^{n},\{t_{k}\})\big\|%
		\lesssim \big\|f|\dot{F}_{p,q}(\mathbb{R}^{n},\{t_{k}\})\big\|.
		\end{equation*}%
	}Define a sequence $\lambda =\{\lambda _{j,h}\}_{j{\in \mathbb{Z},h\in 
			\mathbb{Z}^{n}}}$ by 
	\begin{equation*}
	\lambda _{j,h}=2^{-jn/2}\inf_{y\in Q_{j,h}}|\widetilde{\varphi _{j-\gamma }}%
	\ast f(y)|,\quad j{\in \mathbb{Z},h\in \mathbb{Z}^{n}.}
	\end{equation*}%
	Then for all $0<r<\infty $, any $k{\in \mathbb{Z},m\in \mathbb{Z}^{n}}$ and
	a fixed $\lambda >n$, we have%
	\begin{equation*}
	\text{inf}_{k,m,\gamma }(f)\chi _{k,m}\lesssim 2^{\gamma
		(d/r+n/2)}\sum_{h\in {\mathbb{Z}^{n},Q}_{k+\gamma ,h}\subset Q_{k,m}}\lambda
	_{k+\gamma ,h,r,d}^{\ast }\chi _{k+\gamma ,h}.
	\end{equation*}%
	Picking $r=\min(p,q)$, we obtain%
	\begin{align*}
	& \big\|\text{inf}_{\gamma }(f)|\dot{f}_{p,q}(\mathbb{R}^{n},\{t_{k}\})\big\|
	\\
	& \lesssim 2^{\gamma (d/r+n/2)}\Big\|\Big(\sum_{k=-\infty }^{\infty
	}2^{knq/2}\sum\limits_{h\in \mathbb{Z}^{n}}(t_{k}\lambda _{k+\gamma
		,h,r,d}^{\ast })^{q}\chi _{k+\gamma ,h}\Big)^{1/q}|L_{p}(\mathbb{R}^{n})%
	\Big\| \\
	& =c2^{\gamma d/r}\Big\|\Big(\sum_{k=-\infty }^{\infty
	}2^{knq/2}\sum\limits_{h\in \mathbb{Z}^{n}}(t_{k-\gamma }\lambda
	_{k,h,r,d}^{\ast })^{q}\chi _{k,h}\Big)^{1/q}|L_{p}(\mathbb{R}^{n})\Big\|.
	\end{align*}%
	We apply Lemma {\ref{lamda-equi} }to estimate the last expression by%
	\begin{align*}
	& c2^{\gamma d/r}\Big\|\Big(\sum_{k=-\infty }^{\infty
	}2^{knq/2}\sum\limits_{h\in \mathbb{Z}^{n}}(t_{k-\gamma }\lambda
	_{k,h})^{q}\chi _{k,h}\Big)^{1/q}|L_{p}(\mathbb{R}^{n})\Big\| \\
	& \lesssim \Big\|\Big(\sum_{k=-\infty }^{\infty }t_{k}^{q}\big|(\widetilde{%
		\varphi _{k}}\ast f)\big|^{q}\Big)^{1/q}|L_{p}(\mathbb{R}^{n})\Big\| \\
	& \lesssim \big\|f|\dot{F}_{p,q}(\mathbb{R}^{n},\{t_{k}\})\big\|.
	\end{align*}%
	
	\textit{Step 2.} {We will prove that }%
	\begin{equation}
	\big\|\text{inf}_{\gamma }(f)|\dot{f}_{p,q}(\mathbb{R}^{n},\{t_{k}\})\big\|%
	\approx \big\|f|\dot{F}_{p,q}(\mathbb{R}^{n},\{t_{k}\})\big\|\approx \big\|%
	\sup (f)|\dot{f}_{p,q}(\mathbb{R}^{n},\{t_{k}\})\big\|.  \label{step21}
	\end{equation}%
	{Applying Lemma A.4 of \cite{FJ90}, see also Lemma 8.3 of \cite{M07}, to the
		function }$(\widetilde{\varphi _{k}}\ast f)(2^{-j}x)$ {we obtain }%
	\begin{equation*}
	\text{{inf}}_{\gamma }{(f)_{\min(p,q),d}^{\ast }\approx \sup (f)_{\min(p,q),d}^{\ast }.}
	\end{equation*}%
	{Hence for $\gamma >0$ sufficiently large we obtain by applying Lemma \ref%
		{lamda-equi}, }%
	\begin{equation*}
	\big\|\text{inf}_{\gamma }(f)_{\min(p,q),d}^{\ast }|\dot{f}_{p,q}(\mathbb{R}%
	^{n},\{t_{k}\})\big\|{\approx }\big\|\text{inf}_{\gamma }(f)|\dot{f}_{p,q}(%
	\mathbb{R}^{n},\{t_{k}\})\big\|
	\end{equation*}%
	{\ and}%
	\begin{equation*}
	\big\|\sup (f)_{\min(p,q),d}^{\ast }|\dot{f}_{p,q}(\mathbb{R}^{n},\{t_{k}\})\big\|{%
		\approx }\big\|\sup (f)|\dot{f}_{p,q}(\mathbb{R}^{n},\{t_{k}\})\big\|.
	\end{equation*}%
	{\ Therefore, 
		\begin{equation}
		\big\|\text{inf}_{\gamma }(f)|\dot{f}_{p,q}(\mathbb{R}^{n},\{t_{k}\})\big\|%
		\approx \big\|\sup (f)|\dot{f}_{p,q}(\mathbb{R}^{n},\{t_{k}\})\big\|.
		\label{sup-inf}
		\end{equation}%
		From the definition of the spaces }$\dot{F}_{p,q}(\mathbb{R}^{n},\{t_{k}\})$
	it follows that%
	\begin{equation*}
	\big\|f|\dot{F}_{p,q}(\mathbb{R}^{n},\{t_{k}\})\big\|\lesssim \big\|\sup (f)|%
	\dot{f}_{p,q}(\mathbb{R}^{n},\{t_{k}\})\big\|.
	\end{equation*}%
	Consequently $\mathrm{\eqref{sup-inf}}$ and Step 1 yield $\mathrm{%
		\eqref{step21}.}$
	
	\textit{Step 3}. In this step we prove the boundedness of $S_{\varphi }$ and 
	$T_{\psi }$. We have%
	\begin{equation*}
	|(S_{\varphi }f)_{k,m}|=|\langle f,\varphi _{k,m}\rangle |=2^{-kn/2}|f\ast 
	\widetilde{\varphi _{k}}(2^{-k}m)|\leq \sup_{k,m}(f).
	\end{equation*}%
	Step 2 yields that%
	\begin{equation*}
	\big\|S_{\varphi }f|\dot{f}_{p,q}(\mathbb{R}^{n},\{t_{k}\})\big\|\lesssim %
	\big\|f|\dot{F}_{p,q}(\mathbb{R}^{n},\{t_{k}\})\big\|.
	\end{equation*}%
	To prove the boundedness of $T_{\psi }$ suppose $\lambda =\{\lambda
	_{j,h}\}_{j\in \mathbb{Z},h\in \mathbb{Z}^{n}}$ and 
	\begin{equation*}
	T_{\psi }\lambda =\sum_{j=-\infty }^{\infty }\sum_{h\in \mathbb{Z}%
		^{n}}\lambda _{j,h}\psi _{j,h}.
	\end{equation*}%
	Obviously%
	\begin{equation*}
	\widetilde{\varphi _{k}}\ast T_{\psi }\lambda =\sum_{j=k-1}^{k+1}\sum_{h\in 
		\mathbb{Z}^{n}}\lambda _{j,h}\widetilde{\varphi _{k}}\ast \psi _{j,h}.
	\end{equation*}%
	Since $\widetilde{\varphi }$ and $\psi $ belong to $\mathcal{S}(\mathbb{R}%
	^{n})$ we obtain%
	\begin{equation*}
	|\widetilde{\varphi _{k}}\ast \psi _{j,h}(x)|\lesssim
	2^{jn/2}(1+2^{j}|x-2^{-j}h|)^{-d/\min (1,\min(p,q))},\quad d>n,
	\end{equation*}%
	where the implicit constant is independent of $j,k,h$ and $x$. Therefore, if 
	$x\in Q_{k+1,z}\subset Q_{k,m}\subset Q_{k-1,l}$, $z,l\in \mathbb{Z}^{n}$,
	then we obtain 
	\begin{equation*}
	|\widetilde{\varphi _{k}}\ast T_{\psi }\lambda (x)|\lesssim
	2^{kn/2}\sum_{j=k-1}^{k+1}\sum_{h\in \mathbb{Z}^{n}}\frac{|\lambda _{j,h}|}{%
		(1+2^{j}|x-2^{-j}h|)^{d/\min (1,\min(p,q))}}.
	\end{equation*}%
	Assume that $0<\min(p,q)\leq 1$. Using the inequality%
	\begin{equation*}
	\Big(\sum_{h\in \mathbb{Z}^{n}}|a_{h}|\Big)^{\min(p,q)}\leq \sum_{h\in \mathbb{Z}%
		^{n}}|a_{h}|^{\min(p,q)},\quad \{a_{h}\}_{h\in \mathbb{Z}^{n}}\subset \mathbb{C},
	\end{equation*}%
	we obtain%
	\begin{equation}
	|\widetilde{\varphi _{k}}\ast T_{\psi }\lambda (x)|\lesssim
	2^{kn/2}\sum_{j=k-1}^{k+1}\Big(\sum_{h\in \mathbb{Z}^{n}}\frac{|\lambda
		_{j,h}|^{\min(p,q)}}{(1+2^{j}|x-2^{-j}h|)^{d}}\Big)^{1/\min(p,q)}.  \label{est-T}
	\end{equation}%
	Now if $\min(p,q)>1$, then by the H\"{o}lder inequality and the fact that 
	\begin{equation*}
	\sum_{h\in \mathbb{Z}^{n}}\frac{1}{(1+2^{j}|x-2^{-j}h|)^{d}}\lesssim 1,
	\end{equation*}%
	we also have \eqref{est-T} with $\min(p,q)>1$. Hence if $x\in Q_{k+1,z}\subset
	Q_{k,m}\subset Q_{k-1,l}$, then we have%
	\begin{equation*}
	|\widetilde{\varphi _{k}}\ast T_{\psi }\lambda (x)|\lesssim 2^{kn/2}(\lambda
	_{k-1,l,\min(p,q),d}^{\ast }+\lambda _{k,m,\min(p,q),d}^{\ast }+\lambda _{k+1,z,\min(p,q),d}^{\ast
	}).
	\end{equation*}%
	Consequently%
	\begin{equation*}
	\big\|T_{\psi }\lambda |\dot{F}_{p,q}(\mathbb{R}^{n},\{t_{k}\})\big\|%
	\lesssim \sum_{i=-1}^{1}I_{i},
	\end{equation*}%
	where%
	\begin{equation*}
	I_{i}=\Big\|\Big(\sum_{k=-\infty }^{\infty }2^{knq/2}\big\|\sum\limits_{h\in 
		\mathbb{Z}^{n}}(t_{k+i}\lambda _{k,h,\min(p,q),d}^{\ast })^{q}\chi _{k,h}\Big)%
	^{1/q}|L_{p}(\mathbb{R}^{n})\Big\|,
	\end{equation*}%
	$i=-1,0,1$.	Applying \eqref{Second} we obtain%
	\begin{equation*}
	\big\|T_{\psi }\lambda |\dot{F}_{p,q}(\mathbb{R}^{n},\{t_{k}\})\big\|%
	\lesssim \big\|\lambda |\dot{f}_{p,q}(\mathbb{R}^{n},\{t_{k}\})\big\|.
	\end{equation*}%
	The proof is complete.
\end{proof}

\begin{rem}
	This theorem can then be exploited to obtain a variety of results for the $%
	\dot{F}_{p,q}(\mathbb{R}^{n},\{t_{k}\})$ spaces, where arguments can be
	equivalently transferred to the sequence space, which is often more
	convenient to handle. More precisely, under the same hypothesis of the last
	theorem, 
	\begin{equation*}
	\big\|\{\langle f,\varphi _{k,m}\rangle \}_{k\in \mathbb{Z},m\in \mathbb{Z}%
		^{n}}|\dot{f}_{p,q}(\mathbb{R}^{n},\{t_{k}\})\big\|\approx \big\|\{f|\dot{F}%
	_{p,q}(\mathbb{R}^{n},\{t_{k}\})\big\|.
	\end{equation*}
\end{rem}

From Theorem \ref{phi-tran}, we obtain the next important property of the
function spaces\ $\dot{F}_{p,q}(\mathbb{R}^{n},\{t_{k}\})$.

\begin{cor}
	\label{Indpendent}Let $\alpha =(\alpha _{1},\alpha _{2})\in \mathbb{R}%
	^{2},0<\theta \leq p<\infty $ and $0<q<\infty $. Let $\{t_{k}\}\in \dot{X}%
	_{\alpha ,\sigma ,p}$ be a $p$-admissible weight sequence with $\sigma
	=(\sigma _{1}=\theta \left( p/\theta \right) ^{\prime },\sigma _{2}\geq p)$.
	The definition of the spaces $\dot{F}_{p,q}(\mathbb{R}^{n},\{t_{k}\})$ is
	independent of the choices of $\varphi \in \mathcal{S}(\mathbb{R}^{n})$
	satisfying $\mathrm{\eqref{Ass1}}$\ and\ $\mathrm{\eqref{Ass2}}$.
\end{cor}

\begin{lem}
	Let $0<\theta \leq p<\infty $ and $0<q<\infty $. Let $\{t_{k}\}$ be a $p$%
	-admissible weight sequence. $\dot{f}_{p,q}(\mathbb{R}^{n},\{t_{k}\})$ are
	quasi-Banach spaces. They are Banach spaces if $1\leq p<\infty $ and $1\leq
	q<\infty $.
\end{lem}

\begin{proof}
	We need only to assume that $1\leq p<\infty $ and $1\leq q<\infty $. Let $%
	\{\lambda ^{(i)}\}_{i\in \mathbb{N}_{0}}$ be sequence in $\dot{f}_{p,q}(%
	\mathbb{R}^{n},\{t_{k}\})$\ such that%
	\begin{equation*}
	\sum_{i=0}^{\infty }\big\|\lambda ^{(i)}|\dot{f}_{p,q}(\mathbb{R}%
	^{n},\{t_{k}\})\big\|<\infty .
	\end{equation*}%
	W{e obtain}%
	\begin{equation*}
	\sum_{i=0}^{\infty }\Big\|\Big(\sum_{k=-\infty }^{\infty }\sum\limits_{m\in 
		\mathbb{Z}^{n}}2^{knq/2}t_{k}^{q}|\lambda _{k,m}^{(i)}|^{q}\chi _{k,m}\Big)%
	^{1/q}|L_{p}(\mathbb{R}^{n})\Big\|<\infty .
	\end{equation*}%
	We put%
	\begin{equation*}
	A^{(i)}=\{A_{k}^{(i)}\}_{k\in \mathbb{Z}},\quad
	A_{k}^{(i)}=\sum\limits_{m\in \mathbb{Z}^{n}}2^{kn/2}t_{k}|\lambda
	_{k,m}^{(i)}|\chi _{k,m},\quad i\in \mathbb{N}_{0},
	\end{equation*}%
	it follows that $\sum_{i=0}^{\infty }A^{(i)}$ is a is absolutely convergent
	in $L_{p}(\ell _{q})$, so the sequence $\{A^{(i)}\}_{i\in \mathbb{N}_{0}}$
	converges in $L_{p}(\ell _{q})$ and%
	\begin{equation*}
	\Big\|\sum_{i=0}^{\infty }A^{(i)}|L_{p}(\ell _{q})\Big\|<\infty .
	\end{equation*}%
	Then 
	\begin{equation*}
	\sum\limits_{m\in \mathbb{Z}^{n}}t_{k}\chi _{k,m}\sum_{i=0}^{\infty
	}|\lambda _{k,m}^{(i)}|=\sum_{i=0}^{\infty }\sum\limits_{m\in \mathbb{Z}%
		^{n}}t_{k}|\lambda _{k,m}^{(i)}|\chi _{k,m}<\infty ,\quad k\in 
	\mathbb{Z},a.e.,
	\end{equation*}%
	which yields that%
	\begin{equation*}
	\sum_{i=0}^{\infty }\lambda _{k,m}^{(i)}<\infty ,\quad k\in \mathbb{Z},m\in 
	\mathbb{Z}^{n},a.e.
	\end{equation*}%
	Therefore%
	\begin{equation*}
	\Big\|\sum_{i=0}^{\infty }\lambda ^{(i)}|\dot{f}_{p,q}(\mathbb{R}%
	^{n},\{t_{k}\})\Big\|\leq \sum_{i=0}^{\infty }\Big\|A^{(i)}|L_{p}(\ell _{q})%
	\Big\|<\infty .
	\end{equation*}%
	This completes the proof.
\end{proof}

Applying this lemma and Theorem \ref{phi-tran} we obtain the following
useful properties of these function spaces, see \cite{GJN17}.

\begin{thm}
	Let $\alpha =(\alpha _{1},\alpha _{2})\in \mathbb{R}^{2},0<\theta \leq
	p<\infty $ and $0<q<\infty $. Let $\{t_{k}\}\in \dot{X}_{\alpha ,\sigma
		,p}$ be a $p$-admissible weight sequence with $\sigma =(\sigma _{1}=\theta
	\left( p/\theta \right) ^{\prime },\sigma _{2}\geq p)$. $\dot{F}_{p,q}(%
	\mathbb{R}^{n},\{t_{k}\})$ are quasi-Banach spaces. They are Banach spaces
	if $1\leq p<\infty $ and $1\leq q<\infty $.
\end{thm}

\begin{thm}
	\label{embeddings-S-inf}Let $0<\theta \leq p<\infty $ and $0<q<\infty $.
	Let\ $\{t_{k}\}\in \dot{X}_{\alpha ,\sigma ,p}$ be a $p$-admissible weight
	sequence with $\sigma =(\sigma _{1}=\theta \left( p/\theta \right) ^{\prime
	},\sigma _{2}\geq p)$ and $\alpha =(\alpha _{1},\alpha _{2})\in \mathbb{R}%
	^{2}$.\ \newline
	$\mathrm{(i)}$ We have the embedding%
	\begin{equation}
	\mathcal{S}_{\infty }(\mathbb{R}^{n})\hookrightarrow \dot{F}_{p,q}(\mathbb{R}%
	^{n},\{t_{k}\}).  \label{embedding}
	\end{equation}%
	In addition $\mathcal{S}_{\infty }(\mathbb{R}^{n})$ is dense in $\dot{F}%
	_{p,q}(\mathbb{R}^{n},\{t_{k}\}$.\newline
	$\mathrm{(ii)}$ We have the embedding 
	\begin{equation*}
	\dot{F}_{p,q}(\mathbb{R}^{n},\{t_{k}\})\hookrightarrow \mathcal{S}_{\infty
	}^{\prime }(\mathbb{R}^{n}).
	\end{equation*}
\end{thm}

\begin{proof}
	The proof is a variant of that given for Besov spaces in \cite{D20}{. }For
	the convenience of the reader, we give some details. The embedding $\mathrm{%
		\eqref{embedding}}$ follows by 
	\begin{equation*}
	\mathcal{S}_{\infty }(\mathbb{R}^{n})\hookrightarrow \dot{B}_{p,\min (p,q)}(%
	\mathbb{R}^{n},\{t_{k}\})\hookrightarrow \dot{F}_{p,q}(\mathbb{R}%
	^{n},\{t_{k}\}),
	\end{equation*}%
	see \cite{D20} for the first embedding. Now, we prove the density of $\mathcal{S}%
	_{\infty }(\mathbb{R}^{n})$ in $\dot{F}_{p,q}(\mathbb{R}^{n},\{t_{k}\})$.
	Let $\varphi $, $\psi \in \mathcal{S}(\mathbb{R}^{n})$\ satisfying $\mathrm{%
		\eqref{Ass1}}$ through $\mathrm{\eqref{Ass3}}$ and $f\in \dot{F}_{p,q}(%
	\mathbb{R}^{n},\{t_{k}\})$. Let%
	\begin{equation*}
	f_{N}=\sum\limits_{k=-N}^{N}\tilde{\psi}_{k}\ast \varphi _{k}\ast f,\quad
	N\in \mathbb{N}.
	\end{equation*}%
	Observe that%
	\begin{equation*}
	\varphi _{j}\ast \tilde{\psi}_{k}=0,\quad \text{if\quad }k\notin
	\{j-1,j,j+1\}.
	\end{equation*}%
	Then, by Lemma \ref{key-estimate1},%
	\begin{align*}
	\big\|f_{N}|\dot{F}_{p,q}(\mathbb{R}^{n},\{t_{k}\})\big\|& =\big\|\Big(%
	\sum\limits_{|k|\leq N+1}|t_{k}(\varphi _{k}\ast \tilde{\psi}_{k}\ast \bar{%
		\varphi}_{k}\ast f)|^{q}\Big)^{1/q}|L_{p}(\mathbb{R}^{n})\big\| \\
	& \lesssim \big\|\Big(\sum\limits_{|k|\leq N+1}|t_{k}\mathcal{M}_{\tau
	}(\varphi _{k}\ast f)|^{q}\Big)^{1/q}|L_{p}(\mathbb{R}^{n})\big\| \\
	& \lesssim \big\|f|\dot{F}_{p,q}(\mathbb{R}^{n},\{t_{k}\})\big\|<\infty
	\end{align*}%
	for any $N\in \mathbb{N}$ where $\bar{\varphi}_{k}=\varphi _{k-1}+\varphi
	_{k}+\varphi _{k+1},k\in \mathbb{Z}$ and $0<\tau <\min (1,p,q)$. The first
	inequality follows by Lemma 2.4 of \cite{Dr15}. Consequently,%
	\begin{align*}
	\big\|f-f_{N}|\dot{F}_{p,q}(\mathbb{R}^{n},\{t_{k}\})\big\|& \leq \big\|\Big(%
	\sum\limits_{|k|\geq N+1}|t_{k}(\varphi _{k}\ast \tilde{\psi}_{k}\ast \bar{%
		\varphi}_{k}\ast f)|^{q}\Big)^{1/q}|L_{p}(\mathbb{R}^{n})\big\| \\
	& \lesssim \big\|\Big(\sum\limits_{|k|\geq N+1}|t_{k}\mathcal{M}_{\tau
	}(\varphi _{k}\ast f)|^{q}\Big)^{1/q}|L_{p}(\mathbb{R}^{n})\big\| \\
	& \lesssim \big\|\Big(\sum\limits_{|k|\geq N+1}|t_{k}(\varphi _{k}\ast
	f)|^{q}\Big)^{1/q}|L_{p}(\mathbb{R}^{n})\big\|,
	\end{align*}%
	where we used again Lemma \ref{key-estimate1}. The dominated convergence
	theorem implies that $f_{N}$ approximate $f$ in $\dot{F}_{p,q}(\mathbb{R}%
	^{n},\{t_{k}\})$. But $f_{N}$, $N\in \mathbb{N}$ is not necessary an element
	of $\mathcal{S}_{\infty }(\mathbb{R}^{n})$, so we need to approximate $f_{N}$
	in $\mathcal{S}_{\infty }\left( \mathbb{R}^{n}\right) $. Let $\omega \in 
	\mathcal{S}\left( \mathbb{R}^{n}\right) $ with $\omega (0)=1$ and supp($%
	\mathcal{F}(\omega) )\subset \{\xi :|\xi |\leq 1\}$. Put 
	\begin{equation*}
	f_{N,\delta }:=f_{N}\omega (\delta \cdot ),\quad 0<\delta <1.
	\end{equation*}%
	We have $f_{N,\delta }\in \mathcal{S}_{\infty }\left( \mathbb{R}^{n}\right) $
	see \cite[Lemma 5.3]{YY1}, and 
	\begin{equation*}
	f_{N}-f_{N,\delta }=\sum\limits_{k=-N}^{N}(\tilde{\psi}_{k}\ast \varphi
	_{k}\ast f)(1-\omega (\delta \cdot )).
	\end{equation*}%
	After simple calculation, we obtain 
	\begin{equation*}
	\varphi _{j}\ast \lbrack (\tilde{\psi}_{k}\ast \varphi _{k}\ast f)(\omega
	(\delta \cdot ))](x)=\int_{\mathbb{R}^{n}}\varphi _{k}\ast f(y)\varphi
	_{j}\ast (\tilde{\psi}_{k}\omega (\delta (\cdot +y))(x-y)dy,\quad x\in 
	\mathbb{R}^{n},
	\end{equation*}%
	which together with the fact that 
	\begin{equation*}
	\text{supp}(\mathcal{F}(\tilde{\psi}_{k}\omega (\delta (\cdot +y)))\subset
	\{\xi :2^{k-2}\leq |\xi |\leq 2^{k+1}\},\quad y\in \mathbb{R}^{n},|k|\leq N
	\end{equation*}%
	if $0<\delta <2^{-N-3}$\ yield that%
	\begin{equation*}
	\varphi _{j}\ast \lbrack (\tilde{\psi}_{k}\ast \varphi _{k}\ast f)(\omega
	(\delta \cdot ))]=0\quad \text{if}\quad |j-k|\geq 2.
	\end{equation*}%
	Therefore, we obtain that%
	\begin{equation*}
	\big\|f_{N}-f_{N,\delta }|\dot{F}_{p,q}(\mathbb{R}^{n},\{t_{k}\})\big\|
	\end{equation*}%
	can be estimated by%
	\begin{align*}
	& \Big\|\Big(\sum\limits_{|k|\leq N+2}\big|t_{k}(\varphi _{k}\ast
	\sum\limits_{i=-2}^{2}[(\tilde{\psi}_{k+i}\ast \varphi _{k+i}\ast
	f)(1-\omega (\delta \cdot ))])\big|^{q}\Big)^{1/q}|L_{p}(\mathbb{R}^{n})%
	\Big\| \\
	& \lesssim \sum\limits_{i=-2}^{2}\Big\|\Big(\sum\limits_{|k|\leq N+2}\big|%
	t_{k+i}\big((\tilde{\psi}_{k+i}\ast \varphi _{k+i}\ast f)(1-\omega (\delta
	\cdot ))\big|^{q}\Big)^{1/q}|L_{p}(\mathbb{R}^{n})\Big\|.
	\end{align*}%
	Again, by Lebesgue's dominated convergence theorem $f_{N,\delta }$
	approximate $f_{N}$ in the spaces $\dot{F}_{p,q}(\mathbb{R}^{n},\{t_{k}\})$.
	This prove that $\mathcal{S}_{\infty }(\mathbb{R}^{n})$ is dense in $\dot{F}%
	_{p,q}(\mathbb{R}^{n},\{t_{k}\})$.
	
	The proof of (ii) follows by the embedding 
	\begin{equation*}
	\dot{F}_{p,q}(\mathbb{R}^{n},\{t_{k}\})\hookrightarrow \dot{B}_{p,\max
		(p,q)}(\mathbb{R}^{n},\{t_{k}\})\hookrightarrow \mathcal{S}_{\infty
	}^{\prime }(\mathbb{R}^{n}),
	\end{equation*}%
	where the second embedding is proved in \cite{D20}.
\end{proof}

\subsection{Embeddings}

For our spaces introduced above we want to show Sobolev embedding theorems.
We recall that a quasi-Banach space $A_{1}$ is continuously embedded in
another quasi-Banach space $A_{2}$, $A_{1}\hookrightarrow A_{2}$, if $%
A_{1}\subset A_{2}$ and there is a $c>0$ such that 
\begin{equation*}
\big\|f|A_{2}\big\|\leq c\big\|f|A_{1}\big\|
\end{equation*}
for all $f\in A_{1}$. We begin with the following elementary embeddings.

\begin{thm}
	\label{elem-embedding}Let $0<p<\infty $ and $0<q\leq r<\infty $. Let $%
	\{t_{k}\}\in \dot{X}_{\alpha ,\sigma ,p}$ be a $p$-admissible weight
	sequence with $\sigma =(\sigma _{1}=\theta \left( \frac{p}{\theta }\right)
	^{\prime },\sigma _{2}\geq p)$. We have 
	\begin{equation*}
	\dot{F}_{p,q}(\mathbb{R}^{n},\{t_{k}\})\hookrightarrow \dot{F}_{p,r}(\mathbb{%
		R}^{n},\{t_{k}\}).
	\end{equation*}
\end{thm}

\begin{proof}
	It is a ready consequence of the embeddings between Lebesgue sequence spaces.
\end{proof}

The main result of this subsection is the following Sobolev-type embedding.
In the classical setting this was done in \cite{J77} and \cite{T1}. We set%
\begin{equation*}
w_{k,Q}(p_{1})=\Big(\int_{Q}w_{k}^{p_{1}}(x)dx\Big)^{1/p_{1}}\quad \text{%
	and\quad }t_{k,Q}(p_{0})=\Big(\int_{Q}t_{k}^{p_{0}}(x)dx\Big)^{1/p_{0}},
\end{equation*}%
where $Q\in \mathcal{Q}$ with $\ell (Q)=2^{-k}, k \in \mathbb{Z}.$

\begin{thm}
	\label{Sobolev-embedding-sequence}Let $0<\theta \leq p_{0}<p_{1}<\infty $
	and $0<q,r<\infty $. Let\ $\{t_{k}\}$ be a $p_{0}$-admissible weight
	sequence satisfying $\mathrm{\eqref{Asum1}}$\ with $p=p_{0}$, $\sigma
	_{1}=\theta \left( p_{0}/\theta \right) ^{\prime }$ and $j=k$. Let\ $%
	\{w_{k}\}$ be a $p_{1}$-admissible weight sequence satisfying $\mathrm{%
		\eqref{Asum1}}$\ with $p=p_{1}$, $\sigma _{1}=\theta \left( p_{1}/\theta
	\right) ^{\prime }$ and $j=k$. If $w_{k,Q}(p_{1})\lesssim t_{k,Q}(p_{0})$
	for all $Q\in \mathcal{Q}$ with $\ell (Q)=2^{-k}, k \in \mathbb{Z}$, then we have%
	\begin{equation*}
	\dot{f}_{p_{0},q}(\mathbb{R}^{n},\{t_{k}\})\hookrightarrow \dot{f}_{p_{1},r}(%
	\mathbb{R}^{n},\{w_{k}\}).
	\end{equation*}
\end{thm}

\begin{proof}
	Let $f\in \dot{f}_{p_{0},q}(\mathbb{R}^{n},\{t_{k}\})$. Without loss of
	generality, we may assume that 
	\begin{equation*}
	\big\|\lambda |\dot{f}_{p_{0},q}(\mathbb{R}^{n},\{t_{k}\})\big\|=1.
	\end{equation*}%
	We set%
	\begin{equation*}
	f_{k}(x)=\sum\limits_{m\in \mathbb{Z}%
		^{n}}2^{kn/2}|Q_{k,m}|^{-1/p_{1}}t_{k,m}(p_{0})|\lambda _{k,m}|\chi
	_{k,m}(x),\quad x\in \mathbb{R}^{n},k\in \mathbb{Z}.
	\end{equation*}%
	Using Proposition {\ref{Equi-norm1} and the fact that}%
	\begin{equation*}
	w_{k,m}(p_{1})\lesssim t_{k,m}(p_{0}),\quad k\in \mathbb{Z},m\in \mathbb{Z}%
	^{n},
	\end{equation*}%
	we obtain%
	\begin{equation*}
	\big\|\lambda |\dot{f}_{p_{1},r}(\mathbb{R}^{n},\{w_{k}\})\big\|\lesssim %
	\Big\|\Big(\sum\limits_{k=-\infty }^{\infty }f_{k}^{r}\Big)^{1/r}|L_{p_{1}}(%
	\mathbb{R}^{n})\Big\|.
	\end{equation*}%
	Now we prove our embedding. Let $K\in \mathbb{Z}$. By Lemma {\ref{Lamda-est},%
	} 
	\begin{equation*}
	|\lambda _{k,m}|\lesssim 2^{-kn/2}t_{k,m}^{-1}(p_{0})
	\end{equation*}%
	for any $k\in \mathbb{Z}$ and $m\in \mathbb{Z}^{n}$. Therefore,%
	\begin{equation}
	\sum_{k=-\infty }^{K}f_{k}^{r}(x)\lesssim \sum_{k=-\infty }^{K}2^{knr/p_{1}}\leq C2^{Knr/p_{1}}.  \label{est-sum1}
	\end{equation}%
	On the other hand it follows that%
	\begin{align}
	& \sum_{k=K+1}^{\infty }f_{k}^{r}(x)  \notag \\
	& \lesssim \sup_{k\in \mathbb{Z}}\Big(\sum\limits_{m\in \mathbb{Z}%
		^{n}}2^{kn/2}|Q_{k,m}|^{-1/p_{0}}t_{k,m}(p_{0})|\lambda _{k,m}|\chi _{k,m}(x)%
	\Big)^{r}\sum_{k=K+1}^{\infty }2^{knr(1/p_{1}-1/p_{0})}  \notag \\
	& \lesssim \sup_{k\in \mathbb{Z}}\Big(\sum\limits_{m\in \mathbb{Z}%
		^{n}}2^{kn/2}|Q_{k,m}|^{-1/p_{0}}t_{k,m}(p_{0})|\lambda _{k,m}|\chi _{k,m}(x)%
	\Big)^{r}2^{Knr(1/p_{1}-1/p_{0})}.  \label{est-sum2}
	\end{align}%
	The identity%
	\begin{equation*}
	\big\|g|L_{p_{1}}(\mathbb{R}^{n})\big\|^{p_{1}}=p_{1}\int_{0}^{\infty
	}y^{p_{1}-1}\left\vert \left\{ x\in \mathbb{R}^{n}:\left\vert g\left(
	x\right) \right\vert >y\right\} \right\vert dy,
	\end{equation*}%
	justifies the estimate 
	\begin{equation*}
	\big\|\lambda |\dot{f}_{p_{1},r}(\mathbb{R}^{n},\{w_{k}\})\big\|%
	^{p_{1}}\lesssim \int_{0}^{\infty }y^{p_{1}-1}\Big|\Big\{x\in \mathbb{R}^{n}:%
	\Big(\sum\limits_{k=-\infty }^{\infty }f_{k}^{r}(x)\Big)^{1/r}>y\Big\}\Big|%
	dy.
	\end{equation*}%
	We use $\mathrm{\eqref{est-sum1}}$ with $K$ the largest integer such that 
	\begin{equation*}
	C2^{Knr/p_{1}}\leq \frac{y^{r}}{2}.
	\end{equation*}
	Since $2^{Kn(1/p_{0}-1/p_{1})}y\geq c$ $y^{p_{1}/p_{0}}$,\ using $\mathrm{%
		\eqref{est-sum2}}$, we obtain%
	\begin{equation*}
	\Big|\Big\{x\in \mathbb{R}^{n}:\sum_{k=-\infty }^{\infty }f_{k}^{r}(x)>y^{r}%
	\Big\}\Big|\lesssim \Big|\Big\{x\in \mathbb{R}^{n}:\sum_{k=K+1}^{\infty
	}f_{k}^{r}(x)>\frac{y^{r}}{2}\Big\}\Big|,
	\end{equation*}%
	and hence%
	\begin{equation*}
	\Big|\Big\{x\in \mathbb{R}^{n}:\sum_{k=-\infty }^{\infty }f_{k}^{r}(x)>y^{r}%
	\Big\}\Big|
	\end{equation*}%
	does not exceed%
	\begin{equation*}
	c\Big|\Big\{x\in \mathbb{R}^{n}:\sup_{k\in \mathbb{Z}}\big(%
	2^{kn(1/p_{0}-1/p_{1})}f_{k}(x)\big)>c\text{ }2^{Kn(1/p_{0}-1/p_{1})}y\Big\}%
	\Big|.
	\end{equation*}%
	Therefore%
	\begin{align*}
	& \big\|\lambda |\dot{f}_{p_{1},r}(\mathbb{R}^{n},\{w_{k}\})\big\|^{p_{1}} \\
	& \lesssim \int_{0}^{\infty }y^{p_{1}-1}\Big|\Big\{x\in \mathbb{R}%
	^{n}:\sup_{k\in \mathbb{Z}}\big(2^{kn(1/p_{0}-1/p_{1})}f_{k}(x)\big)>c\text{ 
	}y^{p_{1}/p_{0}}\Big\}\Big|dy.
	\end{align*}%
	After a simple change of variable, we estimate the last term by%
	\begin{align*}
	&c\int_{0}^{\infty }y^{p_{0}-1}\Big|\Big\{x\in \mathbb{R}^{n}:\sup_{k\in 
		\mathbb{Z}}\big(2^{kn(1/p_{0}-1/p_{1})}f_{k}(x)\big)>y\Big\}\Big|dy\\
	&	\lesssim \big\|\lambda |\dot{f}_{p_{0},\infty }(\mathbb{R}^{n},\{t_{k}\})\big\|%
	^{p_{0}},
	\end{align*}%
	where we have used again Proposition {\ref{Equi-norm1}. }Hence the theorem
	is proved.
\end{proof}

From Theorems {\ref{phi-tran} and \ref{Sobolev-embedding-sequence}}, we
infer the following Sobolev-type embedding for $\dot{F}_{p,q}(\mathbb{R}%
^{n},\{t_{k}\})$.

\begin{thm}
	\label{Sobolev-embedding}Let $0<\theta \leq p_{0}<p_{1}<\infty $ and $%
	0<q,r<\infty $. Let $\{t_{k}\}\in \dot{X}_{\alpha _{0},\sigma ,p_{0}}$ be a $%
	p_{0}$-admissible weight sequence with $\sigma =(\sigma _{1}=\theta \left(
	p_{0}/\theta \right) ^{\prime },\sigma _{2}\geq p_{0})$ and $\alpha
	_{0}=(\alpha _{1,0},\alpha _{2,0})\in \mathbb{R}^{2}$. Let\ $\{w_{k}\}$ $\in 
	\dot{X}_{\alpha _{1},\sigma ,p_{1}}$ be a $p_{1}$-admissible weight sequence
	with $\sigma =(\sigma _{1}=\theta \left( p_{1}/\theta \right) ^{\prime
	},\sigma _{2}\geq p_{1})$ and $\alpha _{1}=(\alpha _{1,1},\alpha _{2,1})\in 
	\mathbb{R}^{2}$. Then%
	\begin{equation*}
	\dot{F}_{p_{0},q}(\mathbb{R}^{n},\{t_{k}\})\hookrightarrow \dot{F}_{p_{1},r}(%
	\mathbb{R}^{n},\{w_{k}\}),
	\end{equation*}%
	hold if%
	\begin{equation*}
	w_{k,Q}(p_{1})\lesssim t_{k,Q}(p_{0})
	\end{equation*}%
	for all $Q\in \mathcal{Q}$ and all $k\in \mathbb{Z}$.
\end{thm}

\subsection{Atomic and molecular decompositions}

We will use the notation of \cite{FJ90}. We shall say that an operator $A$
is associated with the matrix $\{a_{Q_{k,m}P_{v,h}}\}_{k,v\in \mathbb{Z}%
	,m,h\in \mathbb{Z}^{n}}$, if for all sequences $\lambda =\{\lambda
_{k,m}\}_{k\in \mathbb{Z},m\in \mathbb{Z}^{n}}\subset \mathbb{C}$,%
\begin{equation*}
A\lambda =\{(A\lambda )_{k,m}\}_{k\in \mathbb{Z},m\in \mathbb{Z}^{n}}=\Big\{%
\sum_{v=-\infty }^{\infty }\sum_{h\in \mathbb{Z}^{n}}a_{Q_{k,m}P_{v,h}}%
\lambda _{v,h}\Big\}_{k\in \mathbb{Z},m\in \mathbb{Z}^{n}}.
\end{equation*}%
We will use the notation%
\begin{equation*}
J=\frac{n}{\min (1,p,q)}.
\end{equation*}%
We say that $A$, with associated matrix $\{a_{Q_{k,m}P_{v,h}}\}_{k,v\in 
	\mathbb{Z},m,h\in \mathbb{Z}^{n}}$, is almost diagonal on $\dot{f}_{p,q}(%
\mathbb{R}^{n},\{t_{k}\})$ if there exists $\varepsilon >0$ such that%
\begin{equation*}
\sup_{k,v\in \mathbb{Z},m,h\in \mathbb{Z}^{n}}\frac{|a_{Q_{k,m}P_{v,h}}|}{%
	\omega _{Q_{k,m}P_{v,h}}(\varepsilon )}<\infty ,
\end{equation*}%
where

\begin{align}
&\omega _{Q_{k,m}P_{v,h}}(\varepsilon ) \notag \\
&=\Big(1+\frac{|x_{Q_{k,m}}-x_{P_{v,h}}|%
}{\max \big(2^{-k},2^{-v}\big)}\Big)^{-J-\varepsilon }\left\{ 
\begin{array}{ccc}
2^{(v-k)(\alpha _{2}+(n+\varepsilon )/2)}, & \text{if} & v\leq k, \\ 
2^{(v-k)(\alpha _{1}-(n+\varepsilon )/2-J+n)}, & \text{if} & v>k.%
\end{array}%
\right.  \label{omega-assumption}
\end{align}

Using Lemma {\ref{key-estimate1.1} }the following theorem is a
generalization of {\cite[Theorem 3.3]{FJ90}.}

\begin{thm}
	\label{almost-diag-est}Let $\alpha _{1},\alpha _{2}\in \mathbb{R}$, $%
	0<\theta \leq p<\infty $ and $0<q<\infty $. Let $\{t_{k}\}_{k}\in \dot{X}%
	_{\alpha ,\sigma ,p}$ be a $p$-admissible weight sequence with $\sigma
	_{1}=\theta \left( p/\theta \right) ^{\prime }$ and\ $\sigma _{2}\geq p$.
	Any almost diagonal operator $A$ on $\dot{f}_{p,q}(\mathbb{R}^{n},\{t_{k}\})$
	is bounded.
\end{thm}

\begin{proof}
	We write $A\equiv A_{0}+A_{1}$ with\ 
	\begin{equation*}
	(A_{0}\lambda )_{k,m}=\sum_{v=-\infty }^{k}\sum_{h\in \mathbb{Z}%
		^{n}}a_{Q_{k,m}P_{v,h}}\lambda _{v,h},\quad k\in \mathbb{Z},m\in \mathbb{Z}%
	^{n}
	\end{equation*}%
	and%
	\begin{equation*}
	(A_{1}\lambda )_{k,m}=\sum_{v=k+1}^{\infty }\sum_{h\in \mathbb{Z}%
		^{n}}a_{Q_{k,m}P_{v,h}}\lambda _{v,h},\quad k\in \mathbb{Z},m\in \mathbb{Z}%
	^{n}.
	\end{equation*}%
	
	\textit{Estimate of }$A_{0}$. From $\mathrm{\eqref{omega-assumption}}$,\ we
	obtain%
	\begin{align*}
	\big|(A_{0}\lambda )_{k,m}\big|& \leq \sum_{v=-\infty }^{k}\sum_{h\in 
		\mathbb{Z}^{n}}2^{(v-k)(\alpha _{2}+(n+\varepsilon )/2)}\frac{|\lambda
		_{v,h}|}{\big(1+2^{v}|x_{k,m}-x_{v,h}|\big)^{J+\varepsilon }} \\
	& =\sum_{v=-\infty }^{k}2^{(v-k)(\alpha _{2}+(n+\varepsilon )/2)}S_{k,v,m}.
	\end{align*}%
	{For each }$j\in \mathbb{N},k\in \mathbb{Z}$\ and $m\in \mathbb{Z}^{n}${\ we
		define }%
	\begin{equation*}
	{\Omega _{j,k,m}=\{h\in \mathbb{Z}^{n}:2^{j-1}<2^{v}|x_{k,m}-x_{v,h}|\leq
		2^{j}\}}
	\end{equation*}%
	{and }%
	\begin{equation*}
	{\Omega _{0,k,m}=\{h\in \mathbb{Z}^{n}:2^{v}|x_{k,m}-x_{v,h}|\leq 1\}.}
	\end{equation*}%
	Let $n/(J+\frac{\varepsilon }{2})<\tau <\min (1,p,q)$. We rewrite $S_{k,v,m}$
	as follows%
	\begin{align*}
	S_{k,v,m}& =\sum\limits_{j=0}^{\infty }\sum\limits_{h\in {\Omega _{j,k,m}}}%
	\frac{|\lambda _{v,h}|}{\big(1+2^{v}|x_{k,m}-x_{v,h}|\big)^{J+\varepsilon }}
	\\
	& \leq \sum\limits_{j=0}^{\infty }2^{-(J+\varepsilon )j}\sum\limits_{h\in {%
			\Omega _{j,k,m}}}|\lambda _{v,h}|.
	\end{align*}%
	By the embedding $\ell _{\tau }\hookrightarrow \ell _{1}$ we deduce that%
	\begin{align*}
	S_{k,v,m} &\leq \sum\limits_{j=0}^{\infty }2^{-(J+\varepsilon )j}\big(%
	\sum\limits_{h\in {\Omega _{j,k,m}}}|\lambda _{v,h}|^{\tau }\big)^{1/\tau } 
	\notag \\
	&=\sum\limits_{j=0}^{\infty }2^{(n/\tau -J-\varepsilon )j}\Big(%
	2^{(v-j)n}\int\limits_{\cup _{z\in {\Omega _{j,k,m}}}Q_{v,z}}\sum\limits_{h%
		\in {\Omega _{j,k,m}}}|\lambda _{v,h}|^{\tau }\chi _{v,h}(y)dy\Big)^{1/\tau
	}.  \label{est-Sv,h}
	\end{align*}%
	Let $y\in \cup _{z\in {\Omega _{j,k,m}}}Q_{v,z}$ and $x\in Q_{k,m}$. It
	follows that $y\in Q_{v,z}$ for some $z\in {\Omega _{j,k,m}}$ and ${2^{j-1}<2%
	}^{v}{|2}^{-k}m{-2}^{-v}{z|\leq 2^{j}}$. From this we obtain that%
	\begin{align*}
	\left\vert y-x\right\vert & \leq \left\vert y-{2}^{-k}m\right\vert
	+\left\vert x-{2}^{-k}m\right\vert \\
	& \lesssim 2^{-v}+2^{j-v}+2^{-k} \\
	& \leq 2^{j-v+\delta _{n}},\quad \delta _{n}\in \mathbb{N},
	\end{align*}%
	which implies that $y$ is located in the ball $B(x,2^{j-v+\delta _{n}})$.
	Consequently%
	\begin{equation*}
	S_{k,v,m}\lesssim \mathcal{M}_{\tau }\big(\sum\limits_{h\in \mathbb{Z}%
		^{n}}\lambda _{v,h}\chi _{v,h}\big)(x)
	\end{equation*}%
	for any $x\in Q_{k,m}$ and any $k\leq v$. Applying Lemma \ref%
	{key-estimate1.1}, we obtain that%
	\begin{equation*}
	{{\big\|}A_{0}\lambda |\dot{f}_{p,q}(\mathbb{R}^{n},\{t_{k}\}){\big\|}}
	\end{equation*}%
	is bounded by 
	\begin{equation*}
	c{{\big\|}\lambda |\dot{f}_{p,q}(\mathbb{R}^{n},\{t_{k}\}){\big\|}.}
	\end{equation*}%
	
	\textit{Estimate of }$A_{1}$. Again from $\mathrm{\eqref{omega-assumption}}$%
	,\ we see that%
	\begin{align*}
	\big|(A_{1}\lambda )_{v,h}\big|& \leq \sum_{v=k+1}^{\infty }\sum_{h\in 
		\mathbb{Z}^{n}}2^{(v-k)(\alpha _{1}-\varepsilon /2-J+n/2)}\frac{|\lambda
		_{v,h}|}{\big(1+2^{k}|x_{k,m}-x_{v,h}|\big)^{J+\varepsilon }} \\
	& =\sum_{v=k+1}^{\infty }2^{(v-k)(\alpha _{1}-\varepsilon
		/2-J+n/2)}T_{k,v,m}.
	\end{align*}%
	We proceed as in the estimate of $A_{0}$ we can prove that%
	\begin{equation*}
	T_{k,v,m}\leq c2^{(v-k)n/\tau }\mathcal{M}_{\tau }\big(\sum\limits_{h\in 
		\mathbb{Z}^{n}}\lambda _{v,h}\chi _{v,h}\big)(x),\quad v>k,x\in Q_{k,m},
	\end{equation*}%
	where $n/(J+\frac{\varepsilon }{2})<\tau <\min (1,p,q)$ and the positive
	constant $c$ is independent of $v$, $k$ and $m$. Again applying Lemma \ref%
	{key-estimate1.1} we obtain%
	\begin{equation*}
	{{\big\|}A_{1}\lambda |\dot{f}_{p,q}(\mathbb{R}^{n},\{t_{k}\}){\big\|}}
	\end{equation*}%
	is bounded by 
	\begin{equation*}
	c{{\big\|}\lambda |\dot{f}_{p,q}(\mathbb{R}^{n},\{t_{k}\}){\big\|}.}
	\end{equation*}%
	Hence the theorem is proved.
\end{proof}

\begin{defn}
	\label{Atom-Def}Let\ $\alpha _{1},\alpha _{2}\in \mathbb{R},0<p<\infty $ and 
	$0<q<\infty $. Let $\{t_{k}\}$ be a $p$-admissible weight sequence. Let $%
	N=\max \{J-n-\alpha _{1},-1\}$ and $\alpha _{2}^{\ast }=\alpha _{2}-\lfloor
	\alpha _{2}\rfloor $.\newline
	$\mathrm{(i)}$\ Let $k\in \mathbb{Z}$ and $m\in \mathbb{Z}^{n}$. A function $%
	\varrho _{Q_{k,m}}$ is called an homogeneous smooth synthesis molecule for $%
	\dot{F}_{p,q}(\mathbb{R}^{n},\{t_{k}\})$\ supported near $Q_{k,m}$ if there
	exist a real number $\delta \in (\alpha _{2}^{\ast },1]$ and a real number $%
	M\in (J,\infty )$ such that%
	\begin{equation}
	\int_{\mathbb{R}^{n}}x^{\beta }\varrho _{Q_{k,m}}(x)dx=0\text{\quad if\quad }%
	0\leq |\beta |\leq N,  \label{mom-cond}
	\end{equation}%
	\begin{equation}
	|\varrho _{Q_{k,m}}(x)|\leq 2^{\frac{kn}{2}}(1+2^{k}|x-x_{Q_{k,m}}|)^{-\max
		(M,M-\alpha _{1})},  \label{cond1}
	\end{equation}%
	\begin{equation}
	|\partial ^{\beta }\varrho _{Q_{k,m}}(x)|\leq 2^{k(|\beta |+\frac{1}{2}%
		)}(1+2^{k}|x-x_{Q_{k,m}}|)^{-M}\quad \text{if}\quad |\beta |\leq \lfloor
	\alpha _{2}\rfloor  \label{cond2}
	\end{equation}%
	and%
	\begin{align}
	&|\partial ^{\beta }\varrho _{Q_{k,m}}(x)-\partial ^{\beta }\varrho
	_{Q_{k,m}}(y)|  \label{cond3} \\
	&\leq 2^{k(|\beta |+\frac{1}{2}+\delta )}|x-y|^{\delta }\sup_{|z|\leq
		|x-y|}(1+2^{k}|x-z-x_{Q_{k,m}}|)^{-M}\text{\quad if\quad }|\beta |=\lfloor
	\alpha _{2}\rfloor .  \notag
	\end{align}%
	A collection $\{\varrho _{Q_{k,m}}\}_{k\in \mathbb{Z},m\in \mathbb{Z}^{n}}$
	is called a family of homogeneous smooth synthesis molecules for $\dot{F}%
	_{p,q}(\mathbb{R}^{n},\{t_{k}\})$, if each $\varrho _{Q_{k,m}}$, $k\in 
	\mathbb{Z},m\in \mathbb{Z}^{n}$, is an homogeneous smooth synthesis molecule
	for $\dot{F}_{p,q}(\mathbb{R}^{n},\{t_{k}\})$ supported near $Q_{k,m}$. 
	\newline
	$\mathrm{(ii)}$\ Let $k\in \mathbb{Z}$ and $m\in \mathbb{Z}^{n}$. A function 
	$b_{Q_{k,m}}$ is called an homogeneous smooth analysis molecule for $\dot{F}%
	_{p,q}(\mathbb{R}^{n},\{t_{k}\})$ supported near $Q_{k,m}$ if there exist a $%
	\kappa \in ((J-\alpha _{2})^{\ast },1]$ and an $M\in (J,\infty )$ such that
	\begin{equation}
	\int_{\mathbb{R}^{n}}x^{\beta }b_{Q_{k,m}}(x)dx=0\text{\quad if\quad }0\leq
	|\beta |\leq \left\lfloor \alpha _{2}\right\rfloor ,  \label{mom-cond2}
	\end{equation}%
	\begin{equation}
	|b_{Q_{k,m}}(x)|\leq 2^{\frac{kn}{2}}(1+2^{k}|x-x_{Q_{k,m}}|)^{-\max
		(M,M+n+\alpha _{2}-J)},  \label{cond1.1}
	\end{equation}%
	\begin{equation}
	|\partial ^{\beta }b_{Q_{k,m}}(x)|\leq 2^{k(|\beta |+\frac{n}{2}%
		)}(1+2^{k}|x-x_{Q_{k,m}}|)^{-M}\quad \text{if}\quad |\beta |\leq N
	\label{cond1.2}
	\end{equation}%
	and%
	\begin{align}
	&|\partial ^{\beta }b_{Q_{k,m}}(x)-\partial ^{\beta }b_{Q_{k,m}}(y)|
	\label{cond1.3} \\
	&\leq 2^{k(|\beta |+\frac{n}{2}+\kappa )}|x-y|^{\kappa }\sup_{|z|\leq
		|x-y|}(1+2^{k}|x-z-x_{Q_{k,m}}|)^{-M}\text{\quad if\quad }|\beta |=N.  \notag
	\end{align}%
	A collection $\{b_{Q_{k,m}}\}_{k\in \mathbb{Z},m\in \mathbb{Z}^{n}}$ is
	called a family of homogeneous smooth analysis molecules for $\dot{F}_{p,q}(%
	\mathbb{R}^{n},\{t_{k}\})$, if each $b_{Q_{k,m}}$, $k\in \mathbb{Z},m\in 
	\mathbb{Z}^{n}$, is an homogeneous smooth synthesis molecule for $\dot{F}%
	_{p,q}(\mathbb{R}^{n},\{t_{k}\})$ supported near $Q_{k,m}$.
\end{defn}

We will use the notation $\{b_{k,m}\}_{k\in \mathbb{Z},m\in \mathbb{Z}^{n}}$
instead of $\{b_{Q_{k,m}}\}_{k\in \mathbb{Z},m\in \mathbb{Z}^{n}}$. The
proof of the following lemma is given in {\cite{D20}.}

\begin{lem}
	\label{matrix-est}Let\ $\alpha _{1},\alpha _{2},J,M,N,\delta ,\kappa ,p\ $%
	and $q$ be as in Definition {\ref{Atom-Def}}. Let $\{t_{k}\}$ be a $p$%
	-admissible weight sequence. Suppose $\{\varrho _{v,h}\}_{v\in \mathbb{Z}%
		,h\in \mathbb{Z}^{n}}$ is a family of smooth synthesis molecules for $\dot{F}%
	_{p,q}(\mathbb{R}^{n},\{t_{k}\})$ and $\{b_{k,m}\}_{k\in \mathbb{Z},m\in 
		\mathbb{Z}^{n}}$\ is a family of homogeneous smooth analysis molecules for\ $%
	\dot{F}_{p,q}(\mathbb{R}^{n},\{t_{k}\})$. Then there exist a positive real
	number $\varepsilon _{1}$ and a positive constant $c$ such that%
	\begin{equation*}
	\left\vert \langle \varrho _{v,h},b_{k,m}\rangle \right\vert \leq c\text{ }%
	\omega _{Q_{k,m}P_{v,h}}(\varepsilon ),\quad k,v\in \mathbb{Z},h,m\in 
	\mathbb{Z}^{n}
	\end{equation*}%
	if $\varepsilon \leq \varepsilon _{1}$.
\end{lem}

As an immediate consequence, we have the following analogues of the
corresponding results on \cite[Corollary\ B.3]{FJ90}.

\begin{cor}
	Let\ $\alpha _{1},\alpha _{2},J,M,N,\delta ,\kappa ,p\ $and $q$ be as in
	Definition {\ref{Atom-Def}}. Let $\{t_{k}\}$ be a $p$-admissible weight
	sequence. Let $\Phi $ and $\varphi $ satisfy, respectively $\mathrm{%
		\eqref{Ass1}}$ and $\mathrm{\eqref{Ass2}}$.\newline
	$\mathrm{(i)}$\ If $\{\varrho _{k,m}\}_{k\in \mathbb{Z},m\in \mathbb{Z}^{n}}$
	is a family of homogeneous synthesis molecules for the Triebel-Lizorkin
	spaces $\dot{F}_{p,q}(\mathbb{R}^{n},\{t_{k}\})$, then the operator $A$ with
	matrix\ $a_{Q_{k,m}P_{v,h}}=\langle \varrho _{v,h},\varphi _{k,m}\rangle $, $%
	k,v\in \mathbb{Z},m,h\in \mathbb{Z}^{n}$, is almost diagonal.\newline
	$\mathrm{(ii)}$ If $\{b_{k,m}\}_{k\in \mathbb{Z},m\in \mathbb{Z}^{n}}$ is a
	family of homogeneous smooth analysis molecules for the Triebel-Lizorkin
	spaces $\dot{F}_{p,q}(\mathbb{R}^{n},\{t_{k}\})$, then \ the operator\ $A$,
	with matrix $a_{Q_{k,m}P_{v,h}}=\langle \varphi _{v,h},b_{Q_{k,m}}\rangle $, 
	$k,v\in \mathbb{Z},m,h\in \mathbb{Z}^{n}$, is almost diagonal.
\end{cor}

Let $f\in \dot{F}_{p,q}(\mathbb{R}^{n},\{t_{k}\})\ $and $\{b_{k,m}\}_{k\in 
	\mathbb{Z},m\in \mathbb{Z}^{n}}$ be a family of homogeneous\ smooth analysis
molecules. To prove that $\langle f,b_{Q_{k,m}}\rangle $, $k\in \mathbb{Z}%
,m\in \mathbb{Z}^{n}$, is well defined for all homogeneous smooth analysis
molecules for $\dot{F}_{p,q}(\mathbb{R}^{n},\{t_{k}\}$, we need the
following result, which proved in {\cite[Lemma 5.4]{BoHo06}}. Suppose that $%
\Phi $ is a smooth analysis (or synthesis) molecule supported near $Q\in 
\mathcal{Q}$ . Then there exists a sequence $\{\varphi _{k}\}_{k\in \mathbb{N%
}}\subset \mathcal{S(}\mathbb{R}^{n}\mathcal{)}$ and $c>0$ such that $%
c\varphi _{k}$ is a smooth analysis (or synthesis) molecule supported near $%
Q $ for every $k$,and $\varphi _{k}(x)\rightarrow \Phi (x)$ uniformly on $%
\mathbb{R}^{n}$ as $k\rightarrow \infty $.

Now we have the following smooth molecular characterization of the spaces $\dot{F}%
_{p,q}(\mathbb{R}^{n},\{t_{k}\})$. We refer the reader to {\cite{D20}} {for
	the corresponding result for Besov spaces.}

\begin{thm}
	\label{molecules-dec}Let $\alpha _{1}$, $\alpha _{2}\in \mathbb{R},0<\theta
	\leq p<\infty $ and $0<q<\infty $. Let $\{t_{k}\}_{k}\in \dot{X}_{\alpha
		,\sigma ,p}$ be a $p$-admissible weight sequence with $\sigma _{1}=\theta
	\left( p/\theta \right) ^{\prime }$ and\ $\sigma _{2}\geq p$. Let $%
	J,M,N,\delta $ and $\kappa $ be as in Definition {\ref{Atom-Def}}. \newline
	$\mathrm{(i)}$\ If $f=\sum_{v=-\infty }^{\infty }\sum_{h\in \mathbb{Z}%
		^{n}}\varrho _{v,h}\lambda _{v,h}$, where $\{\varrho _{v,h}\}_{v\in \mathbb{Z%
		},h\in \mathbb{Z}^{n}}$ is a family of homogeneous smooth synthesis
	molecules for $\dot{F}_{p,q}(\mathbb{R}^{n},\{t_{k}\})$, then for all $%
	\lambda \in \dot{f}_{p,q}(\mathbb{R}^{n},\{t_{k}\})$ 
	\begin{equation*}
	{{\big\|}f|\dot{F}_{p,q}(\mathbb{R}^{n},\{t_{k}\}){\big\|}\lesssim {\big\|}%
		\lambda |\dot{f}_{p,q}(\mathbb{R}^{n},\{t_{k}\}){\big\|}.}
	\end{equation*}%
	$\mathrm{(ii)}$\ Let $\{b_{k,m}\}_{k\in \mathbb{Z},m\in \mathbb{Z}^{n}}$ be
	a family of homogeneous\ smooth analysis molecules.\ Then for all\ $f\in 
	\dot{F}_{p,q}(\mathbb{R}^{n},\{t_{k}\})$%
	\begin{equation*}
	{{\big\|}\{\langle f,b_{k,m}\rangle \}_{k\in \mathbb{Z},m\in \mathbb{Z}^{n}}|%
		\dot{f}_{p,q}(\mathbb{R}^{n},\{t_{k}\}){\big\|}\lesssim {\big\|}f|\dot{F}%
		_{p,q}(\mathbb{R}^{n},\{t_{k}\})\big\|.}
	\end{equation*}
\end{thm}

\begin{proof}
	The proof is a slight variant of \cite{FJ90} and \cite{D20}{. }For the
	convenience of the reader, we give some details.
	
	\textit{Step 1. Proof of }\textrm{(i)}. By\ $\mathrm{\eqref{proc2}}$ we can
	write%
	\begin{equation*}
	\varrho _{v,h}=\sum_{k=-\infty }^{\infty }2^{-kn}\sum_{m\in \mathbb{Z}^{n}}%
	\widetilde{\varphi }_{k}\ast \varrho _{v,h}(2^{-k}m)\psi _{k}(\cdot -2^{-k}m)
	\end{equation*}%
	for\ any\ $v\in \mathbb{Z},h\in \mathbb{Z}^{n}$.\ Therefore,%
	\begin{equation*}
	f=\sum_{k=-\infty }^{\infty }\sum_{m\in \mathbb{Z}^{n}}S_{k,m}\psi
	_{k,m}=T_{\psi }S,
	\end{equation*}%
	where $S=\{S_{k,m}\}_{k\in \mathbb{Z},m\in \mathbb{Z}^{n}}$, with%
	\begin{equation*}
	S_{k,m}=2^{-kn/2}\sum_{v=-\infty }^{\infty }\sum_{h\in \mathbb{Z}^{n}}%
	\widetilde{\varphi }_{k}\ast \varrho _{v,h}(2^{-k}m)\lambda _{v,h}.
	\end{equation*}%
	From Theorem {\ref{phi-tran}, }we have%
	\begin{equation*}
	{{\big\|}f|\dot{F}_{p,q}(\mathbb{R}^{n},\{t_{k}\}){\big\|}={\big\|}T_{\psi
		}S|\dot{F}_{p,q}(\mathbb{R}^{n},\{t_{k}\})\big\|\lesssim {\big\|}S|\dot{f}%
		_{p,q}(\mathbb{R}^{n},\{t_{k}\})\big\|.}
	\end{equation*}%
	But%
	\begin{equation*}
	S_{k,m}=\sum_{v=-\infty }^{\infty }\sum_{h\in \mathbb{Z}%
		^{n}}a_{Q_{k,m}P_{v,h}}\lambda _{v,h},
	\end{equation*}%
	with 
	\begin{equation*}
	a_{Q_{k,m}P_{v,h}}=\langle \varrho _{v,h},\widetilde{\varphi }_{k,m}\rangle
	,\quad k,v\in \mathbb{Z},m,h\in \mathbb{Z}^{n}.
	\end{equation*}%
	Applying Lemma {\ref{matrix-est} and Theorem \ref{almost-diag-est} we find
		that}%
	\begin{equation*}
	{{\big\|}S|\dot{f}_{p,q}(\mathbb{R}^{n},\{t_{k}\}){\big\|}\lesssim {\big\|}%
		\lambda |\dot{f}_{p,q}(\mathbb{R}^{n},\{t_{k}\}){\big\|}.}
	\end{equation*}%
	
	\textit{Step 2. Proof of }$\mathit{\mathrm{(ii)}}$\textit{.} We have%
	\begin{align*}
	\langle f,b_{k,m}\rangle &=\sum_{v=-\infty }^{\infty }2^{-vn}\sum_{m\in 
		\mathbb{Z}^{n}}\langle \psi _{v}(\cdot -2^{-v}h),b_{k,m}\rangle \widetilde{%
		\varphi }_{v}\ast f(2^{-v}h) \\
	&=\sum_{v=-\infty }^{\infty }\sum_{m\in \mathbb{Z}^{n}}\langle \psi
	_{v,h},b_{k,m}\rangle \lambda _{v,h} \\
	&=\sum_{v=-\infty }^{\infty }\sum_{h\in \mathbb{Z}^{n}}a_{Q_{k,m}P_{v,h}}%
	\lambda _{v,h},
	\end{align*}%
	{\ where }%
	\begin{equation*}
	a_{Q_{k,m}P_{v,h}}=\langle \psi _{v,h},b_{k,m}\rangle ,\quad \lambda
	_{v,h}=2^{-vn/2}\widetilde{\varphi }_{v}\ast f(2^{-v}h).
	\end{equation*}%
	Again by Lemma {\ref{matrix-est} and Theorem \ref{almost-diag-est} we find that}%
	\begin{align*}
	{\big\|\{\langle f,b_{k,m}\rangle \}_{k\in \mathbb{Z},m\in \mathbb{Z}^{n}}|%
		\dot{f}_{p,q}(\mathbb{R}^{n},\{t_{k}\})\big\|} &{\lesssim }{\big\|\{\lambda
		_{v,h}\}_{v\in \mathbb{Z},h\in \mathbb{Z}^{n}}|\dot{f}_{p,q}(\mathbb{R}%
		^{n},\{t_{k}\})\big\|} \\
	&=c{\big\|\{(S_{\varphi })_{v,h}\}_{v\in \mathbb{Z},h\in \mathbb{Z}^{n}}|%
		\dot{f}_{p,q}(\mathbb{R}^{n},\{t_{k}\})\big\|.}
	\end{align*}%
	Applying Theorem {\ref{phi-tran} we find that}%
	\begin{equation*}
	{{\big\|}\{\langle f,b_{k,m}\rangle \}_{k\in \mathbb{Z},m\in \mathbb{Z}^{n}}|%
		\dot{f}_{p,q}(\mathbb{R}^{n},\{t_{k}\}){\big\|}\lesssim {\big\|}f|\dot{F}%
		_{p,q}(\mathbb{R}^{n},\{t_{k}\}){\big\|}.}
	\end{equation*}%
	The proof is complete.
\end{proof}

\begin{defn}
	Let $\alpha _{1},\alpha _{2}\in \mathbb{R},0<p<\infty ,0<q<\infty $ and\ $%
	N=\max \{J-n-\alpha _{1},-1\}$. Let $\{t_{k}\}$ be a $p$-admissible weight
	sequence. A function $a_{Q_{v,m}}$ is called an homogeneous smooth atom for $%
	\dot{F}_{p,q}(\mathbb{R}^{n},\{t_{k}\})$ supported near $Q_{k,m}$, $k\in 
	\mathbb{Z}$ and $m\in \mathbb{Z}^{n}$, if%
	\begin{equation}
	\mathrm{supp}(\text{ }a_{Q_{k,m}})\subseteq 3Q_{k,m}  \label{supp-cond}
	\end{equation}
	
	\begin{equation}
	|\partial ^{\beta }a_{Q_{k,m}}(x)|\leq 2^{kn(|\beta |+1/2)}\text{\quad
		if\quad }0\leq |\beta |\leq \max (0,1+\lfloor \alpha _{2}\rfloor ),\quad
	x\in \mathbb{R}^{n}  \label{diff-cond}
	\end{equation}%
	and if%
	\begin{equation}
	\int_{\mathbb{R}^{n}}x^{\beta }a_{Q_{k,m}}(x)dx=0\text{\quad if\quad }0\leq
	|\beta |\leq N\text{ and }k\in \mathbb{Z}.  \label{mom-cond1}
	\end{equation}
\end{defn}

A collection $\{a_{Q_{k,m}}\}_{k\in \mathbb{Z},m\in \mathbb{Z}^{n}}$\ is
called a family of homogeneous smooth atoms for $\dot{F}_{p,q}(\mathbb{R}%
^{n},\{t_{k}\})$, if each $a_{Q_{k,m}}$ is an homogeneous smooth atom for $%
\dot{F}_{p,q}(\mathbb{R}^{n},\{t_{k}\})$ supported near $Q_{v,m}$. We point
out that in the moment condition $\mathrm{\eqref{mom-cond1}}$ can be
strengthened into that%
\begin{equation*}
\int_{\mathbb{R}^{n}}x^{\beta }a_{Q_{k,m}}(x)dx=0\text{\quad if\quad }0\leq
|\beta |\leq \tilde{N}\text{ and }k\in \mathbb{Z}
\end{equation*}%
and the regularity condition $\mathrm{\eqref{diff-cond}}$ can be
strengthened into that%
\begin{equation*}
|\partial ^{\beta }a_{Q_{k,m}}(x)|\leq 2^{kn(|\beta |+1/2)}\text{\quad
	if\quad }0\leq |\beta |\leq \tilde{K},\quad x\in \mathbb{R}^{n},
\end{equation*}%
where $\tilde{K}$ and $\tilde{N}$ are arbitrary fixed integer satisfying $%
\tilde{K}\geq \max (0,1+\lfloor \alpha _{2}\rfloor )$ and $\tilde{N}\geq
\max \{J-n-\alpha _{1},-1\}$. If an atom $a$ is supported near $Q_{v,m}$,
then we denote it by $a_{v,m}$.\vskip5pt

Now we come to the atomic decomposition theorem, see {\cite{D20}} {for Besov
	spaces and the same arguments are true for Triebel-Lizorkin spaces.}

\begin{thm}
	\label{atomic-dec}Let $\alpha _{1}$, $\alpha _{2}\in \mathbb{R}$, $0<\theta
	\leq p<\infty $, $0<q<\infty $. Let $\{t_{k}\}_{k}\in \dot{X}_{\alpha
		,\sigma ,p}$ be a $p$-admissible weight sequence with $\sigma _{1}=\theta
	\left( p/\theta \right) ^{\prime }$ and\ $\sigma _{2}\geq p$. Then for each $%
	f\in \dot{F}_{p,q}(\mathbb{R}^{n},\{t_{k}\})$, there exist a family\ $%
	\{\varrho _{k,m}\}_{k\in \mathbb{Z},m\in \mathbb{Z}^{n}}$ of homogeneous
	smooth atoms for $\dot{F}_{p,q}(\mathbb{R}^{n},\{t_{k}\})$ and $\lambda
	=\{\lambda _{k,m}\}_{k\in \mathbb{Z},m\in \mathbb{Z}^{n}}\in {\dot{f}}_{p,q}(%
	\mathbb{R}^{n},\{t_{k}\})$ such that 
	\begin{equation*}
	f=\sum\limits_{k=-\infty }^{\infty }\sum\limits_{m\in \mathbb{Z}^{n}}\lambda
	_{k,m}\varrho _{k,m},\text{\quad converging in }\mathcal{S}_{\infty
	}^{\prime }(\mathbb{R}^{n})
	\end{equation*}%
	and%
	\begin{equation*}
	{{\big\|}\{\lambda _{k,m}\}_{k\in \mathbb{Z},m\in \mathbb{Z}^{n}}|\dot{f}%
		_{p,q}(\mathbb{R}^{n},\{t_{k}\}){\big\|}\lesssim {\big\|}f|\dot{F}_{p,q}(%
		\mathbb{R}^{n},\{t_{k}\})\big\|.}
	\end{equation*}%
	Conversely, for any family of homogeneous smooth atoms for $\dot{F}_{p,q}(%
	\mathbb{R}^{n},\{t_{k}\})$ and $\lambda =\{\lambda _{k,m}\}_{k\in \mathbb{Z}%
		,m\in \mathbb{Z}^{n}}\in {\dot{f}}_{p,q}(\mathbb{R}^{n},\{t_{k}\})$%
	\begin{equation*}
	{{\big\|}\sum\limits_{k=-\infty }^{\infty }\sum\limits_{m\in \mathbb{Z}%
			^{n}}\lambda _{k,m}\varrho _{k,m}|\dot{F}_{p,q}(\mathbb{R}^{n},\{t_{k}\})%
		\big\|\lesssim {\big\|}\{\lambda _{k,m}\}_{k\in \mathbb{Z},m\in \mathbb{Z}%
			^{n}}|\dot{f}_{p,q}(\mathbb{R}^{n},\{t_{k}\}){\big\|}.}
	\end{equation*}
\end{thm}

\begin{rem}
	$\mathrm{(i)}$ Further results, concerning, for instance, characterizations
	via oscillations, box spline and tensor-product B-spline representations\
	are given\ in {\cite{D7}.}\newline
	$\mathrm{(ii)}$ We mention that the techniques of {\cite{HN07} }are
	incapable of dealing with spaces of variable smoothness. Also our
	assumptions on the weight $\{t_{k}\}$\ play an exceptional role in the paper.%
	\newline
	$\mathrm{(iii)}$ We draw the reader's attention to paper {\cite{LSTDW10} }%
	where generalized Besov-type and Triebel-Lizorkin-type spaces are studied.
	They assumed that the weight sequence $\{t_{k}\}$ lies in some class which
	different from the class $\dot{X}_{\alpha ,\sigma ,p}.$
\end{rem}

\section{The non-homogeneous space\textbf{\ }$F_{p,q}(\mathbb{R}%
	^{n},\{t_{k}\}_{k\in \mathbb{N}_{0}})$}

In this section, we present the inhomogeneous version of our results given
above. Let $\Phi ,\psi ,\varphi $ and $\Psi $ satisfy 
\begin{equation}
\Phi ,\Psi ,\varphi ,\psi \in \mathcal{S}(\mathbb{R}^{n})  \label{Ass1.1}
\end{equation}%
\begin{equation}
\text{supp}(\mathcal{F}(\Phi ))\cup \text{supp}\mathcal{F}((\Psi ))\subset 
\overline{B(0,2)},\text{\quad }|\mathcal{F}(\Phi )(\xi )|,|\mathcal{F}(\Psi
)(\xi )|\geq c,  \label{Ass2.1}
\end{equation}%
if $|\xi |\leq 5/3$ and

\begin{equation}
\text{supp}(\mathcal{F}(\varphi ))\cup \text{supp}(\mathcal{F}(\psi
))\subset \overline{B(0,2)}\backslash B(0,1/2),\text{\quad }|\mathcal{F}%
(\varphi )(\xi )|,|\mathcal{F}(\psi )(\xi )|\geq c,  \label{Ass3.1}
\end{equation}%
if $3/5\leq |\xi |\leq 5/3$, such that%
\begin{equation}
\overline{\mathcal{F}(\Phi )(\xi )}\mathcal{F}(\Psi )(\xi
)+\sum_{k=1}^{\infty }\overline{\mathcal{F(}\varphi )(2^{-k}\xi )}\mathcal{F}%
(\psi )(2^{-k}\xi )=1,\quad \xi \in \mathbb{R}^{n},  \label{Ass4.1}
\end{equation}%
where $c>0$. Let $\Phi ,\varphi \in \mathcal{S}(\mathbb{R}^{n})$\ satisfy,
respectively, $\mathrm{\eqref{Ass2.1}}$ and $\mathrm{\eqref{Ass3.1}}$. We
recall that by \cite[pp. 130--131]{FJ90} or \cite[Lemma 6.9]{FrJaWe01},
there exist functions $\Psi \in \mathcal{S}(\mathbb{R}^{n})$ satisfying $%
\mathrm{\eqref{Ass2.1}}$ and $\psi \in \mathcal{S}(\mathbb{R}^{n})$
satisfying $\mathrm{\eqref{Ass3.1}}$ such that $\mathrm{\eqref{Ass4.1}}$ holds.

The $\varphi $-transform $S_{\varphi }$ is defined by setting 
\begin{equation*}
(S_{\varphi}f)_{0,m}=\langle f,\Psi _{m}\rangle ,
\end{equation*}
with $\Psi _{m}(x)=\Psi (x-m)$ and 
\begin{equation*}
(S_{\varphi }f)_{k,m}=\langle f,\varphi _{k,m}\rangle, 
\end{equation*}
where $\varphi
_{k,m}(x)=2^{kn/2}\varphi (2^{k}x-m)$, $k\in \mathbb{N}$ and $m\in \mathbb{Z}%
^{n}$. The inverse $\varphi $-transform $T_{\psi }$ is defined by 
\begin{equation*}
T_{\psi }\lambda =\sum_{m\in \mathbb{Z}^{n}}\lambda _{0,m}\Psi
_{m}+\sum_{k=1}^{\infty }\sum_{m\in \mathbb{Z}^{n}}\lambda _{k,m}\psi _{k,m},
\end{equation*}%
with $\lambda =\{\lambda _{k,m}\}_{k\in \mathbb{N}_{0},m\in \mathbb{Z}%
	^{n}}\subset \mathbb{C}$, see again \cite{FJ90}.

Now we present the inhomogenous version of Definition \ref{Tyulenev-class}.

\begin{defn}
	\label{Tyulenev-class-inho}Let $\alpha _{1}$, $\alpha _{2}\in \mathbb{R}$, $%
	\sigma _{1}$, $\sigma _{2}$ $\in (0,+\infty ]$, $\alpha =(\alpha _{1},\alpha
	_{2})$\ and let $\sigma =(\sigma _{1},\sigma _{2})$. We let $X_{\alpha
		,\sigma ,p}=X_{\alpha ,\sigma ,p}(\mathbb{R}^{n})$ denote the set of $p$%
	-admissible weight sequences $\{t_{k}\}_{k\in \mathbb{N}_{0}}$ satisfying $%
	\mathrm{\eqref{Asum1}}$ and $\mathrm{\eqref{Asum2}}$ for any $0\leq k\leq j$%
	, with constants $C_{1},C_{2}>0$ are independent of both the indexes $k$ and 
	$j$.
\end{defn}

\begin{ex}
	A sequence $\{\gamma _{j}\}_{j\in \mathbb{N}_{0}}$ of positive real numbers
	is said to be admissible if there exist two positive constants $d_{0}$ and $%
	d_{1}$ such that%
	\begin{equation}
	d_{0}\gamma _{j}\leq \gamma _{j+1}\leq d_{1}\gamma _{j},\quad j\in \mathbb{N}%
	_{0}.  \label{Farkas-Leop}
	\end{equation}%
	For an admissible sequence $\{\gamma _{j}\}_{j\in \mathbb{N}_{0}}$, let%
	\begin{equation*}
	\underline{\gamma }_{j}=\inf_{k\geq 0}\frac{\gamma _{j+k}}{\gamma _{k}}\quad 
	\text{and}\quad \overline{\gamma }_{j}=\sup_{k\geq 0}\frac{\gamma _{j+k}}{%
		\gamma _{k}},\quad j\in \mathbb{N}_{0}.
	\end{equation*}%
	Let%
	\begin{equation*}
	\alpha _{\gamma }=\lim_{j\longrightarrow \infty }\frac{\log \overline{\gamma 
		}_{j}}{j}\quad \text{and}\quad \beta _{\gamma }=\lim_{j\longrightarrow
		\infty }\frac{\log \underline{\gamma }_{j}}{j},
	\end{equation*}%
	be the upper and lower Boyd index of the given sequence $\{\gamma
	_{j}\}_{j\in \mathbb{N}_{0}}$, respectively. Then%
	\begin{equation*}
	\underline{\gamma }_{j}\gamma _{k}\leq \gamma _{j+k}\leq \overline{\gamma }%
	_{j}\gamma _{k},\quad j,k\in \mathbb{N}_{0}
	\end{equation*}%
	and for each $\varepsilon >0$,%
	\begin{equation*}
	c_{1}2^{(\beta _{\gamma }-\varepsilon )j}\leq \underline{\gamma }_{j}\leq 
	\overline{\gamma }_{j}\leq c_{2}2^{(\alpha _{\gamma }+\varepsilon )j},\quad
	j\in \mathbb{N}_{0}
	\end{equation*}%
	for some constants $c_{1}=c_{1}(\varepsilon )>0$ and $c_{2}=c_{2}(%
	\varepsilon )>0$. Also, $\underline{\gamma }_{1}$ and $\overline{\gamma }%
	_{1} $ are the best possible constants $d_{0}$ and $d_{1}$ in %
	\eqref{Farkas-Leop}, respectively. \newline
	Clearly the sequence $\{\gamma _{j}\}_{j\in \mathbb{N}_{0}}$ lies in $%
	X_{\alpha ,\sigma ,p}$ for$\ \alpha _{1}=\beta _{\gamma }-\varepsilon
	,\alpha _{2}=\alpha _{\gamma }+\varepsilon $ and $0<p,\sigma _{1},\sigma
	_{2}\leq \infty $.\newline
	These type of admissible sequences are used in \cite{FL06} to study Besov and Triebel-Lizorkin spaces in terms of a generalized smoothness, see also 	\cite{HaS08}.\newline
	Let us consider some examples of admissible sequences. The sequence $%
	\{\gamma _{j}\}_{j\in \mathbb{N}_{0}}$,%
	\begin{equation*}
	\gamma _{j}=2^{sj}(1+j)^{b}(1+\log (1+j))^{c},\quad j\in \mathbb{N}_{0}
	\end{equation*}%
	with arbitrary fixed real numbers $s,b$ and $c$ is a an admissible sequence
	with 
	\begin{equation*}
	\beta _{\gamma }=\alpha _{\gamma }=s.
	\end{equation*}
\end{ex}

\begin{ex}
	\label{Example1 copy(1)}Let $0<r<p<\infty $, a weight $\omega ^{p}\in A_{%
		\frac{p}{r}}(\mathbb{R}^{n})$ and
	\begin{equation*}
	\{s_{k}\}=\{2^{ks}\omega ^{p}(2^{-k})\}_{k\in \mathbb{N}_{0}},\quad s\in \mathbb{R}.
	\end{equation*}
	Obviously, $\{s_{k}\}_{k\in \mathbb{N}_{0}}$ lies in $X_{\alpha ,\sigma ,p}$ for $\alpha
	_{1}=\alpha _{2}=s$, $\sigma =(r(p/r)^{\prime },p)$.\ 
\end{ex}

Now, we define the spaces under consideration.

\begin{defn}
	\label{B-F-def-inh}Let $0<p< \infty $ and $0<q\leq \infty $. Let $%
	\{t_{k}\}_{k\in \mathbb{N}_{0}}$ be a $p$-admissible weight sequence. Let $%
	\Phi ,\varphi \in \mathcal{S}(\mathbb{R}^{n})$\ satisfy $\mathrm{%
		\eqref{Ass2.1}}$ and $\mathrm{\eqref{Ass3.1}}$, respectively, and we put $%
	\varphi _{k}=2^{kn}\varphi(2^{k}\cdot) ,k\in \mathbb{N}_{0}$. The Triebel-Lizorkin space $F_{p,q}(%
	\mathbb{R}^{n},\{t_{k}\}_{k\in \mathbb{N}_{0}})$\ is the collection of all $%
	f\in \mathcal{S}^{\prime }(\mathbb{R}^{n})$\ such that 
	\begin{equation*}
	\big\|f|F_{p,q}(\mathbb{R}^{n},\{t_{k}\}_{k\in \mathbb{N}_{0}})\big\|=\Big\|%
	\Big(\sum\limits_{k=0}^{\infty }t_{k}^{q}|\varphi _{k}\ast f|^{q}\Big)%
	^{1/q}|L_{p}(\mathbb{R}^{n})\Big\|<\infty,
	\end{equation*}%
	with the usual modifications if $q=\infty $, where $\varphi _{0}$ is
	replaced by $\Phi $.
\end{defn}

Now we introduce the inhomogeneous sequence spaces $f_{p,q}(\mathbb{R}%
^{n},\{t_{k}\}_{k\in \mathbb{N}_{0}})$. Let $0<p< \infty $ and $0<q\leq
\infty $. Let $\{t_{k}\}_{k\in \mathbb{N}_{0}}$ be a $p$-admissible weight
sequence. Then for all complex valued sequences $\lambda =\{\lambda
_{k,m}\}_{k\in \mathbb{N}_{0},m\in \mathbb{Z}^{n}}\subset \mathbb{C}$ we
define%
\begin{equation*}
f_{p,q}(\mathbb{R}^{n},\{t_{k}\}_{k\in \mathbb{N}_{0}})=\Big\{\lambda :\big\|%
\lambda |f_{p,q}(\mathbb{R}^{n},\{t_{k}\}_{k\in \mathbb{N}_{0}})\big\|%
<\infty \Big\},
\end{equation*}%
where%
\begin{equation*}
\big\|\lambda |f_{p,q}(\mathbb{R}^{n},\{t_{k}\}_{k\in \mathbb{N}_{0}})\big\|=%
\Big\|\Big(\sum_{k=0}^{\infty }2^{knq/2}\sum\limits_{m\in \mathbb{Z}%
	^{n}}t_{k}^{q}|\lambda _{k,m}|^{q}\chi _{k,m}\Big)^{1/q}|L_{p}(\mathbb{R}%
^{n})\big\|.
\end{equation*}%
We have the following analogue of Theorem \ref{phi-tran}.

\begin{thm}
	\label{phi-tran-inho}Let $\alpha =(\alpha _{1},\alpha _{2})\in \mathbb{R}%
	^{2},0<\theta \leq p<\infty $ and$\ 0<q< \infty $. Let $\{t_{k}\}_{k\in 
		\mathbb{N}_{0}}\in X_{\alpha ,\sigma ,p}$ be a $p$-admissible weight
	sequence with $\sigma =(\sigma _{1}=\theta \left( p/\theta \right) ^{\prime
	},\sigma _{2}\geq p)$.\ Let $\varphi $, $\psi $ satisfying $\mathrm{%
		\eqref{Ass1.1}}$\ through\ $\mathrm{\eqref{Ass4.1}}$. The operators 
	\begin{equation*}
	S_{\varphi }:F_{p,q}(\mathbb{R}^{n},\{t_{k}\}_{k\in \mathbb{N}%
		_{0}})\rightarrow f_{p,q}(\mathbb{R}^{n},\{t_{k}\}_{k\in \mathbb{N}_{0}})
	\end{equation*}%
	and 
	\begin{equation*}
	T_{\psi }:f_{p,q}(\mathbb{R}^{n},\{t_{k}\}_{k\in \mathbb{N}_{0}})\rightarrow
	F_{p,q}(\mathbb{R}^{n},\{t_{k}\}_{k\in \mathbb{N}_{0}})
	\end{equation*}%
	are bounded. Furthermore, $T_{\psi }\circ S_{\varphi }$ is the identity on $%
	F_{p,q}(\mathbb{R}^{n},\{t_{k}\}_{k\in \mathbb{N}_{0}})$.
\end{thm}

As a consequence the analogues of Corollary \ref{Indpendent} are now clear.
We obtain the following useful properties of these function spaces.

\begin{thm}
	Let $\alpha =(\alpha _{1},\alpha _{2})\in \mathbb{R}^{2},0<\theta \leq
	p<\infty $ and $0<q< \infty $. Let $\{t_{k}\}_{k\in \mathbb{N}_{0}}\in
	X_{\alpha ,\sigma ,p}$ be a $p$-admissible weight sequence with $\sigma
	=(\sigma _{1}=\theta \left( p/\theta \right) ^{\prime },\sigma _{2}\geq p)$. 
	$F_{p,q}(\mathbb{R}^{n},\{t_{k}\}_{k\in \mathbb{N}_{0}})$ are quasi-Banach
	spaces. They are Banach spaces if $1\leq p<\infty $ and $1\leq q<\infty $.
\end{thm}

Let $0<\theta \leq p<\infty $ and $0<q< \infty $. Let\ $\{t_{k}\}_{k\in 
	\mathbb{N}_{0}}\in X_{\alpha ,\sigma ,p}$ be a $p$-admissible weight
sequence with $\sigma =(\sigma _{1}=\theta \left( p/\theta \right) ^{\prime
},\sigma _{2}\geq p)$ and $\alpha =(\alpha _{1},\alpha _{2})\in \mathbb{R}%
^{2}$.\ As in Theorem \ref{embeddings-S-inf} we have the embedding%
\begin{equation*}
\mathcal{S}(\mathbb{R}^{n})\hookrightarrow F_{p,q}(\mathbb{R}%
^{n},\{t_{k}\}_{k\in \mathbb{N}_{0}}).
\end{equation*}%
In addition $\mathcal{S}(\mathbb{R}^{n})$ is dense in $F_{p,q}(\mathbb{R}%
^{n},\{t_{k}\}_{k\in \mathbb{N}_{0}})$. Also if $0<\theta \leq p<\infty $ and $0<q< \infty $, then 
\begin{equation*}
F_{p,q}(\mathbb{R}^{n},\{t_{k}\}_{k\in \mathbb{N}_{0}})\hookrightarrow 
\mathcal{S}^{\prime }(\mathbb{R}^{n}).
\end{equation*}%
All the results in Subsection 3.2 are true for the inhomogeneous case.

We begin with the following elementary embeddings, where the proof can be
obtained by using the properties of sequence Lebesgue spaces.

\begin{thm}
	\label{elem-embedding copy(1)}Let $0<p<\infty $ and $0<q\leq r<\infty $. Let 
	$\{t_{k}\}_{k\in \mathbb{N}_{0}}\in X_{\alpha ,\sigma ,p}$ be a $p$%
	-admissible weight sequence with $\sigma =(\sigma _{1}=\theta \left(
	p/\theta \right) ^{\prime },\sigma _{2}\geq p)$. We have 
	\begin{equation*}
	F_{p,q}(\mathbb{R}^{n},\{t_{k}\}_{k\in \mathbb{N}_{0}})\hookrightarrow
	F_{p,r}(\mathbb{R}^{n},\{t_{k}\}_{k\in \mathbb{N}_{0}}).
	\end{equation*}
\end{thm}

As in Subection 3.2 we obtain the following Sobolev-type embedding. We set%
\begin{equation*}
w_{k,Q}(p_{1})=\Big(\int_{Q}w_{k}^{p_{1}}(x)dx\Big)^{1/p_{1}}\quad \text{%
	and\quad }t_{k,Q}(p_{0})=\Big(\int_{Q}t_{k}^{p_{0}}(x)dx\Big)^{1/p_{0}},%
\end{equation*}%
where $ Q\in \mathcal{Q}$ with $\ell (Q)=2^{-k},k\in \mathbb{N}_{0}.$

\begin{thm}
	\label{Sobolev-embedding-sequence copy(1)}Let $0<\theta \leq
	p_{0}<p_{1}<\infty $ and $0<q,r<\infty $. Let\ $\{t_{k}\}_{k\in \mathbb{N}%
		_{0}}$ be a $p_{0}$-admissible weight sequence satisfying $\mathrm{%
		\eqref{Asum1}}$\ with $p=p_{0}$, $\sigma _{1}=\theta \left( p_{0}/\theta
	\right) ^{\prime }$ and $j=k\geq 0$. Let\ $\{w_{k}\}$ be a $p_{1}$%
	-admissible weight sequence satisfying $\mathrm{\eqref{Asum1}}$\ with $%
	p=p_{1}$, $\sigma _{1}=\theta \left( p_{1}/\theta \right) ^{\prime }$ and $%
	j=k\geq 0$. If $w_{k,Q}(p_{1})\lesssim t_{k,Q}(p_{0})$ for all $Q\in 
	\mathcal{Q}$ with $\ell (Q)=2^{-k},k\in \mathbb{N}_{0}$, then we have%
	\begin{equation*}
	F_{p_{0},q}(\mathbb{R}^{n},\{t_{k}\}_{k\in \mathbb{N}_{0}})\hookrightarrow
	F_{p_{1},r}(\mathbb{R}^{n},\{w_{k}\}_{k\in \mathbb{N}_{0}}).
	\end{equation*}
\end{thm}

From Theorems {\ref{phi-tran-inho} and \ref{Sobolev-embedding-sequence
		copy(1)}}, we have the following Sobolev-type embedding conclusions for $%
F_{p,q}(\mathbb{R}^{n},\{t_{k}\})$.

\begin{thm}
	\label{Sobolev-embedding copy(1)}Let $0<\theta \leq p_{0}<p_{1}<\infty $ and 
	$0<q,r<\infty $. Let $\{t_{k}\}_{k\in \mathbb{N}_{0}}\in X_{\alpha
		_{0},\sigma ,p_{0}}$ be a $p_{0}$-admissible weight sequence with $\sigma
	=(\sigma _{1}=\theta \left( p_{0}/\theta \right) ^{\prime },\sigma _{2}\geq
	p_{0})$ and $\alpha _{0}=(\alpha _{1,0},\alpha _{2,0})\in \mathbb{R}^{2}$.
	Let\ $\{w_{k}\}_{k\in \mathbb{N}_{0}}\in X_{\alpha _{1},\sigma ,p_{1}}$ be a 
	$p_{1}$-admissible weight sequence with $\sigma =(\sigma _{1}=\theta \left(
	p_{1}/\theta \right) ^{\prime },\sigma _{2}\geq p_{1})$ and $\alpha
	_{1}=(\alpha _{1,1},\alpha _{2,1})\in \mathbb{R}^{2}$. Then%
	\begin{equation*}
	F_{p_{0},q}(\mathbb{R}^{n},\{t_{k}\}_{k\in \mathbb{N}_{0}})\hookrightarrow
	F_{p_{1},r}(\mathbb{R}^{n},\{w_{k}\}_{k\in \mathbb{N}_{0}}),
	\end{equation*}%
	hold if%
	\begin{equation*}
	w_{k,Q}(p_{1})\lesssim t_{k,Q}(p_{0})
	\end{equation*}%
	for all $Q\in \mathcal{Q}$ and all $k\in \mathbb{N}_{0}$ with $\ell
	(Q)=2^{-k},k\in \mathbb{N}_{0}$ .
\end{thm}
In the sequel, we shall say that an operator $A$ is associated with the matrix 
\begin{equation*}
 \{a_{Q_{k,m}P_{v,h}}\}_{k,v\in \mathbb{N}_{0},m,h\in \mathbb{Z}^{n}},
 \end{equation*}
 if for all sequences $\lambda =\{\lambda _{k,m}\}_{k\in \mathbb{N}_{0},m\in 
	\mathbb{Z}^{n}}\subset \mathbb{C}$,%
\begin{equation*}
A\lambda =\{(A\lambda )_{k,m}\}_{k\in \mathbb{N}_{0},m\in \mathbb{Z}^{n}}=%
\Big\{\sum_{v=0}^{\infty }\sum_{h\in \mathbb{Z}^{n}}a_{Q_{k,m}P_{v,h}}%
\lambda _{v,h}\Big\}_{k\in \mathbb{N}_{0},m\in \mathbb{Z}^{n}}.
\end{equation*}%
We say that $A$, with associated matrix  $\{a_{Q_{k,m}P_{v,h}}\}_{k,v\in 
	\mathbb{N}_{0},m,h\in \mathbb{Z}^{n}}$, is almost diagonal on $f_{p,q}(%
\mathbb{R}^{n},\{t_{k}\}_{k\in \mathbb{N}_{0}})$ if there exists $%
\varepsilon >0$ such that%
\begin{equation*}
\sup_{k,v\in \mathbb{N}_{0},m,h\in \mathbb{Z}^{n}}\frac{|a_{Q_{k,m}P_{v,h}}|%
}{\omega _{Q_{k,m}P_{v,h}}(\varepsilon )}<\infty ,
\end{equation*}%
where $\omega _{Q_{k,m}P_{v,h}}(\varepsilon )$ as in Section 5. Let $\alpha
_{1},\alpha _{2}\in \mathbb{R},0<\theta \leq p<\infty $ and $0<q\leq \infty $%
. Let $\{t_{k}\}_{k\in \mathbb{N}_{0}}\in X_{\alpha ,\sigma ,p}$ be a $p$%
-admissible weight sequence with $\sigma _{1}=\theta \left( p/\theta \right)
^{\prime }$ and\ $\sigma _{2}\geq p$. It is obvious that an operator $A$ on $%
f_{p,q}(\mathbb{R}^{n},\{t_{k}\}_{k\in \mathbb{N}_{0}})$ given by an almost
diagonal matrix is bounded.

Let $J$ be defined as in Section 5. we present the inhomogeneous versions of
Definition \ref{Atom-Def}.

\begin{defn}
	\label{Atom-Def-inho}Let\ $\alpha _{1},\alpha _{2}\in \mathbb{R},0<p<\infty $
	and $0<q\leq \infty $. Let $\{t_{k}\}_{k\in \mathbb{N}_{0}}$ be a $p$%
	-admissible weight sequence. Let $N=\max \{J-n-\alpha _{1},-1\}$ and $\alpha
	_{2}^{\ast }=\alpha _{2}-\lfloor \alpha _{2}\rfloor $.\newline
	$\mathrm{(i)}$\ We say that $\varrho _{Q_{k,m}}$, $k\in \mathbb{N}_{0},m\in 
	\mathbb{Z}^{n}$, is an inhomogeneous smooth synthesis molecule for $F_{p,q}(%
	\mathbb{R}^{n},\{t_{k}\}_{k\in \mathbb{N}_{0}})$\ supported near $Q_{k,m}$
	if it satisfies, for some real number $\delta \in (\alpha _{2}^{\ast },1]$
	and a real number $M\in (J,\infty )$, $\mathrm{\eqref{mom-cond}}$, $\mathrm{%
		\eqref{cond1}}$\textrm{, }$\mathrm{\eqref{cond2}}$ and $\mathrm{\eqref{cond3}%
	}$ if $k\in \mathbb{N}$. If $k=0$ we assume $\mathrm{\eqref{cond2}}$, $%
	\mathrm{\eqref{cond3}}$ and 
	\begin{equation*}
	|\varrho _{Q_{0,m}}(x)|\leq (1+|x-x_{Q_{0,m}}|)^{-M}.
	\end{equation*}%
	A collection $\{\varrho _{Q_{k,m}}\}_{k\in \mathbb{N}_{0},m\in \mathbb{Z}%
		^{n}}$ is called a family of inhomogeneous smooth synthesis molecules for $%
	F_{p,q}(\mathbb{R}^{n},\{t_{k}\}_{k\in \mathbb{N}_{0}})$, if each $\varrho
	_{Q_{k,m}}$ is an inhomogeneous smooth synthesis molecule for $F_{p,q}(%
	\mathbb{R}^{n},\{t_{k}\}_{k\in \mathbb{N}_{0}})$ supported near $Q_{k,m}$. 
	\newline
	$\mathrm{(ii)}$\ We say that $b_{Q_{k,m}}$, $k\in \mathbb{N}_{0},m\in 
	\mathbb{Z}^{n}$, is an inhomogeneous smooth analysis molecule for $F_{p,q}(%
	\mathbb{R}^{n},\{t_{k}\}_{k\in \mathbb{N}_{0}})$ supported near $Q_{k,m}$ if
	it satisfies, for some $\kappa \in ((J-\alpha _{2})^{\ast },1]$ and an $M\in
	(J,\infty )$, $\mathrm{\eqref{mom-cond2}}$, $\mathrm{\eqref{cond1.1}}$%
	\textrm{, }$\mathrm{\eqref{cond1.2}}$ and $\mathrm{\eqref{cond1.3}}$ if $%
	k\in \mathbb{N}$. If $k=0$ we assume $\mathrm{\eqref{cond1.2}}$, $\mathrm{%
		\eqref{cond1.3}}$ and 
	\begin{equation*}
	|b_{Q_{0,m}}(x)|\leq (1+|x-x_{Q_{0,m}}|)^{-M}.
	\end{equation*}%
	A collection $\{b_{Q_{k,m}}\}_{k\in \mathbb{N}_{0},m\in \mathbb{Z}^{n}}$ is called a family of inhomogeneous smooth analysis molecules for $F_{p,q}(%
	\mathbb{R}^{n},\{t_{k}\}_{k\in \mathbb{N}_{0}})$, if each $b_{Q_{k,m}}$ is an inhomogeneous smooth synthesis molecule for $F_{p,q}(\mathbb{R}%
	^{n},\{t_{k}\}_{k\in \mathbb{N}_{0}})$ supported near $Q_{k,m}$.
\end{defn}

As a consequence, we formulate  the inhomogeneous counterpart of Theorem \ref%
{molecules-dec}.

\begin{thm}
	\label{molecules-dec-inho}Let $\alpha _{1}$, $\alpha _{2}\in \mathbb{R}%
	,0<\theta \leq p<\infty $ and $0<q< \infty $. Let $\{t_{k}\}_{k\in 
		\mathbb{N}_{0}}\in X_{\alpha ,\sigma ,p}$ be a $p$-admissible weight
	sequence with $\sigma _{1}=\theta \left( p/\theta \right) ^{\prime }$ and\ $%
	\sigma _{2}\geq p$. Let $J,M,N,\delta $ and $\kappa $ be as in Definition {\ref{Atom-Def-inho}}. \newline
	$\mathrm{(i)}$\ If $f=\sum_{k=0}^{\infty }\sum_{m\in \mathbb{Z}^{n}}\varrho
	_{k,m}\lambda _{k,m}$, where $\{\varrho _{k,m}\}_{k\in \mathbb{N}_{0},m\in 
		\mathbb{Z}^{n}}$ is a family of inhomogeneous smooth synthesis molecules for 
	$F_{p,q}(\mathbb{R}^{n},\{t_{k}\}_{k\in \mathbb{N}_{0}})$, then for all $%
	\lambda \in f_{p,q}(\mathbb{R}^{n},\{t_{k}\}_{k\in \mathbb{N}_{0}})$ 
	\begin{equation*}
	{{\big\|}f|F_{p,q}(\mathbb{R}^{n},\{t_{k}\}_{k\in \mathbb{N}_{0}}){\big\|}%
		\lesssim {\big\|}\lambda |f_{p,q}(\mathbb{R}^{n},\{t_{k}\}_{k\in \mathbb{N}%
			_{0}}){\big\|}.}
	\end{equation*}%
	$\mathrm{(ii)}$\ Let $\{b_{k,m}\}_{k\in \mathbb{N}_{0},m\in \mathbb{Z}^{n}}$
	be a family of inhomogeneous\ smooth analysis molecules.\ Then for all\ $%
	f\in F_{p,q}(\mathbb{R}^{n},\{t_{k}\}_{k\in \mathbb{N}_{0}})$%
	\begin{equation*}
	{{\big\|}\{\langle f,b_{k,m}\rangle \}_{k\in \mathbb{N}_{0},m\in \mathbb{Z}%
			^{n}}|f_{p,q}(\mathbb{R}^{n},\{t_{k}\}_{k\in \mathbb{N}_{0}}){\big\|}%
		\lesssim {\big\|}f|F_{p,q}(\mathbb{R}^{n},\{t_{k}\}_{k\in \mathbb{N}_{0}})%
		\big\|.}
	\end{equation*}
\end{thm}

Now we present the analogue of smooth atomic decomposition. First we need 
the definition of inhomogeneous smooth.

\begin{defn}
	Let $\alpha _{1},\alpha _{2}\in \mathbb{R},0<p<\infty ,0<q\leq \infty $ and\ 
	$N=\max \{J-n-\alpha _{1},-1\}$. Let $\{t_{k}\}_{k\in \mathbb{N}_{0}}$ be a $%
	p$-admissible weight sequence. A function $a_{Q_{k,m}}$ is called an
	inhomogeneous smooth atom for $F_{p,q}(\mathbb{R}^{n},\{t_{k}\}_{k\in 
		\mathbb{N}_{0}})$ supported near $Q_{k,m}$, $k\in \mathbb{N}_{0}$ and $m\in 
	\mathbb{Z}^{n}$, if it is satisfies $\mathrm{\eqref{supp-cond}}$, $\mathrm{%
		\eqref{diff-cond}}$\ and $\mathrm{\eqref{mom-cond1}\ }$if $k\in \mathbb{N}$.
	If $k=0$ we assume $\mathrm{\eqref{supp-cond}}$ and $\mathrm{%
		\eqref{diff-cond}.}$
\end{defn}

A collection $\{a_{Q_{k,m}}\}_{k\in \mathbb{N}_{0},m\in \mathbb{Z}^{n}}$ is
called a family of inhomogeneous smooth atoms for $F_{p,q}(\mathbb{R}%
^{n},\{t_{k}\}_{k\in \mathbb{N}_{0}})$, if each $a_{Q_{k,m}}$ is an
inhomogeneous smooth atom for $F_{p,q}(\mathbb{R}^{n},\{t_{k}\}_{k\in 
	\mathbb{N}_{0}})$ supported near $Q_{k,m}$.

Now we come to the atomic decomposition theorem.

\begin{thm}
	\label{atomic-dec copy(1)}Let $\alpha _{1}$, $\alpha _{2}\in \mathbb{R}$, $%
	0<\theta \leq p<\infty $, $0<q< \infty $. Let $\{t_{k}\}_{k\in \mathbb{N}%
		_{0}}\in X_{\alpha ,\sigma ,p}$ be a $p$-admissible weight sequence with $%
	\sigma _{1}=\theta \left( p/\theta \right) ^{\prime }$ and\ $\sigma _{2}\geq
	p$. Then for each $f\in F_{p,q}(\mathbb{R}^{n},\{t_{k}\}_{k\in \mathbb{N}%
		_{0}})$, there exist a family\ $\{\varrho _{k,m}\}_{k\in \mathbb{N}_{0},m\in 
		\mathbb{Z}^{n}}$ of inhomogeneous smooth atoms for the spaces $F_{p,q}(\mathbb{R}%
	^{n},\{t_{k}\}_{k\in \mathbb{N}_{0}})$ and $\lambda =\{\lambda
	_{k,m}\}_{k\in \mathbb{N}_{0},m\in \mathbb{Z}^{n}}\in f_{p,q}(\mathbb{R}%
	^{n},\{t_{k}\}_{k\in \mathbb{N}_{0}})$ such that 
	\begin{equation*}
	f=\sum\limits_{k=0}^{\infty }\sum\limits_{m\in \mathbb{Z}^{n}}\lambda
	_{k,m}\varrho _{k,m},\text{\quad converging in }\mathcal{S}^{\prime }(%
	\mathbb{R}^{n})
	\end{equation*}%
	and%
	\begin{equation*}
	{{\big\|}\{\lambda _{k,m}\}_{k\in \mathbb{N}_{0},m\in \mathbb{Z}%
			^{n}}|f_{p,q}(\mathbb{R}^{n},\{t_{k}\}_{k\in \mathbb{N}_{0}}){\big\|}%
		\lesssim {\big\|}f|F_{p,q}(\mathbb{R}^{n},\{t_{k}\}_{k\in \mathbb{N}_{0}})%
		\big\|.}
	\end{equation*}%
	Conversely, for any family of inhomogeneous smooth atoms for the spaces $F_{p,q}(%
	\mathbb{R}^{n},\{t_{k}\}_{k\in \mathbb{N}_{0}})$ and $\lambda =\{\lambda
	_{k,m}\}_{k\in \mathbb{N}_{0},m\in \mathbb{Z}^{n}}\in f_{p,q}(\mathbb{R}%
	^{n},\{t_{k}\}_{k\in \mathbb{N}_{0}})$ we have
	\begin{align*}
	& {\big\|}\sum\limits_{k=0}^{\infty }\sum\limits_{m\in \mathbb{Z}^{n}}\lambda
	_{k,m}\varrho _{k,m}|F_{p,q}(\mathbb{R}^{n},\{t_{k}\}_{k\in \mathbb{N}_{0}})%
	\big\|\\
	& \lesssim {\big\|}\{\lambda _{k,m}\}_{k\in \mathbb{N}_{0},m\in \mathbb{Z%
		}^{n}}|f_{p,q}(\mathbb{R}^{n},\{t_{k}\}_{k\in \mathbb{N}_{0}}){\big\|}.
	\end{align*}
\end{thm}

\begin{rem}
	One of the applications of atomic and molecular decompositions of the spaces 
	$\dot{F}_{p,q}(\mathbb{R}^{n},\{t_{k}\})$ and $F_{p,q}(\mathbb{R}%
	^{n},\{t_{k}\}_{k\in \mathbb{N}_{0}})$ is studying the continuity of
	singular integral operators of non convolution type in such function spaces,
	where it is enough to show that it maps every family of smooth atoms into a
	family of smooth molecules. For classical Triebel-Lizorkin spaces we refer
	the reader to, e.g., {\cite{FJHW88} and \cite{Torres91}.}\\
	It is also interesting to study other properties and characterizations of these function spaces such as the wavelet characterization, see, e.g., Lemari\'{e} and Meyer {\cite{LM88}}, and Triebel  {\cite{T4}}.
\end{rem}

\subsection*{Acknowledgment}

The author would like to thank W. Sickel and A. Tyulenev for valuable
discussions and suggestions. \\
We thank the referees for carefully reading our paper and for their several useful
suggestions and comments, which improved the exposition of the paper substantially.\\
This work is found by the General Direction of Higher Education and Training under\ Grant No. C00L03UN280120220004 and by
The General Directorate of Scientific Research and Technological
Development, Algeria.

\end{document}